# RIGOROUS ANALYSIS FOR THE DIRAC SYSTEM ON THE QUARTER-PLANE

Hassan Babaei [1] , Jerry L. Bona [1] , Andreas Chatziafratis [2,*]

[1] Department of Mathematics, Statistics and Computer Science, University of Illinois at Chicago, USA

[2] Department of Mathematics, National and Kapodistrian University of Athens, Greece

June 15, 2026

**Abstract**. Considered and analyzed below are fully non-homogeneous initial-boundary-value problems for the celebrated Dirac system, formulated on the spatial half-line. Analytical solution formulae are derived formally via suitable implementation of the well-known Fokas' unified transform methodology, and rigorously verified a posteriori. The latter substantial task relies on complex-analytic tools and careful interpretation of the obtained integral representations. These valid solutions are then used for investigating qualitative properties. These include boundary behavior near the axes of the domain as well as long-range asymptotics and long-time (eventual) periodicity. Notably, smoothness of the solution, both within and upto the boundary of the domain, depends heavily on certain compatibility conditions between initial, boundary and forcing data. Further results pertaining to solution's regularity and uniqueness are thence established based on the qualitative theory. The closed-form expressions reported here are also useful in the study of non-linear counterparts.

## 1. Introduction

The Dirac equation was first introduced by Paul Dirac as a relativistic wave equation for the electron [1]. It is a (3+1)-dimensional equation that appears in quantum field theory and describes the dynamics of spin-1/2 particles; see, for example, [2-4]. In its free one-space-dimensional form, with normalized physical constants, it can be written as

$$i\partial_t \Psi = -i\alpha\partial_x \Psi + m\beta\Psi ,$$

where $m \geq 0$ is the mass of a free particle, $\Psi = \begin{pmatrix} \psi_1 \\ \psi_2 \end{pmatrix} \in \mathbb{C}^2$ is a two-component complex vector called spinor field, and $\alpha$, $\beta$ are $2\times 2$ Hermitian and traceless matrices, known as the Dirac matrices, satisfying

$$\alpha^2 = \beta^2 = I \text{ and } \alpha\beta + \beta\alpha = 0 .$$

There is, however, no unique preferred $2\times 2$ representation in one spatial dimension. Different choices of matrices $\alpha$ and $\beta$ lead to different-looking systems, but these systems are equivalent under unitary transformations of spinor variables. Thus the representation changes the form of the component equations, and hence the form of the boundary-value problem and of the associated transform formulas, but not the underlying free Dirac dynamics [3,4].

One convenient form is

$$i\partial_t \begin{pmatrix} \psi_1 \\ \psi_2 \end{pmatrix} = -i\begin{pmatrix} 0 & 1 \\ 1 & 0 \end{pmatrix}\partial_x \begin{pmatrix} \psi_1 \\ \psi_2 \end{pmatrix} + m\begin{pmatrix} 0 & 1 \\ 1 & -1 \end{pmatrix}\begin{pmatrix} \psi_1 \\ \psi_2 \end{pmatrix},$$

---

*Acknowledgement.* H. Babaei is partially supported by the UIC's Graduate College Award for Graduate Research.

* corresponding author, e-mail address: chatziafrati@math.uoa.gr

which is known as the diagonal mass representation. Simplifying it yields

$$\begin{cases} \partial_t \psi_1 = -\partial_x \psi_2 - im\psi_1 \\ \partial_t \psi_2 = -\partial_x \psi_1 + im\psi_2 . \end{cases}$$

This is the Dirac system we study in the present work. In contrast, we note that in diagonal transport variables, the two components propagate along the characteristic families $x - t = constant$ and $x + t = constant$, while the mass term couples the two components. On the half-line $x > 0$, only one characteristic family enters through the boundary, and hence one prescribes a single boundary trace rather than both component traces.

We herein study the initial-boundary-value problem for this Dirac system on the quarter-plane

$$Q = \{(x,t) : x > 0, t > 0\},$$

with initial conditions

$$\psi_1(x,0) = g_1(x),\ \psi_2(x,0) = g_2(x),\ x > 0,$$

and a single boundary condition

$$\psi_1(0,t) = h_1(t),\ x > 0,$$

which turns out to be a sufficiently determined, solvable problem.

The purpose of this paper is to rigorously solve and analyze this quarter-plane problem through the *unified transform method* (UTM). This modern PDE method was introduced by Fokas a few decades ago as a systematic approach to boundary-value problems for linear and integrable nonlinear evolution equations as well as for elliptic PDE [5-10]; see also, e.g., [11-15], in chronological order.

A previous application of the UTM to the one-dimensional Dirac equation derived formulas for the massive and massless problems on the half-line and on a finite interval [16]. That work is useful as an early UTM treatment of the Dirac equation, but several important analytical issues were not addressed and many questions still remain open. In particular, the contour formulas were obtained formally, the convergence of the representation was not justified term by term, while *a posteriori* verification of the formula, including recovery of the initial and boundary data, together with corner compatibility, were not considered. The present work provides a rigorous treatment of those issues and establishes many qualitative properties of the problems under consideration and the associated solutions as well.

A few works have also addressed nonlinear Dirac equations in half-line and quarter-plane settings, including cubic Dirac equations, Thirring-type models, and low-regularity problems with vector self-interaction [17-22]. These studies focus on questions such as well-posedness, global existence, asymptotic behavior, and, in integrable cases, inverse-scattering or Riemann–Hilbert formulations. However, the fundamental linear quarter-plane analysis remained largely unexplored, and is thus developed here.

The analysis undertaken here is close in spirit to an established research program, led by Chatziafratis with a large international network of collaborators, on rigorous studies of exact transform solutions and qualitative theory for linear evolution equations on semi-unbounded domains, including analyses on: Benjamin-Bona-Mahony, higher-order Cahn-Hilliard-Kuramoto-Sivashinsky, Korteweg-de Vries / Airy, hyperbolic Maxwell-Cattaneo-Vernotte (damped wave / telegrapher), Rubinstein-Aifantis (double-diffusivity), Schrödinger, Sobolev-Barenblatt (pseudo-parabolic), Whitham-Broer-Kaup and other (systems of) PDE [23-35], in alphabetical order. Foundational ideas towards this ever-expanding line of rigorous UTM-based investigations were set forth in the seminal work [36]. For very recent advances, the reader is referred to, for instance, [37-40].

In the present work, the UTM for our Dirac problem leads to spectral functions involving the "dispersion" quantities

$$\rho(\lambda) = \sqrt{\lambda^2 + m^2},\ \omega_1(\lambda) = i\rho(\lambda),\ \omega_2(\lambda) = -i\rho(\lambda),$$

with a suitably chosen branch structure in the complex λ-plane. This implies that, for example, certain oscillatory integrals of the resulting closed-form representations are not automatically absolutely convergent on their original contours, rendering their analytical handling delicate; their rigorous interpretation requires appropriate, non-trivial analyticity arguments, contour deformations, Riemann-Lebesgue and Jordan-type lemmas, and other tools from real, complex and asymptotic analysis.

One of the main contributions of this work is to give a complete *a posteriori* verification of the obtained UTM-type representation formulae for the Dirac quarter-plane problem, and then to utilize these rigorously justified expressions to study the solution's qualitative behavior. Specifically, we prove that the integral representations define solutions in the appropriate generalized sense, that the differential equations hold away from the characteristic line x = t, and that the prescribed initial and boundary data are recovered. We also analyze the boundary behavior near the axes of the quarter-plane as well as the diagonal. Importantly, compatibility conditions are brought forth which are critical for smoothness across the whole domain including the axes and corner. Moreover, far-field asymptotics and large-time behavior (eventual periodicity) are examined. Finally, the corresponding inhomogeneous (forced) problem, further regularity as well as uniqueness properties are treated within our framework. Results pertaining to nonlinear counterparts will be announced elsewhere.

We therefore here begin by addressing, for the first time, the following:

***Problem*** *Solve, rigorously, the system of partial differential equations*

$$\begin{cases}\partial_t\psi_1 = -\partial_x\psi_2 - im\psi_1 \\ \partial_t\psi_2 = -\partial_x\psi_1 + im\psi_2,\end{cases} \tag{1.1}$$

*for* $\psi_1(x,t)$ *and* $\psi_2(x,t)$, *with* $(x,t)\in\mathbb{R}^+\times\mathbb{R}^+$, *subject to the following initial and boundary conditions:*

$$\begin{cases}\lim_{t\to0^+}\psi_1(x,t) = g_1(x),\ \ x\in\mathbb{R}^+, \\ \lim_{t\to0^+}\psi_2(x,t) = g_2(x),\ \ x\in\mathbb{R}^+, \\ \lim_{x\to0^+}\psi_1(x,t) = h_1(t),\ \ t\in\mathbb{R}^+,\end{cases} \tag{1.2}$$

*and investigate the solution's qualitative properties.*

***Assumptions*** Throughtout this paper, we make the following assumptions:

$$g_1(x),\, g_2(x)\in\mathcal{S}([0,\infty)),\, h_1(t)\in C^\infty([0,\infty)). \tag{1.3}$$

***The UTM-based solution*** The solution we propose for the system $(1.1)\&(1.2)$ is the following:
For $x>0$, $t>0$ and $x\neq t$,

$$\begin{aligned}(1.4)\ \ 2\pi\psi_1(x,t) = &\int_{-\infty}^{\infty} e^{i\lambda x-\omega_1(\lambda)t}\{\lambda\hat{g}_2(\lambda)+[i\omega_2(\lambda)+m]\hat{g}_1(\lambda)\}\frac{d\lambda}{2\rho(\lambda)} \\ &+\int_{-\infty}^{\infty} e^{i\lambda x-\omega_1(\lambda)t}\{\lambda\hat{g}_2(-\lambda)-[i\omega_2(\lambda)+m]\hat{g}_1(-\lambda)\}\frac{d\lambda}{2\rho(\lambda)} \\ &-\int_{-\infty}^{\infty} e^{i\lambda x-\omega_2(\lambda)t}\{\lambda\hat{g}_2(\lambda)+[i\omega_1(\lambda)+m]\hat{g}_1(\lambda)\}\frac{d\lambda}{2\rho(\lambda)}\end{aligned}$$

$$-\int_{-\infty}^{\infty} e^{i\lambda x-\omega_2(\lambda)t}\{\lambda\hat{g}_2(-\lambda)-[i\omega_1(\lambda)+m]\hat{g}_1(-\lambda)\}\frac{d\lambda}{2\rho(\lambda)}$$

$$+\int_{-\infty}^{\infty} e^{i\lambda x-\omega_1(\lambda)t}\left[\int_{\tau=0}^{t} h_1(\tau)e^{\omega_1(\lambda)\tau}d\tau\right]\frac{\lambda d\lambda}{\rho(\lambda)}-\int_{-\infty}^{\infty} e^{i\lambda x-\omega_2(\lambda)t}\left[\int_{\tau=0}^{t} h_1(\tau)e^{\omega_2(\lambda)\tau}d\tau\right]\frac{\lambda d\lambda}{\rho(\lambda)}$$

and

$$(1.5)\quad 2\pi\psi_2(x,t)=\int_{-\infty}^{\infty} e^{i\lambda x-\omega_1(\lambda)t}\{\lambda\hat{g}_1(\lambda)-[i\omega_1(\lambda)+m]\hat{g}_2(\lambda)\}\frac{d\lambda}{2\rho(\lambda)}$$

$$-\int_{-\infty}^{\infty} e^{i\lambda x-\omega_1(\lambda)t}\{\lambda\hat{g}_1(-\lambda)+[i\omega_1(\lambda)+m]\hat{g}_2(-\lambda)\}\frac{d\lambda}{2\rho(\lambda)}$$

$$-\int_{-\infty}^{\infty} e^{i\lambda x-\omega_2(\lambda)t}\{\lambda\hat{g}_1(\lambda)-[i\omega_2(\lambda)+m]\hat{g}_2(\lambda)\}\frac{d\lambda}{2\rho(\lambda)}$$

$$+\int_{-\infty}^{\infty} e^{i\lambda x-\omega_2(\lambda)t}\{\lambda\hat{g}_1(-\lambda)+[i\omega_2(\lambda)+m]\hat{g}_2(-\lambda)\}\frac{d\lambda}{2\rho(\lambda)}$$

$$-\int_{-\infty}^{\infty} e^{i\lambda x-\omega_1(\lambda)t}\left[\int_{\tau=0}^{t} h_1(\tau)e^{\omega_1(\lambda)\tau}d\tau\right]\frac{[i\omega_1(\lambda)+m]d\lambda}{\rho(\lambda)}$$

$$+\int_{-\infty}^{\infty} e^{i\lambda x-\omega_2(\lambda)t}\left[\int_{\tau=0}^{t} h_1(\tau)e^{\omega_2(\lambda)\tau}d\tau\right]\frac{[i\omega_2(\lambda)+m]d\lambda}{\rho(\lambda)}.$$

The above solution is derived by the unified transform method. We will see that the solution extends also for $x=t$, under the additional assumptions: $h_1(0)=g_1(0)$ and $h_1'(0)=-g_2'(0)$. More precisely, we will show that the so extended functions $\psi_1(x,t)$ and $\psi_2(x,t)$ are $C^1$ for $(x,t)\in\mathbb{R}^+\times\mathbb{R}^+$.

***Notation* (1)** $\hat{g}_j(\lambda)=\int_0^{\infty} e^{-i\lambda y}g_j(y)dy$, $j=1,2$, defined for $\lambda\in\mathbb{C}$, with $\operatorname{Im}\lambda\le 0$.

**(2)** $\omega_1(\lambda)=i\rho(\lambda)=i\sqrt{\lambda^2+m^2}$ and $\omega_2(\lambda)=-i\rho(\lambda)=-i\sqrt{\lambda^2+m^2}$.

Here we use the following notation concerning the square-root function $\rho(\lambda)=\sqrt{\lambda^2+m^2}$ : We set

$$\Theta:=\mathbb{C}-\{\lambda=\xi+i\eta\in\mathbb{C}:\xi=0\ \&\ |\eta|\ge m\}=\mathbb{C}-\{(-i\infty,-mi]\cup[mi,i\infty)\},$$

and define $\rho:\Theta\to\mathbb{C}$ by

$$\rho(\lambda)=\sqrt{|\lambda-mi|\cdot|\lambda+mi|}\exp\{i[\theta_{mi}(\lambda)+\theta_{-mi}(\lambda)-\pi]/2\},\ \text{for}\ \lambda\in\Theta\ \ (\text{see fig. 1}).$$

($\sqrt{|\lambda-mi|\cdot|\lambda+mi|}$ denotes the positive square root of the positive number $|\lambda-mi|\cdot|\lambda+mi|$, for $\lambda\ne\pm mi$.)

Let us note that

$$\lambda-mi=|\lambda-mi|e^{i[\theta_{mi}(\lambda)-\frac{\pi}{2}]}\ \text{and}\ \lambda+mi=|\lambda+mi|e^{i[\theta_{-mi}(\lambda)-\frac{\pi}{2}]},$$

so that

$$\lambda^2 + m^2 = (\lambda - mi)(\lambda + mi) = |\lambda - mi||\lambda + mi| e^{i[\theta_{mi}(\lambda) + \theta_{-mi}(\lambda) - \pi]}, \text{ i.e., } [\rho(\lambda)]^2 = \lambda^2 + m^2 .$$

With this choice of the angles $\theta_{mi}(\lambda)$ and $\theta_{-mi}(\lambda)$, the function $\rho(\lambda)$ is analytic for $\lambda \in \Theta$ and $\rho(-\lambda) = \rho(\lambda)$. Also if $\lambda \in \mathbb{R}$ then $\rho(\lambda)$ is equal to the positive square root of the positive number $\lambda^2 + m^2$.

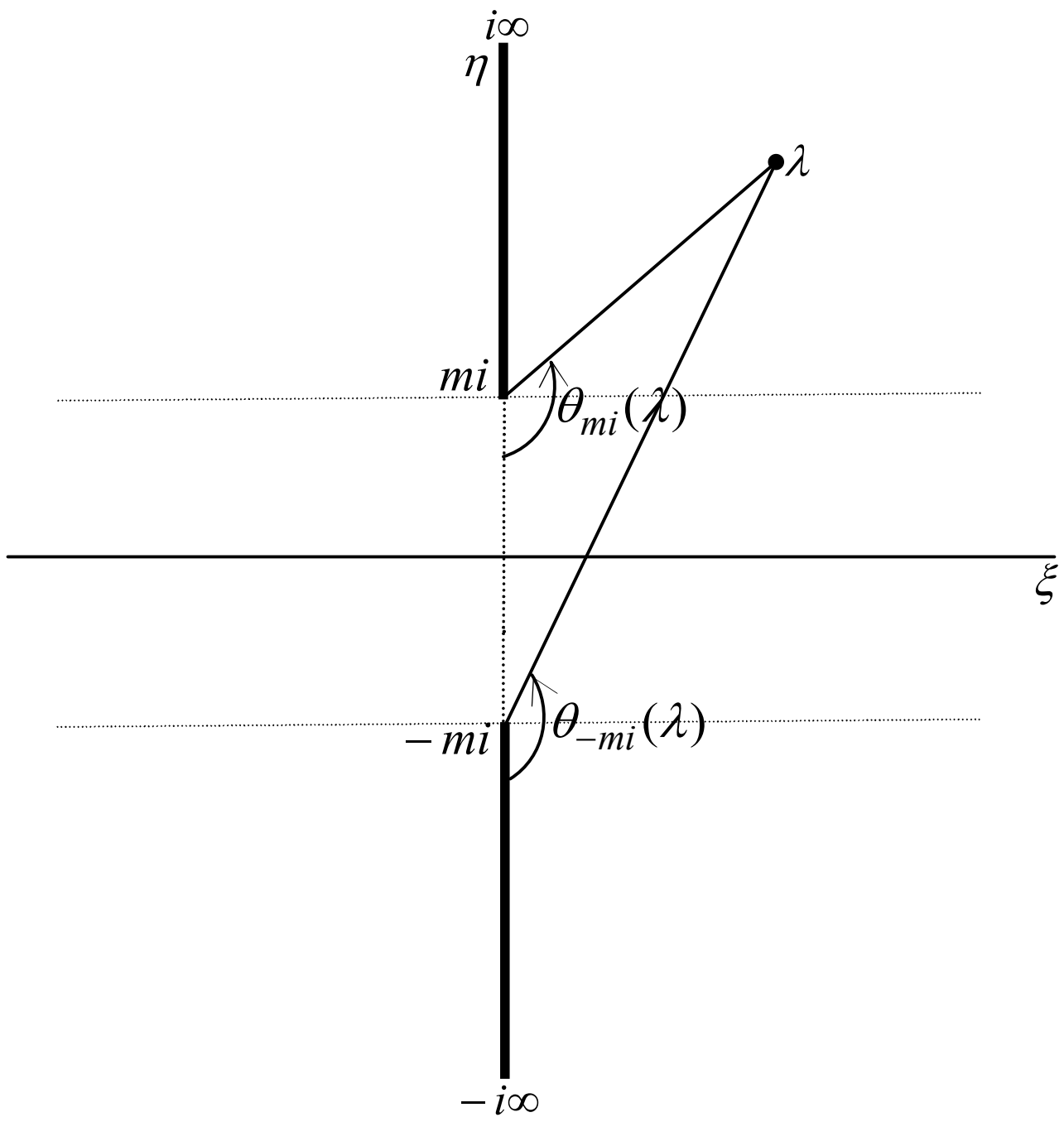

**Fig. 1** The choice of the angles $\theta_{mi}(\lambda)$ and $\theta_{-mi}(\lambda): -\pi < \theta_{mi}(\lambda) < \pi$ and $0 < \theta_{-mi}(\lambda) < 2\pi$.

**(3)** We will denote by

$$\mathcal{G}_{1,1}(x,t),\ \mathcal{G}_{1,2}(x,t),\mathcal{G}_{1,3}(x,t),\ \mathcal{G}_{1,4}(x,t),\ \mathcal{H}_{1,1}(x,t),\ \mathcal{H}_{1,2}(x,t),$$

the integrals which define $\psi_1(x,t)$, in the order they appear in the RHS of (1.4). With this notation, (1.4) is written: $2\pi\psi_1 = \mathcal{G}_{1,1} + \mathcal{G}_{1,2} - \mathcal{G}_{1,3} - \mathcal{G}_{1,4} + \mathcal{H}_{1,1} - \mathcal{H}_{1,2}$.

Similarly, we will denote by

$$\mathcal{G}_{2,1}(x,t),\ \mathcal{G}_{2,2}(x,t),\mathcal{G}_{2,3}(x,t),\ \mathcal{G}_{2,4}(x,t),\ \mathcal{H}_{2,1}(x,t),\ \mathcal{H}_{2,2}(x,t),$$

the integrals which define $\psi_2(x,t)$, in the order they appear in the RHS of (1.5). With this notation, (1.5) is written: $2\pi\psi_2 = \mathcal{G}_{2,1} - \mathcal{G}_{2,2} - \mathcal{G}_{2,3} + \mathcal{G}_{2,4} - \mathcal{H}_{2,1} + \mathcal{H}_{2,2}$.

Finally, for an integral of the form $I = \int_{\lambda=-\infty}^{\infty} \cdots d\lambda$, we will denote by $I^+$ and $I^-$, the following parts of $I$:

$$I^+ = \int_0^\infty \cdots d\lambda \text{ and } I^- = \int_{-\infty}^0 \cdots d\lambda .$$

For example,

$$\mathcal{G}_{1,2}^-(x,t) = \int_{-\infty}^0 e^{i\lambda x - \omega_1(\lambda)t} \{\lambda \hat{g}_2(-\lambda) - [i\omega_2(\lambda) + m]\hat{g}_1(-\lambda)\} \frac{d\lambda}{2\rho(\lambda)} .$$

We will see that all the integrals

$$\mathcal{G}^{+}_{j,l}(x,t)\,,\ \mathcal{G}^{-}_{j,l}(x,t)\,,\ \mathcal{H}^{+}_{j,s}(x,t)\,,\ \mathcal{H}^{-}_{j,s}(x,t)\,,\ j\in\{1,2\}\,,\ l\in\{1,2,3,4\}\,,\ s\in\{1,2\}\,,$$

which are parts of the solutions $\psi_1$ and $\psi_2$, exist in the generalized sense, i.e., as limits of the form

$$\lim_{A\to\infty}\int_{-A}^{0}\ ,\ \lim_{A\to\infty}\int_{0}^{A}\ ,$$

provided that $x\neq t$.

**Theorem 1** *Let $Q:=\mathbb{R}^{+}\times\mathbb{R}^{+}$. Under the assumptions (1.3), the following hold:*

*1st The integrals in (1.4) and (1.5) exist in the generalized sence, provided that $(x,t)\in Q-\{x=t\}$.*

*2nd The functions $\psi_1(x,t)$ and $\psi_2(x,t)$, defined by (1.4) and (1.5), are $C^{\infty}$ for $(x,t)$ in the sets $Q-\{x<t\}$ and $Q-\{x>t\}$.*

*3rd $\psi_1(x,t)$ and $\psi_2(x,t)$ satisfy the differential equations (1.1) in $Q-\{x=t\}$.*

*4th $\lim_{t\to 0^{+}}\psi_1(x,t)=g_1(x)$ and $\lim_{t\to 0^{+}}\psi_2(x,t)=g_2(x)$, uniformly for $x$ in compact sets of $\mathbb{R}^{+}$.*

*5th $\lim_{x\to 0^{+}}\psi_1(x,t)=h_1(t)$, uniformly for $t$ in compact sets of $\mathbb{R}^{+}$.*

*6th If, in addition, $h_1(0)=g_1(0)$ then the functions $\psi_1(x,t)$ and $\psi_2(x,t)$ extend continuously for $(x,t)\in Q$.*

*7th If, in addition, $h_1(0)=g_1(0)$ and $h_1'(0)=-g_2'(0)$ then the extentions of the functions $\psi_1(x,t)$ and $\psi_2(x,t)$ to $Q$, are continuously differentiable ($C^1$) and satisfy the differential equations (1.1) in $Q$.*

## 2. Derivation of the solution

We start with the *assumption* that there exist functions $\psi_1(x,t)$ and $\psi_2(x,t)$, solving the problem (1.1)&(1.2), and having properties which ***justify*** the steps of the following construction. (After we derive formulas for $\psi_1(x,t)$ and $\psi_2(x,t)$, we will *rigorously* verify, *a posteriori*, that, indeed, this pair of functions, $(\psi_1(x,t),\ \psi_2(x,t))$, is the solution of problem (1.1)&(1.2), provided that, in addition to (1.3), the data satisfy the compatibity conditions at the origin: $h_1(0)=g_1(0)$ and $h_1'(0)=-g_2'(0)$. This will be done in section 3, in the proof of Theorem 1.)

For $\lambda\in\mathbb{C}$ with $\operatorname{Im}\lambda\leq 0$, we define the half-line Fourier transform

$$\text{(2.1)}\qquad \hat{\psi}_j(\lambda,t)=\int_0^{\infty}e^{-i\lambda y}\psi_j(y,t)dy\,,\ \partial_t\hat{\psi}_j(\lambda,t)=\int_0^{\infty}e^{-i\lambda y}\partial_t\psi_j(y,t)dy\,,\ j=1,2\,.$$

Then, the differential equations (1.1) give the following:

$$\begin{aligned}\text{(2.2)}\qquad \partial_t\hat{\psi}_1(\lambda,t)&=-\int_0^{\infty}e^{-i\lambda y}\partial_y\psi_2(y,t)dy-im\hat{\psi}_1(\lambda,t)\\ &=-[e^{-i\lambda y}\psi_2(y,t)]\Big|_{y=0}^{y=\infty}-i\lambda\int_0^{\infty}e^{-i\lambda y}\psi_2(y,t)dy-im\hat{\psi}_1(\lambda,t)\\ &=h_2(t)-i\lambda\hat{\psi}_2(\lambda,t)-im\hat{\psi}_1(\lambda,t)\end{aligned}$$

(where we set $h_2(t) := \psi_2(0,t)$ ) and

(2.3) $$\partial_t \hat{\psi}_2(\lambda,t) = -\int_0^\infty e^{-i\lambda y} \partial_y \psi_1(y,t)dy + im\hat{\psi}_2(\lambda,t)$$

$$= -[e^{-i\lambda y}\psi_1(y,t)]\Big|_{y=0}^{y=\infty} - i\lambda \int_0^\infty e^{-i\lambda y}\psi_2(y,t)dy + im\hat{\psi}_2(\lambda,t)$$

$$= h_1(t) - i\lambda\hat{\psi}_1(\lambda,t) + im\hat{\psi}_2(\lambda,t) .$$

Thus, we are led to solve the system of the ODEs

(2.4) $$\begin{cases} \partial_t \hat{\psi}_1(\lambda,t) = -i\lambda\hat{\psi}_2(\lambda,t) - im\hat{\psi}_1(\lambda,t) + h_2(t) \\ \partial_t \hat{\psi}_2(\lambda,t) = -i\lambda\hat{\psi}_1(\lambda,t) + im\hat{\psi}_2(\lambda,t) + h_1(t) \end{cases}$$

with initial conditions

(2.5) $$\hat{\psi}_1(\lambda,t)\Big|_{t=0} = \hat{g}_1(\lambda) \text{ and } \hat{\psi}_2(\lambda,t)\Big|_{t=0} = \hat{g}_2(\lambda),$$

in the unknown functions $\hat{\psi}_1(\lambda,t)$ and $\hat{\psi}_2(\lambda,t)$, of the variable $t$, with $\lambda$ being a parameter, restricted to $\lambda \in \mathbb{C}$ with $\mathrm{Im}\,\lambda \le 0$.

Solving the homogeneous system of ODEs, corresponding to (1.6), namely

$$\begin{cases} \partial_t u_1(\lambda,t) = -i\lambda u_2(\lambda,t) - imu_1(\lambda,t) \\ \partial_t u_2(\lambda,t) = -i\lambda u_1(\lambda,t) + imu_2(\lambda,t), \end{cases}$$

we find (for a fixed $\lambda$ ) the general solution

$$u_1(\lambda,t) = A\lambda e^{-\omega_1(\lambda)t} + B\lambda e^{-\omega_2(\lambda)t}, \quad u_2(\lambda,t) = -A[i\omega_1(\lambda)+m]e^{-\omega_1(\lambda)t} - B[i\omega_2(\lambda)+m]e^{-\omega_2(\lambda)t},$$

with $A, B \in \mathbb{C}$.

Now, letting $A$ and $B$ depend on $t$, we see that if the functions

(2.6) $$\hat{\psi}_1(\lambda,t) = A_0(t)\lambda e^{-\omega_1(\lambda)t} + B_0(t)\lambda e^{-\omega_2(\lambda)t}$$

and

(2.7) $$\hat{\psi}_2(\lambda,t) = -A_0(t)[i\omega_1(\lambda)+m]e^{-\omega_1(\lambda)t} - B_0(t)[i\omega_2(\lambda)+m]e^{-\omega_2(\lambda)t}$$

solve the system (2.4) then

$$A_0'(t)\lambda e^{-\omega_1(\lambda)t} + B_0'(t)\lambda e^{-\omega_2(\lambda)t} = h_2(t)$$

and

$$A_0'(t)[i\omega_1(\lambda)+m]e^{-\omega_1(\lambda)t} + B_0'(t)[i\omega_2(\lambda)+m]e^{-\omega_2(\lambda)t} = -h_1(t) .$$

Solving the above equations for $A_0'(t)$ and $B_0'(t)$, we find

$$A_0'(t) = \frac{\lambda h_1(t) + [i\omega_2(\lambda)+m]h_2(t)}{2\lambda\rho(\lambda)} e^{\omega_1(\lambda)t} \text{ and } B_0'(t) = -\frac{\lambda h_1(t) + [i\omega_1(\lambda)+m]h_2(t)}{2\lambda\rho(\lambda)} e^{\omega_2(\lambda)t} .$$

Thus, with the choice

$$A_0(t)=\int_{\tau=0}^{t}\frac{\lambda h_1(\tau)+[i\omega_2(\lambda)+m]h_2(\tau)}{2\lambda\rho(\lambda)}e^{\omega_1(\lambda)\tau}d\tau$$

and

$$B_0(t)=-\int_{\tau=0}^{t}\frac{\lambda h_1(\tau)+[i\omega_1(\lambda)+m]h_2(\tau)}{2\lambda\rho(\lambda)}e^{\omega_2(\lambda)\tau}d\tau\,,$$

the functions (2.6)&(2.7) is a particular solution of the system (2.4).
It follows that the general solution of (2.4) is given by the formulas

$$\text{(2.8)}\qquad \begin{cases}\hat{\psi}_1(\lambda,t)=[A+A_0(t)]\lambda e^{-\omega_1(\lambda)t}+[B+B_0(t)]\lambda e^{-\omega_2(\lambda)t}\\ \hat{\psi}_2(\lambda,t)=-[A+A_0(t)][i\omega_1(\lambda)+m]e^{-\omega_1(\lambda)t}-[B+B_0(t)][i\omega_2(\lambda)+m]e^{-\omega_2(\lambda)t},\end{cases}$$

with $A,B\in\mathbb{C}$.
In order to determine the solution of (2.4), which satisfies also the initial conditions (2.5), it suffices to set $t=0$ in (2.8) and solve the equations

$$A\lambda+B\lambda=\hat{g}_1(\lambda)\ \ \&\ \ A[i\omega_1(\lambda)+m]+B[i\omega_2(\lambda)+m]=-\hat{g}_2(\lambda),$$

for $A,B$. This gives

$$A=\frac{\lambda\hat{g}_2(\lambda)+[i\omega_2(\lambda)+m]\hat{g}_1(\lambda)}{2\lambda\rho(\lambda)}\ \text{ and }\ B=-\frac{\lambda\hat{g}_2(\lambda)+[i\omega_1(\lambda)+m]\hat{g}_1(\lambda)}{2\lambda\rho(\lambda)}\,.$$

Substituting in (2.8), we find

$$\text{(2.9)}\ \hat{\psi}_1(\lambda,t)=\left\{\frac{\lambda\hat{g}_2(\lambda)+[i\omega_2(\lambda)+m]\hat{g}_1(\lambda)}{2\rho(\lambda)}+\left[\int_{\tau=0}^{t}\frac{\lambda h_1(\tau)+[i\omega_2(\lambda)+m]h_2(\tau)}{2\rho(\lambda)}e^{\omega_1(\lambda)\tau}d\tau\right]\right\}e^{-\omega_1(\lambda)t}$$
$$-\left\{\frac{\lambda\hat{g}_2(\lambda)+[i\omega_1(\lambda)+m]\hat{g}_1(\lambda)}{2\rho(\lambda)}+\left[\int_{\tau=0}^{t}\frac{\lambda h_1(\tau)+[i\omega_1(\lambda)+m]h_2(\tau)}{2\rho(\lambda)}e^{\omega_2(\lambda)\tau}d\tau\right]\right\}e^{-\omega_2(\lambda)t}$$

and

$$\text{(2.10)}\ \hat{\psi}_2(\lambda,t)=\left\{\frac{\lambda\hat{g}_1(\lambda)-[i\omega_1(\lambda)+m]\hat{g}_2(\lambda)}{2\rho(\lambda)}+\left[\int_{\tau=0}^{t}\frac{\lambda h_2(\tau)-[i\omega_1(\lambda)+m]h_1(\tau)}{2\rho(\lambda)}e^{\omega_1(\lambda)\tau}d\tau\right]\right\}e^{-\omega_1(\lambda)t}$$
$$-\left\{\frac{\lambda\hat{g}_1(\lambda)-[i\omega_2(\lambda)+m]\hat{g}_2(\lambda)}{2\rho(\lambda)}+\left[\int_{\tau=0}^{t}\frac{\lambda h_2(\tau)-[i\omega_2(\lambda)+m]h_1(\tau)}{2\rho(\lambda)}e^{\omega_2(\lambda)\tau}d\tau\right]\right\}e^{-\omega_2(\lambda)t}.$$

These are the solutions of (2.4) & (2.5). Let us also keep in mind that the domain of analyticity of $\rho(\lambda)$ is the set $\Theta=\mathbb{C}-\{(-i\infty,-mi]\cup[mi,i\infty)\}$. Thus, the validity of the equations (2.9) and (2.10, is restriced, for $\lambda$, to $\Theta\cap\{\operatorname{Im}\lambda\le 0\}$. The equations (2.9) & (2.10) are the global relations (GRs), associated to problem (1.1)&(1.2).

Next we proceed to eliminate, from the GRs, the unknown quantity $h_2(\tau)$, which is not part of the data. In order to achieve this, we set $-\lambda$ in (2.9), in place of $\lambda$, and we subtract the resulting equation from (2.9), obtaining, for $\lambda\in\mathbb{R}$,

$$\text{(2.11)}\ \hat{\psi}_1(\lambda,t)-\hat{\psi}_1(-\lambda,t)$$
$$=\left\{\frac{\lambda\hat{g}_2(\lambda)+[i\omega_2(\lambda)+m]\hat{g}_1(\lambda)}{2\rho(\lambda)}-\frac{-\lambda\hat{g}_2(-\lambda)+[i\omega_2(\lambda)+m]\hat{g}_1(-\lambda)}{2\rho(\lambda)}\right\}e^{-\omega_1(\lambda)t}$$

$$-\left\{\frac{\lambda\hat{g}_2(\lambda)+[i\omega_1(\lambda)+m]\hat{g}_1(\lambda)}{2\rho(\lambda)}-\frac{-\lambda\hat{g}_2(-\lambda)+[i\omega_1(\lambda)+m]\hat{g}_1(-\lambda)}{2\rho(\lambda)}\right\}e^{-\omega_2(\lambda)t}$$

$$+\frac{\lambda}{\rho(\lambda)}\left[\int_{\tau=0}^{t}h_1(\tau)e^{\omega_1(\lambda)\tau}d\tau\right]e^{-\omega_1(\lambda)t}-\frac{\lambda}{\rho(\lambda)}\left[\int_{\tau=0}^{t}h_1(\tau)e^{\omega_2(\lambda)\tau}d\tau\right]e^{-\omega_2(\lambda)t}.$$

Similarly, setting $-\lambda$ in (2.10), in place of $\lambda$, and adding the resulting equation to (2.10), we find, for $\lambda\in\mathbb{R}$,

(2.12) $\hat{\psi}_2(\lambda,t)+\hat{\psi}_2(-\lambda,t)$

$$=\left\{\frac{\lambda\hat{g}_1(\lambda)-[i\omega_1(\lambda)+m]\hat{g}_2(\lambda)}{2\rho(\lambda)}+\frac{-\lambda\hat{g}_1(-\lambda)-[i\omega_2(\lambda)+m]\hat{g}_2(-\lambda)}{2\rho(\lambda)}\right\}e^{-\omega_1(\lambda)t}$$

$$-\left\{\frac{\lambda\hat{g}_1(\lambda)-[i\omega_2(\lambda)+m]\hat{g}_2(\lambda)}{2\rho(\lambda)}+\frac{-\lambda\hat{g}_1(-\lambda)-[i\omega_2(\lambda)+m]\hat{g}_2(-\lambda)}{2\rho(\lambda)}\right\}e^{-\omega_2(\lambda)t}$$

$$-\frac{i\omega_1(\lambda)+m}{\rho(\lambda)}\left[\int_{\tau=0}^{t}h_1(\tau)e^{\omega_1(\lambda)\tau}d\tau\right]e^{-\omega_1(\lambda)t}+\frac{i\omega_2(\lambda)+m}{\rho(\lambda)}\left[\int_{\tau=0}^{t}h_1(\tau)e^{\omega_2(\lambda)\tau}d\tau\right]e^{-\omega_2(\lambda)t}.$$

Since, by Fourier's inversion formula,

$$\int_{-\infty}^{\infty}e^{i\lambda x}\hat{\psi}_j(\lambda,t)d\lambda=2\pi\psi_j(x,t),$$

($j=1,2$) and, by Cauchy's integral formula and Jordan's lemma,

$$\int_{-\infty}^{\infty}e^{i\lambda x}\hat{\psi}_j(-\lambda,t)d\lambda=0,$$

equations (2.11) and (2.12) yield the solution formulas (1.4) and (1.5).

## 3. The integrals of the solution

***Preliminary comments***

The integrals which are used to express the solution with the formulas (1.15) and (1.16), in general, do not converge absolutely. Indeed, on the one hand, when $\lambda\in\mathbb{R}$ then $\left|e^{i\lambda x}\right|=\left|e^{i\lambda x-\omega_j(\lambda)t}\right|=1$ ($j=1,2$) and, on the other hand,

$$\hat{g}_j(\lambda)=\mathrm{O}(1/\lambda)\ \text{and}\ e^{-\omega_j(\lambda)t}\left[\int_0^t h_1(\tau)e^{\omega_j(\lambda)\tau}d\tau\right]=\mathrm{O}(1/\lambda),\ \text{as}\ \mathbb{R}\ni\lambda\to\pm\infty.$$

Moreover, we point out that, as far as their convergence is concerned, these integrals become worse (*more oscillatory*), if we differentiate them with respect to $x$ or $t$. For example, switching the order of differentiation and integration in the following case

$$\frac{\partial^{n+k}}{\partial x^n\partial t^k}\left[\int_{-\infty}^{\infty}e^{i\lambda x-\omega_1(\lambda)t}\{\lambda\hat{g}_1(\lambda)-[i\omega_1(\lambda)+m]\hat{g}_2(\lambda)\}\frac{d\lambda}{2\rho(\lambda)}\right]$$

would lead us to the integral

$$\int_{-\infty}^{\infty} (i\lambda)^n \left[-i\sqrt{\lambda^2+m^2}\right]^k e^{i\lambda x-\omega_1(\lambda)t}\{\lambda\hat{g}_1(\lambda)-[i\omega_1(\lambda)+m]\hat{g}_2(\lambda)\}\frac{d\lambda}{2\rho(\lambda)},$$

which is highly oscillatory.

In order to address this problem, we will deform appropriately the contours in the above integrals – more precisely in certain parts of these integrals – so that Lebesgue's dominated convergence theorem is applicable. The various deformation processes will be based on Cauchy's integral formula and Jordan's lemma.

Also, there arise additional difficulties in dealing with the functions of $(x,t)$ defined by the above integrals, if we let $x\to 0^+$ or $t\to 0^+$, i.e., switching the order of these limiting processes and integration. Addressing these difficulties is needed in the process of verification of the initial and boundary conditions (1.2) and – more generally – the study of the boundary behaviour of the derivatives of the solutions $\psi_1(x,t)$ and $\psi_2(x,t)$.

***Preliminary remarks*** We point out that the study and interpretation of the various integrals rely on the elementary fact that

$$\left|e^{i\lambda x}\right| = e^{-x(\operatorname{Im}\lambda)} \le e^{-\alpha x} \text{ for } \lambda\in\mathbb{C} \text{ with } \operatorname{Im}\lambda \ge \alpha > 0, \text{ and for } x>0.$$

We will use also this in the following forms:

$$\left|e^{i\lambda x-i\lambda t}\right| = e^{-(x-t)(\operatorname{Im}\lambda)} \le e^{-\alpha(x-t)} \text{ for } \lambda\in\mathbb{C} \text{ with } \operatorname{Im}\lambda \ge \alpha > 0, \text{ and for } x>t,$$

and

$$\left|e^{i\lambda x-i\lambda t}\right| = e^{-(x-t)(\operatorname{Im}\lambda)} \le e^{-\alpha(t-x)} \text{ for } \lambda\in\mathbb{C} \text{ with } \operatorname{Im}\lambda \le -\alpha < 0, \text{ and for } x<t.$$

Thus, several integrals in the representation formulas (1.4) and (1.5) (more precisely, certain parts of these integrals), which do not converge absolutely, will be found to be equal to absolutely convergent integrals. This is done by appropriately deforming their contours, and the main tools, to carry out these deformation processes, are Cauchy's integral theorem and Jordan's lemma.

For example,

$$\int_0^{\infty} e^{i\lambda x-i\lambda t}\varphi(x,t,\lambda)d\lambda = \int_{L^{+,+}} e^{i\lambda x-i\lambda t}\varphi(x,t,\lambda)d\lambda, \text{ for } x>t,$$

(for the contour $L^{+,+}$, see fig. 6, below) if the function $\varphi(x,t,\lambda)$ is continuous and, in addition, analytic in $\{\lambda\in\mathbb{C}:\operatorname{Im}\lambda>0\}$ and $\varphi(x,t,\lambda)=\mathrm{O}(1/\lambda)$, as $\lambda\to\infty$ (with $\lambda\in\mathbb{C}$, $\operatorname{Im}\lambda>0$).

Moreover,

$$\frac{\partial^{n,l}}{\partial x^n\partial t^l}\left[\int_0^{\infty} e^{i\lambda x-i\lambda t}\varphi(x,t,\lambda)d\lambda\right] = \int_{L^{+,+}} \frac{\partial^{n,l}}{\partial x^n\partial t^l}[e^{i\lambda x-i\lambda t}\varphi(x,t,\lambda)]d\lambda, \text{ for } x>t,$$

provided that, in addition, the function $\varphi(x,t,\lambda)$ is "appropriately" smooth.

Similarly,

$$\frac{\partial^{n,l}}{\partial x^n\partial t^l}\left[\int_{-\infty}^{0} e^{i\lambda x-i\lambda t}\varphi(x,t,\lambda)d\lambda\right] = \int_{L^{-,-}} \frac{\partial^{n,l}}{\partial x^n\partial t^l}[e^{i\lambda x-i\lambda t}\varphi(x,t,\lambda)]d\lambda, \text{ for } x<t,$$

with the function $\varphi(x,t,\lambda)$ being "appropriately" smooth and, in addition, analytic in the set $\{\lambda\in\mathbb{C}:\operatorname{Im}\lambda<0\}$ and $\varphi(x,t,\lambda)=\mathrm{O}(1/\lambda)$, as $\lambda\to\infty$ (with $\lambda\in\mathbb{C}$, $\operatorname{Im}\lambda<0$). (For the contour $L^{-,-}$, see fig.7, below.)

***Preliminary computations***

**(1)** With $N \in \mathbb{N}$, integrating by parts, we obtain

$$\hat{g}(\lambda) = \sum_{k=1}^{N} \frac{g^{(k-1)}(0)}{(i\lambda)^k} + \upsilon_{g,N}(\lambda), \text{ for } \lambda \in \mathbb{C} - \{0\} \text{ with } \operatorname{Im}\lambda \le 0, \tag{3.1}$$

where

$$\upsilon_{g,N}(\lambda) = \frac{1}{(i\lambda)^N} \int_0^\infty e^{-i\lambda y} \frac{d^N[g(y)]}{dt^N} dy = \frac{1}{(i\lambda)^N} (g^{(N)})\hat{}(\lambda).$$

Also, we have

$$\hat{g}(\lambda) = \mathrm{O}(1/\lambda) \text{ and } \upsilon_{g,N}(\lambda) = \mathrm{O}(1/\lambda^{N+1}), \text{ as } \lambda \to \infty, \text{ with } \operatorname{Im}\lambda \le 0. \tag{3.2}$$

**(2)** With $N \in \mathbb{N}$ and $\omega = \omega(\lambda) \in \{\omega_1(\lambda), \omega_2(\lambda)\}$, integrating by parts, we obtain

$$e^{-\omega t} \int_0^t e^{\omega\tau} h(\tau) d\tau = \sum_{k=1}^{N} (-1)^{k-1} \frac{h^{(k-1)}(t)}{\omega^k} - e^{-\omega t} \sum_{k=1}^{N} (-1)^{k-1} \frac{h^{(k-1)}(0)}{\omega^k} - U_{h,N}(\lambda,t) \tag{3.3}$$

where

$$U_{h,N}(\lambda,t) = (-1)^{N-1} \frac{1}{\omega^N} e^{-\omega t} \int_0^t e^{\omega\tau} h^{(N)}(\tau) d\tau.$$

Also, uniformly for $t$ in compact sets,

$$e^{-\omega(\lambda)t} \int_0^t e^{\omega(\lambda)\tau} h(\tau) d\tau = \mathrm{O}(1/\lambda) \text{ and } U_{h,N}(\lambda,t) = \mathrm{O}(1/\lambda^{N+1}), \text{ as } \lambda \to \infty, \text{ with } \operatorname{Re}\omega(\lambda) \ge 0. \tag{3.4}$$

Keeping the above observations in mind, we procced with the intepretation of the integrals in the formulas (1.4) and (1.5).

***Interpretation of the integrals of the solution***

**(1)** *The integral* $\mathcal{G}_{1,1}$ Let us write it as follows:

$$\begin{aligned} \mathcal{G}_{1,1}(x,t) &= \int_{-\infty}^{\infty} e^{i\lambda x - \omega_1(\lambda)t} \{\lambda \hat{g}_2(\lambda) + [i\omega_2(\lambda) + m]\hat{g}_1(\lambda)\} \frac{d\lambda}{2\rho(\lambda)} \\ &= \int_{\lambda=0}^{1} + \int_{\lambda=1}^{\infty} + \int_{\lambda=-1}^{0} + \int_{\lambda=-\infty}^{-1} = \mathcal{G}_{1,1}^{0,+}(x,t) + \mathcal{G}_{1,1}^{1,+}(x,t) + \mathcal{G}_{1,1}^{0,-}(x,t) + \mathcal{G}_{1,1}^{1,-}(x,t). \end{aligned} \tag{3.5}$$

Then

$$\begin{aligned} \mathcal{G}_{1,1}^{1,+}(x,t) &= \int_1^\infty e^{i\lambda x - i\rho(\lambda)t} \{\lambda \hat{g}_2(\lambda) + [\rho(\lambda) + m]\hat{g}_1(\lambda)\} \frac{d\lambda}{2\rho(\lambda)} \\ &= \int_1^\infty e^{i\lambda x - i\lambda t} \sigma_1(\lambda,t) \{\lambda \hat{g}_2(\lambda) + [\rho(\lambda) + m]\hat{g}_1(\lambda)\} \frac{d\lambda}{2\rho(\lambda)}, \end{aligned}$$

where

$$\sigma_1(\lambda,t) := \exp[i\lambda t - i\rho(\lambda)t] = \exp\left\{\frac{-im^2 t}{\lambda+\rho(\lambda)}\right\}, \text{ defined for } \lambda\in\mathbb{C} \text{ with } \operatorname{Re}\lambda>0,$$

and

$$\mathcal{G}_{1,1}^{1,-}(x,t) = \int_{-\infty}^{-1} e^{i\lambda x - i\rho(\lambda)t}\{\lambda\hat{g}_2(\lambda)+[\rho(\lambda)+m]\hat{g}_1(\lambda)\}\frac{d\lambda}{2\rho(\lambda)}$$
$$= \int_{-\infty}^{-1} e^{i\lambda x + i\lambda t}\sigma_2(\lambda,t)\{\lambda\hat{g}_2(\lambda)+[\rho(\lambda)+m]\hat{g}_1(\lambda)\}\frac{d\lambda}{2\rho(\lambda)}$$

where

$$\sigma_2(\lambda,t) := \exp[-i\lambda t - i\rho(\lambda)t] = \exp\left\{\frac{im^2 t}{\lambda-\rho(\lambda)}\right\}, \text{ defined for } \lambda\in\mathbb{C} \text{ with } \operatorname{Re}\lambda<0.$$

As far as the integrals

$$\mathcal{G}_{1,1}^{0,+}(x,t) = \int_0^1 \text{ and } \mathcal{G}_{1,1}^{0,-}(x,t) = \int_{-1}^0,$$

defined by (3.5), are concerned, they are easy to deal with, since their contours are finite.

*Note* It is crucial that

$$\lambda+\rho(\lambda) \cong 2|\lambda|e^{i\arg\lambda}, \text{ for } \lambda\to\infty \text{ with } \lambda\in\mathbb{C},\ \operatorname{Re}\lambda>0,$$

and

$$\lambda-\rho(\lambda) \cong 2|\lambda|e^{i\arg\lambda}, \text{ for } \lambda\to\infty \text{ with } \lambda\in\mathbb{C},\ \operatorname{Re}\lambda<0.$$

It follows that

$$\lim_{\substack{|\lambda|\to+\infty\\ \operatorname{Re}\lambda>0}} |\lambda+\rho(\lambda)| = +\infty \text{ and } \lim_{\substack{|\lambda|\to+\infty\\ \operatorname{Re}\lambda<0}} |\lambda-\rho(\lambda)| = +\infty,$$

and, as a consequence,

$$\lim_{\substack{|\lambda|\to+\infty\\ \operatorname{Re}\lambda>0}} \sigma_1(\lambda,t) = \lim_{\substack{|\lambda|\to+\infty\\ \operatorname{Re}\lambda>0}} \exp\left\{\frac{-im^2t}{\lambda+\rho(\lambda)}\right\} = 1 \text{ and } \lim_{\substack{|\lambda|\to+\infty\\ \operatorname{Re}\lambda<0}} \sigma_2(\lambda,t) = \lim_{\substack{|\lambda|\to+\infty\\ \operatorname{Re}\lambda<0}} \exp\left\{\frac{im^2t}{\lambda-\rho(\lambda)}\right\} = 1.$$

In particular, the quantities $\sigma_1(\lambda,t)$ and $\sigma_2(\lambda,t)$ are bounded for $\operatorname{Re}\lambda>0$ and $\operatorname{Re}\lambda<0$, respectively:

$$\sup_{\operatorname{Re}\lambda>0} |\sigma_1(\lambda,t)| < +\infty \text{ and } \sup_{\operatorname{Re}\lambda<0} |\sigma_2(\lambda,t)| < +\infty. \tag{3.6}$$

With (3.6) in mind, we continue with the study of the integral $\mathcal{G}_{1,1}^{+}(x,t)$ and we split it as follows:

$$\mathcal{G}_{1,1}^{1,+}(x,t) = J_1(x,t) + J_2(x,t),$$

where

$$J_1(x,t) = \int_1^\infty e^{i\lambda x - i\lambda t}\sigma_1(\lambda,t)\hat{g}_2(\lambda)\frac{\lambda d\lambda}{2\rho(\lambda)} \text{ and } J_2(x,t) = \int_1^\infty e^{i\lambda x - i\lambda t}\sigma_1(\lambda,t)\hat{g}_1(\lambda)\}\frac{[\rho(\lambda)+m]d\lambda}{2\rho(\lambda)}.$$

We distinguish between the cases « $x>t$ » and « $x<t$ ». (The case « $x=t$ » will be considered later.)

*The case* $x>t$.

In view of (3.1), Cauchy's theorem and Jordan's lemma, we have

(3.7) $$J_1(x,t)=\sum_{k=1}^{N} g_2^{(k)}(0)\int_{\Gamma^{+,+}} e^{i\lambda x-i\lambda t}\sigma_1(\lambda,t)\frac{1}{(i\lambda)^k}\frac{\lambda d\lambda}{2\rho(\lambda)} + \int_1^{\infty} e^{i\lambda x-i\lambda t}\sigma_1(\lambda,t)\frac{(g_2^{(N)})\hat{}(\lambda)}{(i\lambda)^N}\frac{\lambda d\lambda}{2\rho(\lambda)}, \text{ for } x>t,$$

where $\Gamma^{+,+}=[1,1+i]+\{\lambda=\xi+i\eta:\eta=\xi \text{ and } \xi\geq 1\}$, i.e., $\Gamma^{+,+}$ consists of the line segment $[1,1+i]$ and the half-line which starts at the point $1+i$ and its slope is $\pi/4$. (See fig. 2.)

*Note* In writing (3.7) we chose $\Gamma^{+,+}$, as described above. However, there are many other choices. For example one can choose any half-line in the first quarter starting at $1+i$, with positive slope, including the half-line $[1,1+\infty i)=\{\lambda=\xi+i\eta:\xi=1 \text{ and } \eta\geq 0\}$.

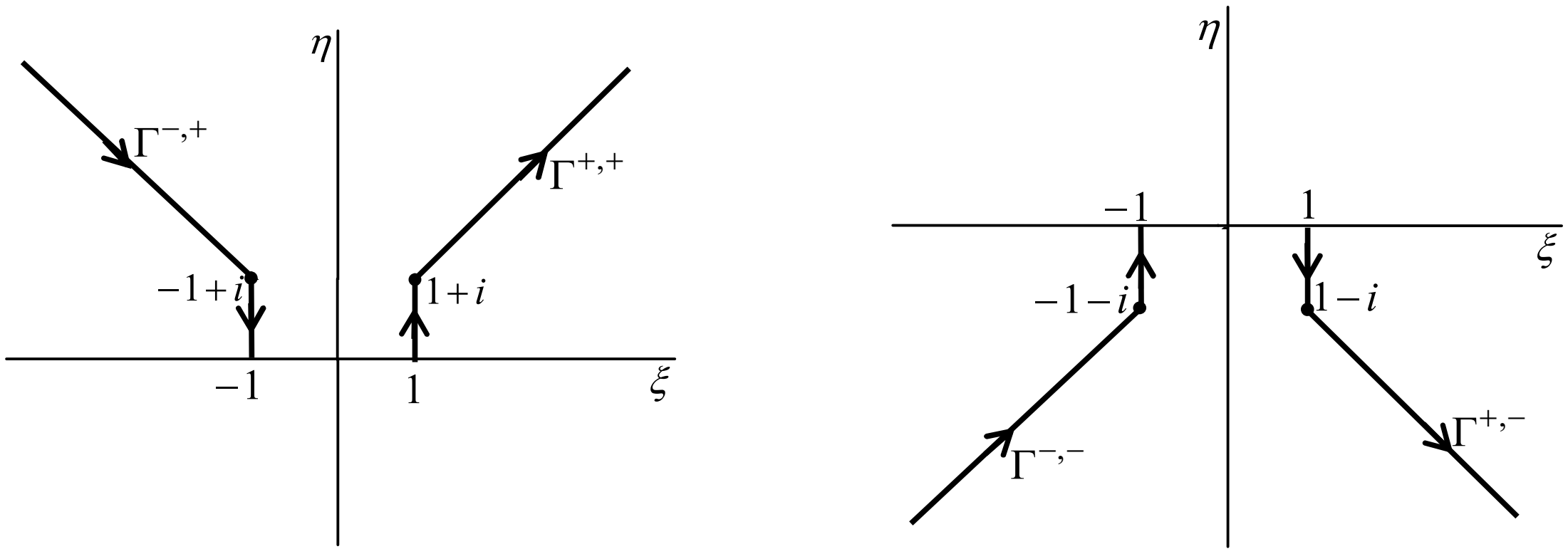

**Fig.2** *The contours* $\Gamma^{-,+}$ *and* $\Gamma^{+,+}$ **Fig.3** *The contours* $\Gamma^{-,-}$ *and* $\Gamma^{+,-}$

Similarly,

(3.8) $$J_2(x,t)=\sum_{k=1}^{N} g_1^{(k)}(0)\int_{\Gamma^{+,+}} e^{i\lambda x-i\lambda t}\sigma_1(\lambda,t)\frac{1}{(i\lambda)^k}\frac{[\rho(\lambda)+m]d\lambda}{2\rho(\lambda)} + \int_1^{\infty} e^{i\lambda x-i\lambda t}\sigma_1(\lambda,t)\frac{(g_1^{(N)})\hat{}(\lambda)}{(i\lambda)^N}\frac{[\rho(\lambda)+m]d\lambda}{2\rho(\lambda)}, \text{ for } x>t.$$

*The case* $x<t$.

We have

(3.9) $$\mathcal{G}_{1,1}^{1,+}(x,t)=\int_{\Gamma^{+,-}} e^{i\lambda x-i\lambda t}\sigma_1(\lambda,t)\{\lambda\hat{g}_2(\lambda)+[\rho(\lambda)+m]\hat{g}_1(\lambda)\}\frac{d\lambda}{2\rho(\lambda)}, \text{ for } x<t,$$

where $\Gamma^{+,-}=[1,1-i]+\{\lambda=\xi+i\eta:\eta=-\xi \text{ and } \xi\geq 1\}$. (See fig. 3.)

*In any case («* $x>t$ *», «* $x<t$ *» or «* $x=t$ *»).*

We have

(3.10) $$\mathcal{G}_{1,1}^{1,-}(x,t)=\int_{\Gamma^{-,+}} e^{i\lambda x+i\lambda t}\sigma_2(\lambda,t)\{\lambda\hat{g}_2(\lambda)+[\rho(\lambda)+m]\hat{g}_1(\lambda)\}\frac{d\lambda}{2\rho(\lambda)} \quad (x>0\,, t>0),$$

where $\Gamma^{-,+}=\{\lambda=\xi+i\eta:\eta=-\xi \text{ and } \xi\leq -1\}+[-1+i,-1]$. (See fig. 2.)

*In conclusion*: In view of (3.5), (3.6), (3.7), (3.8), (3.9) and (3.10), the integral $\mathcal{G}_{1,1}(x,t)$ can be written as a finite sum $\sum_{\ell}\mathcal{G}_{1,1,\ell}(x,t)$, of absolutely convergent integrals $\mathcal{G}_{1,1,\ell}(x,t)$, provided that $x\neq t$. In particular, from this analysis of $\mathcal{G}_{1,1}(x,t)$, it follows that the integrals $\mathcal{G}_{1,1}^{+}(x,t)$ and $\mathcal{G}_{1,1}^{-}(x,t)$ exist in the generalized sense, for $x \neq t$.

Furthermore, if we differentiate the integrands of $\mathcal{G}_{1,1,\ell}(x,t)$, we respect to $x$, $t$, any number of times, the resulting integrals remain absolutely and uniformly convergent, for $(x,t)$ in compact subsets of $Q-\{x=t\}$. Therefore, the function defined by the integral $\mathcal{G}_{1,1}(x,t)$ is $C^{\infty}$ for $(x,t) \in Q-\{x=t\}$.

**(2)** *The integrals* $\mathcal{G}_{1,2}, \mathcal{G}_{1,3}, \mathcal{G}_{1,4}, \mathcal{G}_{2,1}, \mathcal{G}_{2,2}, \mathcal{G}_{2,3}, \mathcal{G}_{2,4}$. The interpretation of these integrals is quite analogous to the interpretation of the integral $\mathcal{G}_{1,1}$, which we gave in the previous paragraph. As a further example, we study the integral

$$\mathcal{G}_{2,4}(x,t) = \int_{-\infty}^{\infty} e^{i\lambda x - \omega_2(\lambda)t} \{\lambda \hat{g}_1(-\lambda) + [\rho(\lambda)+m]\hat{g}_2(-\lambda)\} \frac{d\lambda}{2\rho(\lambda)}.$$

To do this we write it as

$$\mathcal{G}_{2,4}(x,t) = \int_{\lambda=-1}^{1} + \int_{\lambda=-\infty}^{-1} + \int_{\lambda=1}^{\infty} = \int_{\lambda=-1}^{1} + J_3(x,t) + J_4(x,t). \tag{3.11}$$

Now

$$J_3(x,t) = \int_{-\infty}^{-1} e^{i\lambda x + i\rho(\lambda)t} \{\lambda \hat{g}_1(-\lambda) + [\rho(\lambda)+m]\hat{g}_2(-\lambda)\} \frac{d\lambda}{2\rho(\lambda)}$$
$$= \int_{-\infty}^{-1} e^{i\lambda x - i\lambda t} \sigma_3(\lambda,t) \{\lambda \hat{g}_1(-\lambda) + [\rho(\lambda)+m]\hat{g}_2(-\lambda)\} \frac{d\lambda}{2\rho(\lambda)},$$

where $\sigma_3(\lambda,t) = e^{i\lambda t + i\rho(\lambda)t}$, and

$$J_4(x,t) = \int_{1}^{\infty} e^{i\lambda x + i\rho(\lambda)t} \{\lambda \hat{g}_1(-\lambda) + [\rho(\lambda)+m]\hat{g}_2(-\lambda)\} \frac{d\lambda}{2\rho(\lambda)}$$
$$= \int_{1}^{\infty} e^{i\lambda x + i\lambda t} \sigma_4(\lambda,t) \{\lambda \hat{g}_1(-\lambda) + [\rho(\lambda)+m]\hat{g}_2(-\lambda)\} \frac{d\lambda}{2\rho(\lambda)}$$

where $\sigma_4(\lambda,t) = e^{-i\lambda t + i\rho(\lambda)t}$.

Then

$$J_3(x,t) = \int_{\Gamma^{-,-}} e^{i\lambda x - i\lambda t} \sigma_3(\lambda,t) \{\lambda \hat{g}_1(-\lambda) + [\rho(\lambda)+m]\hat{g}_2(-\lambda)\} \frac{d\lambda}{2\rho(\lambda)}, \text{ for } x<t, \tag{3.12}$$

where $\Gamma^{-,-} = \{\lambda = \xi + i\eta : \eta = \xi \ \ and \ \ \xi \leq -1\} + [-1-i, -1]$. (See fig.2.), while

$$J_3(x,t) = \sum_{k=1}^{N} \int_{\Gamma^{-,+}} e^{i\lambda x - i\lambda t} \sigma_3(\lambda,t) \frac{1}{(-i\lambda)^k} \{\lambda g_1^{(k)}(0) + [\rho(\lambda)+m] g_2^{(k)}(0)\} \frac{d\lambda}{2\rho(\lambda)} \tag{3.13}$$
$$+ \int_{-\infty}^{-1} e^{i\lambda x - i\lambda t} \sigma_3(\lambda,t) \frac{1}{(-i\lambda)^N} \{\lambda (g_1^{(N)})\hat{}(-\lambda) + [\rho(\lambda)+m](g_2^{(N)})\hat{}(-\lambda)\} \frac{d\lambda}{2\rho(\lambda)}, \text{ for } x>t,$$

and

$$J_4(x,t) = \int_{\Gamma^{+,+}} e^{i\lambda x + i\lambda t} \sigma_4(\lambda,t) \{\lambda \hat{g}_1(-\lambda) + [\rho(\lambda)+m]\hat{g}_2(-\lambda)\} \frac{d\lambda}{2\rho(\lambda)}, \text{ for } x>0, t>0. \tag{3.14}$$

It follows from (3.11) – (3.14) that the function defined by the integral $\mathcal{G}_{2,4}(x,t)$ is $C^{\infty}$ for $(x,t)\in Q-\{x=t\}$.

**(3)** *The integral* $\mathcal{H}_{1,1}$ In view of (3.3) we have: For $x>0$ and $t>0$,

$$(3.15)\quad \mathcal{H}_{1,1}(x,t)=\int_{-\infty}^{\infty} e^{i\lambda x-\omega_1(\lambda)t}\left[\int_{\tau=0}^{t} h_1(\tau)e^{\omega_1(\lambda)\tau}d\tau\right]\frac{\lambda d\lambda}{\rho(\lambda)}$$

$$=\int_{\Delta} e^{i\lambda x}\left[\sum_{k=1}^{N}(-1)^{k-1}\frac{h_1^{(k-1)}(t)}{[\omega_1(\lambda)]^k}\right]\frac{\lambda d\lambda}{\rho(\lambda)}-J_5(x,t)-(-1)^{N-1}\int_{-\infty}^{\infty}\frac{e^{i\lambda x-\omega_1(\lambda)t}}{[\omega_1(\lambda)]^N}\left[\int_{\tau=0}^{t}e^{\omega_1(\lambda)\tau}h_1^{(N)}(\tau)d\tau\right]\frac{\lambda d\lambda}{\rho(\lambda)},$$

where

$$J_5(x,t)=\int_{-\infty}^{\infty} e^{i\lambda x-\omega_1(\lambda)t}\left[\sum_{k=1}^{N}(-1)^{k-1}\frac{h_1^{(k-1)}(0)}{[\omega_1(\lambda)]^k}\right]\frac{\lambda d\lambda}{\rho(\lambda)},$$

and $\Delta=\{\lambda=\xi+i\eta:\eta=|\xi|\ and\ -\infty<\xi<+\infty\}$. (See fig. 4; other choices of $\Delta$ are shown in fig. 5, where it suffices the angles $\phi$ and $\theta$ to be $<\pi/2$; regarding the restriction to these angles, recall that the positive $\eta-axis$ contains the branch cut $[mi,i\infty)$ of the square root $\rho(\lambda)=\sqrt{\lambda^2+m^2}$.)

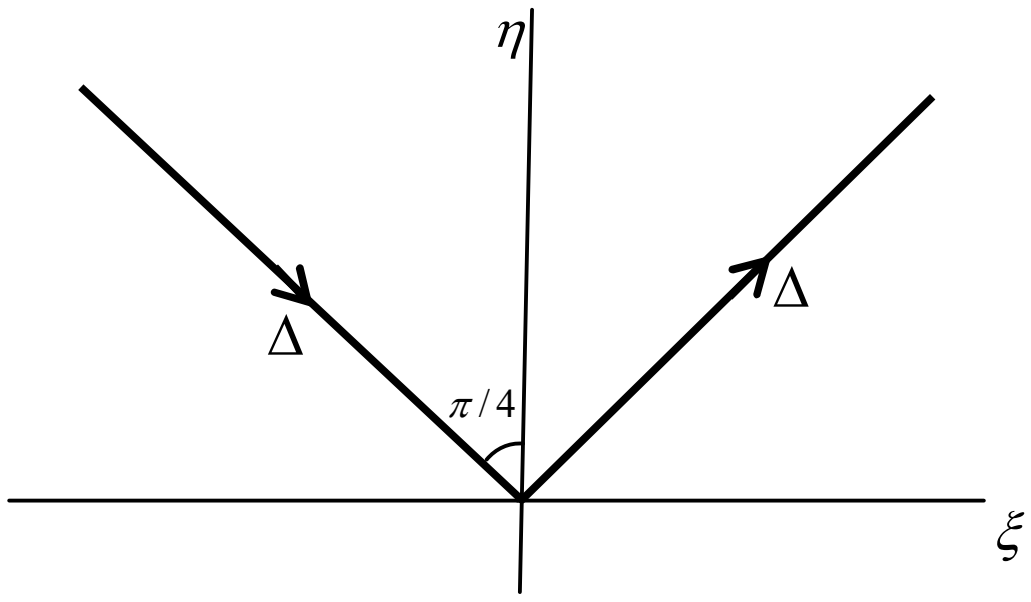

**Fig.4** *The contour* $\Delta$

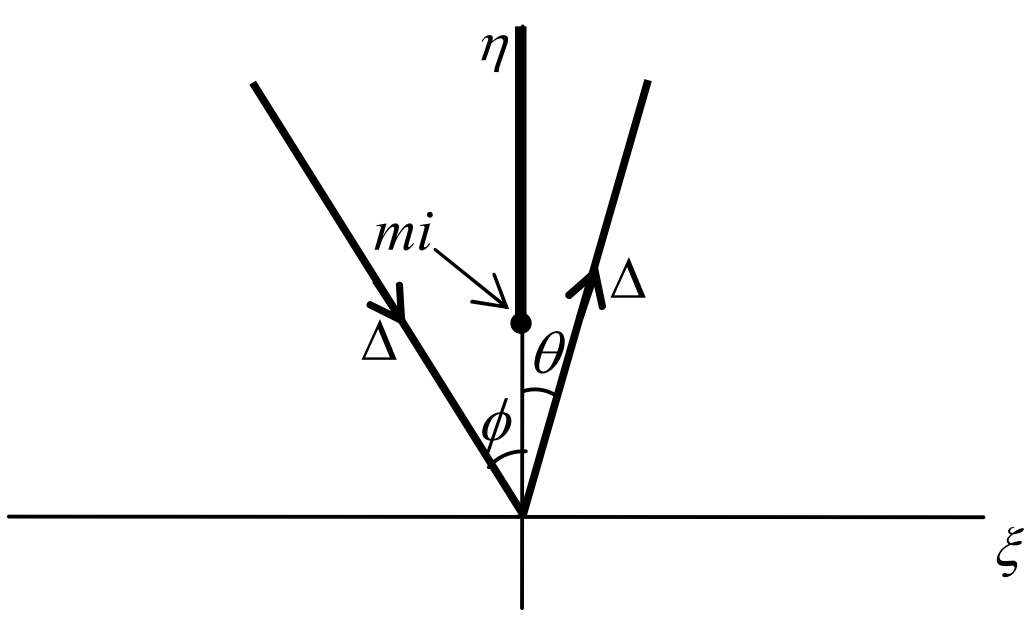

**Fig.5** *Other choices of the contour* $\Delta$

The integral $J_5(x,t)$ requires further interpretation. More precisely we write

$$(3.16)\qquad J_5(x,t)=\int_{\lambda=-\infty}^{0}+\int_{\lambda=0}^{\infty}=J_5^-(x,t)+J_5^+(x,t).$$

Then

$$(3.17)\qquad J_5^-(x,t)=\int_{-\infty}^{0} e^{i\lambda x+i\lambda t}e^{-i\lambda t-\omega_1(\lambda)t}\left[\sum_{k=1}^{N}(-1)^{k-1}\frac{h_1^{(k-1)}(0)}{[\omega_1(\lambda)]^k}\right]\frac{\lambda d\lambda}{\rho(\lambda)}$$

$$=\int_{L^{-,+}} e^{i\lambda x+i\lambda t}\left[e^{-i\lambda t-\omega_1(\lambda)t}\sum_{k=1}^{N}(-1)^{k-1}\frac{h_1^{(k-1)}(0)}{[\omega_1(\lambda)]^k}\right]\frac{\lambda d\lambda}{\rho(\lambda)}\ (x>0,\,t>0),$$

where $L^{-,+}=\{\lambda=\xi+i\eta:\eta=-\xi\ and\ -\infty<\xi\le 0\}$ (see fig. 6), while

$$(3.18)\ J_5^+(x,t)=\int_{0}^{\infty}e^{i\lambda x-i\lambda t}\left[e^{i\lambda t-\omega_1(\lambda)t}\sum_{k=1}^{N}(-1)^{k-1}\frac{h_1^{(k-1)}(0)}{[\omega_1(\lambda)]^k}\right]\frac{\lambda d\lambda}{\rho(\lambda)}=\begin{cases}\displaystyle\int_{L^{+,+}} e^{i\lambda x-i\lambda t}[\cdots]\frac{\lambda d\lambda}{\rho(\lambda)} & for\ x>t\\ \displaystyle\int_{L^{+,-}} e^{i\lambda x-i\lambda t}[\cdots]\frac{\lambda d\lambda}{\rho(\lambda)} & for\ x<t,\end{cases}$$

where

$$L^{+,+} = \{\lambda = \xi + i\eta : \eta = \xi \ \ and \ \ 0 \le \xi < +\infty\} \text{ (see fig. 6)}$$

and

$$L^{+,-} = \{\lambda = \xi + i\eta : \eta = -\xi \ \ and \ \ 0 \le \xi < +\infty\} \text{ (see fig. 7).}$$

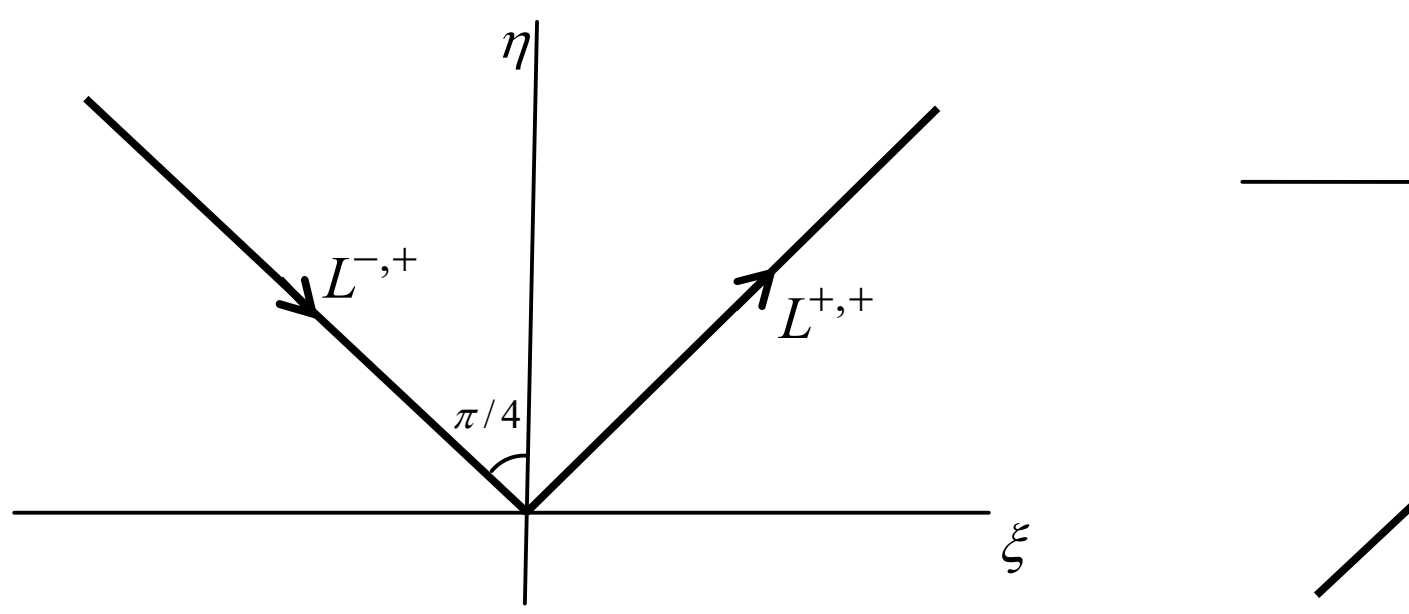

**Fig.6** *The contours* $L^{-,+}$ *and* $L^{+,+}$

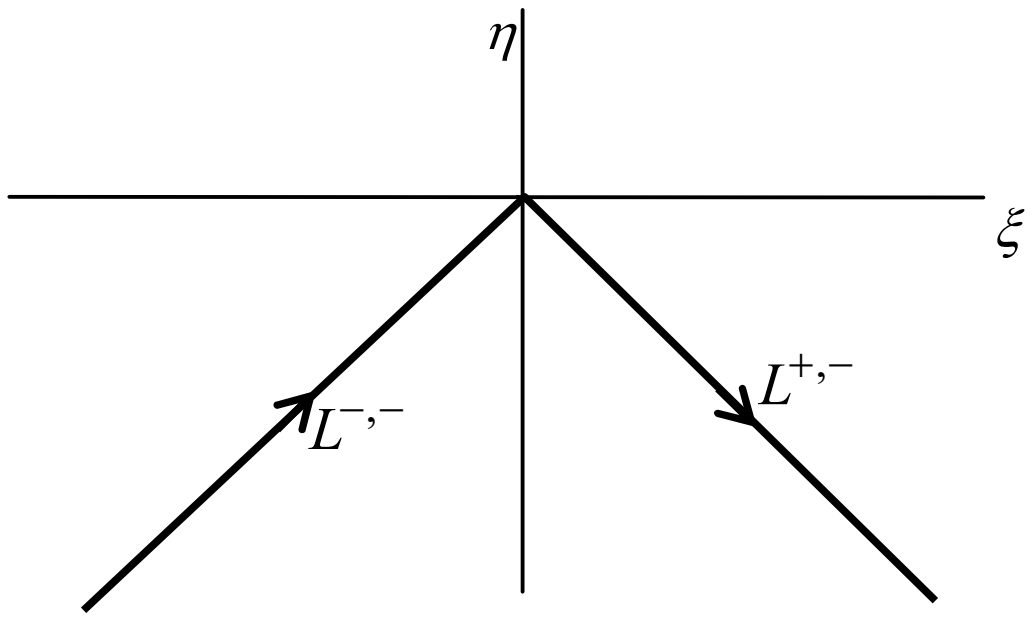

**Fig.7** *The contours* $L^{-,-}$ *and* $L^{+,-}$

It follows from (3.15) – (3.18) that the function defined by the integral $\mathcal{H}_{1,1}(x,t)$ is $C^{\infty}$ for $(x,t) \in Q - \{x = t\}$.

**(4)** *The integrals* $\mathcal{H}_{1,2}, \mathcal{H}_{2,1}, \mathcal{H}_{2,2}$ These can be analysed similarly to $\mathcal{H}_{1,1}$.

## 4. Proof of Theorem 1

**Step 1** *Verification of the differential equations of the system*

**(1)** *Basic exponential solutions* The following pairs $(\mathrm{y}_1(x,t), \mathrm{y}_2(x,t))$ satisfy the differential equations (1.1):

$$(4.1) \qquad \begin{cases} \left(\lambda e^{i\lambda x - \omega_1(\lambda)t}, -[i\omega_1(\lambda) + m]e^{i\lambda x - \omega_1(\lambda)t}\right), \left([i\omega_2(\lambda) + m]e^{i\lambda x - \omega_1(\lambda)t}, \lambda e^{i\lambda x - \omega_1(\lambda)t}\right), \\ \left(-\lambda e^{i\lambda x - \omega_2(\lambda)t}, [i\omega_2(\lambda) + m]e^{i\lambda x - \omega_2(\lambda)t}\right), \left([i\omega_1(\lambda) + m]e^{i\lambda x - \omega_2(\lambda)t}, -\lambda e^{i\lambda x - \omega_2(\lambda)t}\right), \end{cases}$$

i.e., if $(\mathrm{y}_1, \mathrm{y}_2)$ is any of the above pairs then

$$(4.2) \qquad \begin{cases} \partial_t \mathrm{y}_1 = -\partial_x \mathrm{y}_2 - im\mathrm{y}_1 \\ \partial_t \mathrm{y}_2 = -\partial_x \mathrm{y}_1 + im\mathrm{y}_2. \end{cases}$$

**(2)** *The* $\mathcal{G}-$*part of the solution* The part of the solution $(\psi_1(x,t), \psi_2(x,t))$, which involves the integrals

$$\mathcal{G}_{1,1}(x,t),\ \mathcal{G}_{1,2}(x,t), \mathcal{G}_{1,3}(x,t),\ \mathcal{G}_{1,4}(x,t),\ \mathcal{G}_{2,1}(x,t),,\ \mathcal{G}_{2,2}(x,t), \mathcal{G}_{2,3}(x,t),\ \mathcal{G}_{2,4}(x,t),$$

is a "combination" (*combination given by* $d\lambda-$*integrals*) of the pairs (4.1), and, therefore, this part satisfies the system (4.2). We point out that this holds true also for the integrals with the deformed contours, that we used for their interpretation in section 3.

It remains to prove also the corresponding result for the part of the solution which involves the integrals $\mathcal{H}_{1,1}(x,t)$, $\mathcal{H}_{1,2}(x,t)$, $\mathcal{H}_{2,1}(x,t)$, $\mathcal{H}_{2,2}(x,t)$. The point here is that, in the latter case, the integrals involve, in their integrands, apart from the pairs (4.1), also terms of the form $\int_{\tau=0}^{t} \cdots d\tau$, and this adds an extra dependence of the solution on $t$.

**(3)** *The $\mathcal{H}$ – part of the solution* We claim that, for $x>t$, the pair $(\mathrm{y}_1,\mathrm{y}_2)$, defined by the formulas

(4.3) $$\mathrm{y}_1(x,t)=\mathcal{H}_{1,1}(x,t)-\mathcal{H}_{1,2}(x,t) \text{ and } \mathrm{y}_2(x,t)=-\mathcal{H}_{2,1}(x,t)+\mathcal{H}_{2,2}(x,t),$$

satisfies (4.2).

Since we cannot differentiate the above representations of $\mathrm{y}_1$, $\mathrm{y}_2$, under the integral signs, in order to compute the derivatives $\partial_t\mathrm{y}_1,\partial_x\mathrm{y}_2,\partial_t\mathrm{y}_2,\partial_x\mathrm{y}_1$, we will use the following – equivalent – form of the integrals defining $\mathrm{y}_1$, $\mathrm{y}_2$: For $x>t$,

$$\begin{aligned}
\mathrm{y}_1(x,t)=&\int_\Delta\left\{\frac{h_1(t)}{\omega_1(\lambda)}e^{i\lambda x}-\frac{h_1(0)}{\omega_1(\lambda)}e^{i\lambda x-\omega_1(\lambda)t}-\frac{h_1'(t)}{[\omega_1(\lambda)]^2}e^{i\lambda x}+\frac{h_1'(0)}{[\omega_1(\lambda)]^2}e^{i\lambda x-\omega_1(\lambda)t}\right\}\frac{\lambda d\lambda}{\rho(\lambda)}\\
&+\int_{-\infty}^{\infty}\frac{e^{i\lambda x-\omega_1(\lambda)t}}{[\omega_1(\lambda)]^2}\left[\int_{\tau=0}^{t}h_1''(\tau)e^{\omega_1(\lambda)\tau}d\tau\right]\frac{\lambda d\lambda}{\rho(\lambda)}\\
&-\int_\Delta\left\{\frac{h_1(t)}{\omega_2(\lambda)}e^{i\lambda x}-\frac{h_1(0)}{\omega_2(\lambda)}e^{i\lambda x-\omega_2(\lambda)t}-\frac{h_1'(t)}{[\omega_2(\lambda)]^2}e^{i\lambda x}+\frac{h_1'(0)}{[\omega_2(\lambda)]^2}e^{i\lambda x-\omega_2(\lambda)t}\right\}\frac{\lambda d\lambda}{\rho(\lambda)}\\
&-\int_{-\infty}^{\infty}\frac{e^{i\lambda x-\omega_2(\lambda)t}}{[\omega_2(\lambda)]^2}\left[\int_{\tau=0}^{t}h_1''(\tau)e^{\omega_2(\lambda)\tau}d\tau\right]\frac{\lambda d\lambda}{\rho(\lambda)},
\end{aligned}$$

$$\begin{aligned}
\mathrm{y}_2(x,t)=&-\int_\Delta\left\{\frac{h_1(t)}{\omega_1(\lambda)}e^{i\lambda x}-\frac{h_1(0)}{\omega_1(\lambda)}e^{i\lambda x-\omega_1(\lambda)t}-\frac{h_1'(t)}{[\omega_1(\lambda)]^2}e^{i\lambda x}+\frac{h_1'(0)}{[\omega_1(\lambda)]^2}e^{i\lambda x-\omega_1(\lambda)t}\right\}\frac{[i\omega_1(\lambda)+m]d\lambda}{\rho(\lambda)}\\
&-\int_{0}^{\infty}\frac{e^{i\lambda x-\omega_1(\lambda)t}}{[\omega_1(\lambda)]^2}\left[\int_{\tau=0}^{t}h_1''(\tau)e^{\omega_1(\lambda)\tau}d\tau\right]\frac{[i\omega_1(\lambda)+m]d\lambda}{\rho(\lambda)}\\
&+\int_\Delta\left\{\frac{h_1(t)}{\omega_2(\lambda)}e^{i\lambda x}-\frac{h_1(0)}{\omega_2(\lambda)}e^{i\lambda x-\omega_2(\lambda)t}-\frac{h_1'(t)}{[\omega_2(\lambda)]^2}e^{i\lambda x}+\frac{h_1'(0)}{[\omega_2(\lambda)]^2}e^{i\lambda x-\omega_2(\lambda)t}\right\}\frac{[i\omega_2(\lambda)+m]d\lambda}{\rho(\lambda)}\\
&-\int_{0}^{\infty}\frac{e^{i\lambda x-\omega_1(\lambda)t}}{[\omega_2(\lambda)]^2}\left[\int_{\tau=0}^{t}h_1''(\tau)e^{\omega_2(\lambda)\tau}d\tau\right]\frac{[i\omega_2(\lambda)+m]d\lambda}{\rho(\lambda)}.
\end{aligned}$$

Now we can compute: For $x>t$,

$$\begin{aligned}
\frac{\partial\mathrm{y}_1(x,t)}{\partial t}=&\int_\Delta\left\{\frac{h_1'(t)}{\omega_1(\lambda)}e^{i\lambda x}+h_1(0)e^{i\lambda x-\omega_1(\lambda)t}-\frac{h_1''(t)}{[\omega_1(\lambda)]^2}e^{i\lambda x}-\frac{h_1'(0)}{\omega_1(\lambda)}e^{i\lambda x-\omega_1(\lambda)t}\right\}\frac{\lambda d\lambda}{\rho(\lambda)}\\
&-\int_{-\infty}^{\infty}\frac{e^{i\lambda x-\omega_1(\lambda)t}}{\omega_1(\lambda)}\left[\int_{0}^{t}h_1''(\tau)e^{\omega_1(\lambda)\tau}d\tau\right]\frac{\lambda d\lambda}{\rho(\lambda)}+\int_\Delta\frac{h_1''(t)}{[\omega_1(\lambda)]^2}e^{i\lambda x}\frac{\lambda d\lambda}{\rho(\lambda)}\\
&-\int_\Delta\left\{\frac{h_1'(t)}{\omega_2(\lambda)}e^{i\lambda x}+h_1(0)e^{i\lambda x-\omega_2(\lambda)t}-\frac{h_1''(t)}{[\omega_2(\lambda)]^2}e^{i\lambda x}-\frac{h_1'(0)}{\omega_2(\lambda)}e^{i\lambda x-\omega_2(\lambda)t}\right\}\frac{\lambda d\lambda}{\rho(\lambda)}\\
&+\int_{-\infty}^{\infty}\frac{e^{i\lambda x-\omega_2(\lambda)t}}{\omega_2(\lambda)}\left[\int_{0}^{t}h_1''(\tau)e^{\omega_2(\lambda)\tau}d\tau\right]\frac{\lambda d\lambda}{\rho(\lambda)}-\int_\Delta\frac{h_1''(t)}{\omega_2(\lambda)}e^{i\lambda x}\frac{\lambda d\lambda}{\rho(\lambda)},
\end{aligned}$$

$$\frac{\partial\mathrm{y}_1(x,t)}{\partial x}=\int_\Delta i\lambda\left\{\frac{h_1(t)}{\omega_1(\lambda)}e^{i\lambda x}-\frac{h_1(0)}{\omega_1(\lambda)}e^{i\lambda x-\omega_1(\lambda)t}-\frac{h_1'(t)}{[\omega_1(\lambda)]^2}e^{i\lambda x}+\frac{h_1'(0)}{[\omega_1(\lambda)]^2}e^{i\lambda x-\omega_1(\lambda)t}\right\}\frac{\lambda d\lambda}{\rho(\lambda)}$$

$$+\int_{-\infty}^{\infty} i\lambda \frac{e^{i\lambda x-\omega_1(\lambda)t}}{[\omega_1(\lambda)]^2}\left[\int_{\tau=0}^{t} h_1''(\tau)e^{\omega_1(\lambda)\tau}d\tau\right]\frac{\lambda d\lambda}{\rho(\lambda)}$$

$$-\int_{\Delta} i\lambda\left\{\frac{h_1(t)}{\omega_2(\lambda)}e^{i\lambda x}-\frac{h_1(0)}{\omega_2(\lambda)}e^{i\lambda x-\omega_2(\lambda)t}-\frac{h_1'(t)}{[\omega_2(\lambda)]^2}e^{i\lambda x}+\frac{h_1'(0)}{[\omega_2(\lambda)]^2}e^{i\lambda x-\omega_2(\lambda)t}\right\}\frac{\lambda d\lambda}{\rho(\lambda)}$$

$$-\int_{-\infty}^{\infty} i\lambda \frac{e^{i\lambda x-\omega_2(\lambda)t}}{[\omega_2(\lambda)]^2}\left[\int_{\tau=0}^{t} h_1''(\tau)e^{\omega_2(\lambda)\tau}d\tau\right]\frac{\lambda d\lambda}{\rho(\lambda)},$$

and similar formulas for the derivatives $\dfrac{\partial y_2(x,t)}{\partial t}$ and $\dfrac{\partial y_2(x,t)}{\partial x}$.

Using the above formulas, it is straightforward to verify that

$$\frac{\partial y_1(x,t)}{\partial t}+\frac{\partial y_2(x,t)}{\partial x}+imy_1(x,t)=0$$

and

$$\frac{\partial y_2(x,t)}{\partial t}+\frac{\partial y_1(x,t)}{\partial x}-imy_2(x,t)=\int_{\Delta} 2h_1(t)e^{i\lambda x}d\lambda=0.$$

It follows that, indeed, ($y_1, y_2$), defined by (4.3), satisfies (4.2), in the case $x>t$. The case $x<t$ is similar.

**Step 2** *The limits* $\lim_{t\to 0^+}\psi_1(x,t)=g_1(x)$ *and* $\lim_{t\to 0^+}\psi_2(x,t)=g_2(x)$

Now we work with fixed $x>0$ and for $t<x$.

**(1)** We claim that

(4.3) $$\lim_{t\to 0^+}\mathcal{G}_{1,1}(x,t)=\mathcal{G}_{1,1}(x,t)\Big|_{t=0}.$$

In order to prove this, we write, for $t<x$,

$$\mathcal{G}_{1,1}(x,t)=\int_{-1}^{1} e^{i\lambda x-\omega_1(\lambda)t}\{\lambda\hat{g}_2(\lambda)+[i\omega_2(\lambda)+m]\hat{g}_1(\lambda)\}\frac{d\lambda}{2\rho(\lambda)}$$

$$+\int_{\Gamma^{+,+}+\Gamma^{-,+}} e^{i\lambda x-\omega_1(\lambda)t}\left\{\lambda\frac{g_2(0)}{i\lambda}+[i\omega_2(\lambda)+m]\frac{g_1(0)}{i\lambda}\right\}\frac{d\lambda}{2\rho(\lambda)}$$

$$+\left(\int_{-\infty}^{-1}+\int_{-\infty}^{-1}\right)e^{i\lambda x-\omega_1(\lambda)t}\left\{\lambda\frac{(g_2')\hat{}(\lambda)}{i\lambda}+[i\omega_2(\lambda)+m]\frac{(g_2')\hat{}(\lambda)}{i\lambda}\right\}\frac{d\lambda}{2\rho(\lambda)},$$

whence we obtain

$$\lim_{t\to 0^+}\mathcal{G}_{1,1}(x,t)=\int_{-1}^{1} e^{i\lambda x}\{\lambda\hat{g}_2(\lambda)+[i\omega_2(\lambda)+m]\hat{g}_1(\lambda)\}\frac{d\lambda}{2\rho(\lambda)}$$

$$+\int_{\Gamma^{+,+}+\Gamma^{-,+}} e^{i\lambda x}\left\{\lambda\frac{g_2(0)}{i\lambda}+[i\omega_2(\lambda)+m]\frac{g_1(0)}{i\lambda}\right\}\frac{d\lambda}{2\rho(\lambda)}$$

$$+\left(\int_{-\infty}^{-1}+\int_{-\infty}^{-1}\right)e^{i\lambda x}\left\{\lambda\frac{(g_2')\hat{}(\lambda)}{i\lambda}+[i\omega_2(\lambda)+m]\frac{(g_2')\hat{}(\lambda)}{i\lambda}\right\}\frac{d\lambda}{2\rho(\lambda)}$$

$$= \int_{-\infty}^{\infty} e^{i\lambda x}\{\lambda\hat{g}_2(\lambda)+[i\omega_2(\lambda)+m]\hat{g}_1(\lambda)\}\frac{d\lambda}{2\rho(\lambda)}.$$

This proves (4.3).

(The point of using the analysis of $\hat{g}_2(\lambda)$ and $\hat{g}_1(\lambda)$, as this is given by (2.1), is that the integral $\mathcal{G}_{1,1}(x,t)\big|_{t=0}$, although it exists in the generalized sence, it does not – in general – converge absolutely. Thus, Lebesgue's dominated convergence theorem could not be applied immediately.)

**(2)** Similarly, $\lim_{t\to 0^+}\mathcal{G}_{1,3}(x,t)=\mathcal{G}_{1,3}(x,t)\big|_{t=0}$, and, therefore

(4.4) $\lim_{t\to 0^+}[\mathcal{G}_{1,1}(x,t)-\mathcal{G}_{1,3}(x,t)]=[\mathcal{G}_{1,1}(x,t)-\mathcal{G}_{1,3}(x,t)]\big|_{t=0}$

$$= \int_{-\infty}^{\infty} e^{i\lambda x}\{\lambda\hat{g}_2(\lambda)+[i\omega_2(\lambda)+m]\hat{g}_1(\lambda)\}\frac{d\lambda}{2\rho(\lambda)} - \int_{-\infty}^{\infty} e^{i\lambda x}\{\lambda\hat{g}_2(\lambda)+[i\omega_1(\lambda)+m]\hat{g}_1(\lambda)\}\frac{d\lambda}{2\rho(\lambda)}$$

$$= \int_{-\infty}^{\infty} e^{i\lambda x}[i\omega_2(\lambda)-i\omega_1(\lambda)]\hat{g}_1(\lambda)\}\frac{d\lambda}{2\rho(\lambda)} = \int_{-\infty}^{\infty} e^{i\lambda x}\hat{g}_1(\lambda)d\lambda = 2\pi g_1(x).$$

**(3)** Working similarly we see that

(4.5) $\lim_{t\to 0^+}[\mathcal{G}_{1,2}(x,t)-\mathcal{G}_{1,4}(x,t)]=[\mathcal{G}_{1,2}(x,t)-\mathcal{G}_{1,4}(x,t)]\big|_{t=0}$

$$= \int_{-\infty}^{\infty} e^{i\lambda x}\{\lambda\hat{g}_1(-\lambda)+[i\omega_1(\lambda)+m]\hat{g}_2(-\lambda)\}\frac{d\lambda}{2\rho(\lambda)} - \int_{-\infty}^{\infty} e^{i\lambda x}\{\lambda\hat{g}_1(-\lambda)+[i\omega_2(\lambda)+m]\hat{g}_2(-\lambda)\}\frac{d\lambda}{2\rho(\lambda)}$$

$$= \int_{-\infty}^{\infty} e^{i\lambda x}[i\omega_1(\lambda)-i\omega_2(\lambda)]\hat{g}_2(-\lambda)\}\frac{d\lambda}{2\rho(\lambda)} = -\int_{-\infty}^{\infty} e^{i\lambda x}\hat{g}_2(-\lambda)d\lambda = 0.$$

(The last equation follows from Cauchy's theorem and Jordan's lemma.)

**(4)** We claim that

(4.6) $$\lim_{t\to 0^+}[\mathcal{H}_{1,1}(x,t)-\mathcal{H}_{1,2}(x,t)]=0.$$

This follows immediately from the formula

$$\mathcal{H}_{1,1}(x,t)-\mathcal{H}_{1,2}(x,t)=\int_{\Delta}\left\{\frac{h_1(t)}{\omega_1(\lambda)}e^{i\lambda x}-\frac{h_1(0)}{\omega_1(\lambda)}e^{i\lambda x-\omega_1(\lambda)t}\right\}\frac{\lambda d\lambda}{\rho(\lambda)}$$

$$+\int_{-\infty}^{\infty}\frac{e^{i\lambda x-\omega_1(\lambda)t}}{\omega_1(\lambda)}\left[\int_{\tau=0}^{t}h_1'(\tau)e^{\omega_1(\lambda)\tau}d\tau\right]\frac{\lambda d\lambda}{\rho(\lambda)}-\int_{\Delta}\left\{\frac{h_1(t)}{\omega_2(\lambda)}e^{i\lambda x}-\frac{h_1(0)}{\omega_2(\lambda)}e^{i\lambda x-\omega_2(\lambda)t}\right\}\frac{\lambda d\lambda}{\rho(\lambda)}$$

$$-\int_{-\infty}^{\infty}\frac{e^{i\lambda x-\omega_2(\lambda)t}}{\omega_2(\lambda)}\left[\int_{\tau=0}^{t}h_1'(\tau)e^{\omega_2(\lambda)\tau}d\tau\right]\frac{\lambda d\lambda}{\rho(\lambda)}.$$

Thus $\lim_{t\to 0^+}\psi_1(x,t)=g_1(x)$ follows from (4.3), (4.4), (4.5) and (4.6).

The proof of $\lim_{t\to 0^+}\psi_2(x,t)=g_2(x)$ is similar.

**Step 3** *The limit* $\lim_{x\to 0^+}\psi_1(x,t)=h_1(t)$

Now we work with fixed $t>0$ and for $x<t$.

**(1)** We claim that

(4.7) $$\lim_{x\to 0^+}[\mathcal{G}_{1,1}(x,t)+\mathcal{G}_{1,2}(x,t)]=0 \text{ and } \lim_{x\to 0^+}[\mathcal{G}_{1,3}(x,t)+\mathcal{G}_{1,4}(x,t)]=0.$$

Indeed, working with $x<t$ and using (3.1), we obtain

$$\lim_{x\to 0^+}[\mathcal{G}_{1,1}(x,t)+\mathcal{G}_{1,2}(x,t)]=\lim_{x\to 0^+}\int_{-1}^{1} e^{i\lambda x-\omega_1(\lambda)t}\{\lambda\hat{g}_2(\lambda)+[i\omega_2(\lambda)+m]\hat{g}_1(\lambda)\}\frac{d\lambda}{2\rho(\lambda)}$$

$$+\lim_{x\to 0^+}\int_{\Gamma^{+,-}+\Gamma^{-,+}} e^{i\lambda x-\omega_1(\lambda)t}\left\{\lambda\frac{g_2(0)}{i\lambda}+[i\omega_2(\lambda)+m]\frac{g_1(0)}{i\lambda}\right\}\frac{d\lambda}{2\rho(\lambda)}$$

$$+\lim_{x\to 0^+}\int_{-\infty}^{\infty} e^{i\lambda x-\omega_1(\lambda)t}\left\{\lambda\frac{(g_2')\hat{}(\lambda)}{i\lambda}+[i\omega_2(\lambda)+m]\frac{(g_1')\hat{}(\lambda)}{i\lambda}\right\}\frac{d\lambda}{2\rho(\lambda)}$$

$$+\lim_{x\to 0^+}\int_{-1}^{1} e^{i\lambda x-\omega_1(\lambda)t}\{\lambda\hat{g}_2(-\lambda)-[i\omega_2(\lambda)+m]\hat{g}_1(-\lambda)\}\frac{d\lambda}{2\rho(\lambda)}$$

$$+\lim_{x\to 0^+}\int_{\Gamma^{+,-}+\Gamma^{-,+}} e^{i\lambda x-\omega_1(\lambda)t}\left\{\lambda\frac{g_2(0)}{-i\lambda}-[i\omega_2(\lambda)+m]\frac{g_1(0)}{-i\lambda}\right\}\frac{d\lambda}{2\rho(\lambda)}$$

$$+\lim_{x\to 0^+}\int_{-\infty}^{\infty} e^{i\lambda x-\omega_1(\lambda)t}\left\{\lambda\frac{(g_2')\hat{}(-\lambda)}{-i\lambda}-[i\omega_2(\lambda)+m]\frac{(g_1')\hat{}(-\lambda)}{-i\lambda}\right\}\frac{d\lambda}{2\rho(\lambda)}$$

$$=\int_{-1}^{1} e^{-\omega_1(\lambda)t}\{\lambda\hat{g}_2(\lambda)+[i\omega_2(\lambda)+m]\hat{g}_1(\lambda)\}\frac{d\lambda}{2\rho(\lambda)}$$

$$+\int_{\Gamma^{+,-}+\Gamma^{-,+}} e^{-\omega_1(\lambda)t}\left\{\lambda\frac{g_2(0)}{i\lambda}+[i\omega_2(\lambda)+m]\frac{g_1(0)}{i\lambda}\right\}\frac{d\lambda}{2\rho(\lambda)}$$

$$+\int_{-\infty}^{\infty} e^{-\omega_1(\lambda)t}\left\{\lambda\frac{(g_2')\hat{}(\lambda)}{i\lambda}+[i\omega_2(\lambda)+m]\frac{(g_1')\hat{}(\lambda)}{i\lambda}\right\}\frac{d\lambda}{2\rho(\lambda)}$$

$$+\int_{-1}^{1} e^{-\omega_1(\lambda)t}\{\lambda\hat{g}_2(-\lambda)-[i\omega_2(\lambda)+m]\hat{g}_1(-\lambda)\}\frac{d\lambda}{2\rho(\lambda)}$$

$$+\int_{\Gamma^{+,-}+\Gamma^{-,+}} e^{-\omega_1(\lambda)t}\left\{\lambda\frac{g_2(0)}{-i\lambda}-[i\omega_2(\lambda)+m]\frac{g_1(0)}{-i\lambda}\right\}\frac{d\lambda}{2\rho(\lambda)}$$

$$+\int_{-\infty}^{\infty} e^{-\omega_1(\lambda)t}\left\{\lambda\frac{(g_2')\hat{}(-\lambda)}{-i\lambda}-[i\omega_2(\lambda)+m]\frac{(g_1')\hat{}(-\lambda)}{-i\lambda}\right\}\frac{d\lambda}{2\rho(\lambda)}$$

$$=\int_{-\infty}^{\infty} e^{-\omega_1(\lambda)t}\{\lambda\hat{g}_2(\lambda)+[i\omega_2(\lambda)+m]\hat{g}_1(\lambda)\}\frac{d\lambda}{2\rho(\lambda)}$$

$$+\int_{-\infty}^{\infty} e^{-\omega_1(\lambda)t}\{\lambda\hat{g}_2(-\lambda)-[i\omega_2(\lambda)+m]\hat{g}_1(-\lambda)\}\frac{d\lambda}{2\rho(\lambda)}=0.$$

The proof of the second equation in (4.7) is similar.

**(2)** We claim that

(4.8) $$\lim_{x\to 0^+}[\mathcal{H}_{1,1}(x,t)-\mathcal{H}_{1,2}(x,t)]=2\pi h_1(t).$$

In this part we work with $x<t$. For $x>0$, we have

(4.9) $$\mathcal{H}_{1,1}(x,t)-\mathcal{H}_{1,2}(x,t)=\int_{-\infty}^{\infty}\frac{h_1(t)}{\omega_1(\lambda)}e^{i\lambda x}\frac{\lambda d\lambda}{\rho(\lambda)}-\int_{-\infty}^{\infty}\frac{h_1(0)}{\omega_1(\lambda)}e^{i\lambda x-\omega_1(\lambda)t}\frac{\lambda d\lambda}{\rho(\lambda)}$$
$$-\int_{-\infty}^{\infty}\frac{e^{i\lambda x-\omega_1(\lambda)t}}{\omega_1(\lambda)}\left[\int_{\tau=0}^{t}h_1'(\tau)e^{\omega_1(\lambda)\tau}d\tau\right]\frac{\lambda d\lambda}{\rho(\lambda)}$$
$$-\int_{-\infty}^{\infty}\frac{h_1(t)}{\omega_2(\lambda)}e^{i\lambda x}\frac{\lambda d\lambda}{\rho(\lambda)}+\int_{-\infty}^{\infty}\frac{h_1(0)}{\omega_2(\lambda)}e^{i\lambda x-\omega_2(\lambda)t}\frac{\lambda d\lambda}{\rho(\lambda)}+\int_{-\infty}^{\infty}\frac{e^{i\lambda x-\omega_2(\lambda)t}}{\omega_2(\lambda)}\left[\int_{\tau=0}^{t}h_1'(\tau)e^{\omega_2(\lambda)\tau}d\tau\right]\frac{\lambda d\lambda}{\rho(\lambda)}.$$

(We recall the each of the above integrals exist – at least – in the generalized sence.)
Now,

(4.10) $$\lim_{x\to 0^+}\int_{-\infty}^{\infty}\frac{h_1(t)}{\omega_1(\lambda)}e^{i\lambda x}\frac{\lambda d\lambda}{\rho(\lambda)}=\lim_{x\to 0^+}\int_{-\infty}^{\infty}\frac{h_1(t)}{i(\lambda^2+m^2)}e^{i\lambda x}\lambda d\lambda$$
$$=(with\ \varepsilon<m)\lim_{x\to 0^+}\int_{|\lambda-mi|=\varepsilon}\frac{h_1(t)}{i(\lambda^2+m^2)}e^{i\lambda x}\lambda d\lambda=-ih_1(t)\int_{|\lambda-mi|=\varepsilon}\frac{\lambda}{\lambda^2+m^2}d\lambda,$$
$$=-ih_1(t)\left[2\pi i Res\left(\frac{\lambda}{\lambda^2+m^2},\lambda=mi\right)\right]=\pi h_1(t),$$

and, similarly,

(4.11) $$\lim_{x\to 0^+}\int_{-\infty}^{\infty}\frac{h_1(t)}{\omega_2(\lambda)}e^{i\lambda x}\frac{\lambda d\lambda}{\rho(\lambda)}=-\pi h_1(t).$$

On the other hand,

(4.12) $$\lim_{x\to 0^+}\int_{-\infty}^{\infty}\frac{h_1(0)}{\omega_1(\lambda)}e^{i\lambda x-\omega_1(\lambda)t}\frac{\lambda d\lambda}{\rho(\lambda)}$$
$$=\lim_{x\to 0^+}\int_{-\infty}^{0}\frac{h_1(0)}{\omega_1(\lambda)}e^{i\lambda x-\omega_1(\lambda)t}\frac{\lambda d\lambda}{\rho(\lambda)}+\lim_{x\to 0^+}\int_{0}^{\infty}\frac{h_1(0)}{\omega_1(\lambda)}e^{i\lambda x-\omega_1(\lambda)t}\frac{\lambda d\lambda}{\rho(\lambda)}$$
$$=\lim_{x\to 0^+}\int_{L^{-,+}}\frac{h_1(0)}{\omega_1(\lambda)}e^{i\lambda x-\omega_1(\lambda)t}\frac{\lambda d\lambda}{\rho(\lambda)}+\lim_{x\to 0^+}\int_{L^{+,-}}\frac{h_1(0)}{\omega_1(\lambda)}e^{i\lambda x-\omega_1(\lambda)t}\frac{\lambda d\lambda}{\rho(\lambda)}$$
$$=\int_{L^{-,+}}\frac{h_1(0)}{\omega_1(\lambda)}e^{-\omega_1(\lambda)t}\frac{\lambda d\lambda}{\rho(\lambda)}+\int_{L^{+,-}}\frac{h_1(0)}{\omega_1(\lambda)}e^{-\omega_1(\lambda)t}\frac{\lambda d\lambda}{\rho(\lambda)}$$
$$=\int_{-\infty}^{0}\frac{h_1(0)}{\omega_1(\lambda)}e^{-\omega_1(\lambda)t}\frac{\lambda d\lambda}{\rho(\lambda)}+\int_{0}^{\infty}\frac{h_1(0)}{\omega_1(\lambda)}e^{-\omega_1(\lambda)t}\frac{\lambda d\lambda}{\rho(\lambda)}=0,$$

and, similarly,

(4.13) $$\lim_{x\to 0^+}\int_{-\infty}^{\infty}\frac{h_1(0)}{\omega_2(\lambda)}e^{i\lambda x-\omega_2(\lambda)t}\frac{\lambda d\lambda}{\rho(\lambda)}=0.$$

Finally, it is easier to see that

(4.14) $$\lim_{x\to 0^+}\int_{-\infty}^{\infty}\frac{e^{i\lambda x-\omega_j(\lambda)t}}{\omega_j(\lambda)}\left[\int_{\tau=0}^{t}h_1'(\tau)e^{\omega_j(\lambda)\tau}d\tau\right]\frac{\lambda d\lambda}{\rho(\lambda)}=0,$$

since the integrand, this time, is $\mathrm{O}(1/\lambda^2)$, as $\lambda\to\pm\infty$.

Thus, the equation $\lim_{x\to 0^+}\psi_1(x,t)=h_1(t)$ follows from (4.8) – (4.14).

*Comment* We have been persistant in justifying the various interchanges of limiting and differentiating processes with integration processes, since otherwise one can be led to erroneous conclusions. For example, if one wrote the limit

$$\lim_{x\to 0^+}\int_{-\infty}^{\infty}\frac{h_1(t)}{\omega_1(\lambda)}e^{i\lambda x}\frac{\lambda d\lambda}{\rho(\lambda)}$$

as

$$\int_{-\infty}^{\infty}\lim_{x\to 0^+}\frac{h_1(t)}{\omega_1(\lambda)}e^{i\lambda x}\frac{\lambda d\lambda}{\rho(\lambda)}$$

(swithing the order of the limit and the integral), then this would be equal to

$$\int_{-\infty}^{\infty}\frac{h_1(t)}{\omega_1(\lambda)}\frac{\lambda d\lambda}{\rho(\lambda)}=0$$

(since the integrand is an odd function of $\lambda$), which, in view of (4.10), we know to be wrong.

**Step 4** *The limits* $\lim_{\substack{(x,t)\to(a,a)\\ x\neq t}}\psi_1(x,t)$ *and* $\lim_{\substack{(x,t)\to(a,a)\\ x\neq t}}\psi_2(x,t)$

We study these limits with fixed $a>0$.

**(1)** We claim that if $g_1(0)=h_1(0)$ then the limits

(4.15) $$\lim_{\substack{(x,t)\to(a,a)\\ x\neq t}}[\mathcal{G}_{1,1}(x,t)+\mathcal{G}_{1,2}(x,t)+\mathcal{H}_{1,1}(x,t)] \text{ and } \lim_{\substack{(x,t)\to(a,a)\\ x\neq t}}[\mathcal{G}_{1,3}(x,t)+\mathcal{G}_{1,4}(x,t)+\mathcal{H}_{1,2}(x,t)]$$

exist, for every $a>0$.

**(2)** *Notation* For two functions $\mathcal{A}(x,t)$ and $\mathcal{B}(x,t)$, defined for $x\neq t$, we will write

$$\mathcal{A}(x,t)\approx\mathcal{B}(x,t)$$

if and only if the limit $\lim_{\substack{(x,t)\to(a,a)\\ x\neq t}}[\mathcal{A}(x,t)-\mathcal{B}(x,t)]$ exists for every $a>0$.

Equivalently, $\mathcal{A}(x,t)\approx\mathcal{B}(x,t)$ if and only if the limits

$$\lim_{\substack{(x,t)\to(a,a)\\ x>t}}[\mathcal{A}(x,t)-\mathcal{B}(x,t)] \text{ and } \lim_{\substack{(x,t)\to(a,a)\\ x<t}}[\mathcal{A}(x,t)-\mathcal{B}(x,t)]$$

exist and they are equal, for every $a>0$.

For example,

$$\int_{-1}^{1}e^{i\lambda x-\omega_1(\lambda)t}\{\lambda\hat{g}_2(-\lambda)-[i\omega_2(\lambda)+m]\hat{g}_1(-\lambda)\}\frac{d\lambda}{2\rho(\lambda)}\approx 0$$

and

$$\int_0^\infty e^{i\lambda x-\omega_1(\lambda)t}\{\lambda\hat{g}_2(\lambda)+[i\omega_2(\lambda)+m]\hat{g}_1(\lambda)\}\frac{d\lambda}{2\rho(\lambda)}\approx\int_0^\infty e^{i\lambda x-\omega_1(\lambda)t}\{\lambda\hat{g}_2(\lambda)+i\omega_2(\lambda)\hat{g}_1(\lambda)\}\frac{d\lambda}{2\rho(\lambda)}.$$

**(3)** *Proof of the existence of the limits (4.15)* Expanding $\hat{g}_2(\lambda)$ and $\hat{g}_1(\lambda)$, according to the formula (3.1), we easily find that

$$\begin{aligned}\mathcal{G}_{1,1}^+(x,t)+\mathcal{G}_{1,2}^+(x,t)&\approx\int_1^\infty e^{i\lambda x-\omega_1(\lambda)t}\{\lambda\hat{g}_2(\lambda)+[i\omega_2(\lambda)+m]\hat{g}_1(\lambda)\}\frac{d\lambda}{2\rho(\lambda)}\\&\quad+\int_1^\infty e^{i\lambda x-\omega_1(\lambda)t}\{\lambda\hat{g}_2(-\lambda)-[i\omega_2(\lambda)+m]\hat{g}_1(-\lambda)\}\frac{d\lambda}{2\rho(\lambda)}\\&\approx\int_1^\infty e^{i\lambda x-\omega_1(\lambda)t}\left\{i\omega_2(\lambda)\frac{g_1(0)}{i\lambda}-i\omega_2(\lambda)\frac{g_1(0)}{-i\lambda}\right\}\frac{d\lambda}{2\rho(\lambda)}\approx-ig_1(0)\int_1^\infty e^{i\lambda x-\omega_1(\lambda)t}\frac{d\lambda}{\rho(\lambda)}.\end{aligned}$$

(In proving the above relation, in the process of $(x,t)$ tending toward the point $(a,a)$, we have to distinguish between the cases " $x>t$ " and " $x<t$ ".)
On the other hand, using formula (3.3), we find

$$\begin{aligned}\mathcal{H}_{1,1}^+(x,t)&=\int_0^\infty e^{i\lambda x-\omega_1(\lambda)t}\left[\int_0^t h_1(\tau)e^{\omega_1(\lambda)\tau}d\tau\right]\frac{\lambda d\lambda}{\rho(\lambda)}\\&\approx-\int_0^\infty\frac{h_1(0)}{\omega_1(\lambda)}e^{i\lambda x-\omega_1(\lambda)t}\frac{\lambda d\lambda}{\rho(\lambda)}\approx ih_1(0)\int_1^\infty e^{i\lambda x-\omega_1(\lambda)t}\frac{d\lambda}{\rho(\lambda)}.\end{aligned}$$

Therefore,

$$\mathcal{G}_{1,1}^+(x,t)+\mathcal{G}_{1,2}^+(x,t)+\mathcal{H}_{1,1}^+(x,t)\approx[ih_1(0)-ig_1(0)]\int_1^\infty e^{i\lambda x-\omega_1(\lambda)t}\frac{d\lambda}{\rho(\lambda)}.\tag{4.16}$$

Similar computations show that

$$\mathcal{G}_{1,1}^-(x,t)+\mathcal{G}_{1,2}^-(x,t)+\mathcal{H}_{1,1}^-(x,t)\approx 0.\tag{4.17}$$

Also

$$\mathcal{G}_{1,3}^+(x,t)+\mathcal{G}_{1,4}^+(x,t)+\mathcal{H}_{1,2}^+(x,t)\approx 0\tag{4.18}$$

and

$$\mathcal{G}_{1,3}^-(x,t)+\mathcal{G}_{1,4}^-(x,t)+\mathcal{H}_{1,2}^-(x,t)\approx[ih_1(0)-ig_1(0)]\int_{-\infty}^{-1} e^{i\lambda x-\omega_2(\lambda)t}\frac{d\lambda}{\rho(\lambda)}.\tag{4.19}$$

Thus, the existence of the limits (4.15), follow from (4.16) – (4.19), and this proves that

$$g_1(0)=h_1(0)\quad\Rightarrow\quad\lim_{\substack{(x,t)\to(a,a)\\x\neq t}}\psi_1(x,t)\ \text{ exists for every } a>0.\tag{4.20}$$

**(4)** Similar computations show that

$$\mathcal{G}_{2,1}^+(x,t)-\mathcal{G}_{2,2}^+(x,t)-\mathcal{H}_{2,1}^+(x,t)\approx[-ig_1(0)+ih_1(0)]\int_1^\infty e^{i\lambda x-\omega_1(\lambda)t}\frac{d\lambda}{\rho(\lambda)},$$

$$\mathcal{G}_{2,1}^{-}(x,t)-\mathcal{G}_{2,2}^{-}(x,t)-\mathcal{H}_{2,1}^{-}(x,t)\approx 0,$$

$$-\mathcal{G}_{2,3}^{+}(x,t)+\mathcal{G}_{2,4}^{+}(x,t)+\mathcal{H}_{2,2}^{+}(x,t)\approx 0,$$

$$-\mathcal{G}_{2,3}^{-}(x,t)+\mathcal{G}_{2,4}^{-}(x,t)+\mathcal{H}_{2,2}^{-}(x,t)\approx[-ig_1(0)+ih_1(0)]\int_{-\infty}^{-1} e^{i\lambda x-\omega_2(\lambda)t}\frac{d\lambda}{\rho(\lambda)},$$

and, therefore,

(4.21) $$g_1(0)=h_1(0)\ \Rightarrow \lim_{\substack{(x,t)\to(a,a)\\ x\neq t}}\psi_2(x,t)\ \text{exists for every}\ a>0.$$

***Step 5*** *The limits* $\displaystyle\lim_{\substack{(x,t)\to(a,a)\\ x\neq t}}\frac{\partial\psi_j(x,t)}{\partial x}$ *and* $\displaystyle\lim_{\substack{(x,t)\to(a,a)\\ x\neq t}}\frac{\partial\psi_j(x,t)}{\partial t}$ $(j=1,2)$

We study these limits with fixed $a>0$.

Computations as in *Step 4* show that if $h_1(0)=g_1(0)$ then

$$\frac{\partial}{\partial x}[\mathcal{G}_{1,1}^{+}(x,t)+\mathcal{G}_{1,2}^{+}(x,t)+\mathcal{H}_{1,1}^{+}(x,t)]\approx[-ig_2'(0)-ih_1'(0)]\int_{1}^{\infty} e^{i\lambda x-\omega_1(\lambda)t}\frac{d\lambda}{\rho(\lambda)},$$

$$\frac{\partial}{\partial x}[\mathcal{G}_{1,1}^{-}(x,t)+\mathcal{G}_{1,2}^{-}(x,t)+\mathcal{H}_{1,1}^{-}(x,t)]\approx 0,$$

$$\frac{\partial}{\partial x}[-\mathcal{G}_{1,3}^{+}(x,t)-\mathcal{G}_{1,4}^{+}(x,t)-\mathcal{H}_{1,2}^{+}(x,t)]\approx 0,$$

$$\frac{\partial}{\partial x}[\mathcal{G}_{1,3}^{-}(x,t)+\mathcal{G}_{1,4}^{-}(x,t)+\mathcal{H}_{1,2}^{-}(x,t)]\approx[-ig_2'(0)-ih_1'(0)]\int_{-\infty}^{-1} e^{i\lambda x-\omega_2(\lambda)t}\frac{d\lambda}{\rho(\lambda)}.$$

It follows that

(4.22) $$g_1(0)=h_1(0)\ \&\ h_1'(0)=-g_2'(0)\ \Rightarrow \lim_{\substack{(x,t)\to(a,a)\\ x\neq t}}\frac{\partial\psi_1(x,t)}{\partial x}\ \text{exists for every}\ a>0.$$

Similarly, for every $a>0$,

(4.23) $$g_1(0)=h_1(0)\ \&\ h_1'(0)=-g_2'(0)$$
$$\Rightarrow \lim_{\substack{(x,t)\to(a,a)\\ x\neq t}}\frac{\partial\psi_1(x,t)}{\partial t},\ \lim_{\substack{(x,t)\to(a,a)\\ x\neq t}}\frac{\partial\psi_2(x,t)}{\partial x},\ \lim_{\substack{(x,t)\to(a,a)\\ x\neq t}}\frac{\partial\psi_2(x,t)}{\partial t}\ \text{exist.}$$

**Step 6** *Completion of the proof of the theorem*

**(1)** The 1st assertion follows from the analysis of the integrals carried out in section 3. In this regard, we point out that, although these integrals do not converge absolutely, their equivalent ones, obtained by appropriately deforming their contours, as this was done in section 3, do converge absolutely.

**(2)** It follows from (3.5) – (3.10) that $\mathcal{G}_{1,1}(x,t)\in C^{\infty}(Q-\{x=t\})$. For example, (3.7) implies that for $J_1(x,t)$, which is part of $\mathcal{G}_{1,1}^{1,+}(x,t)$, we have, for nonnegative integers $n$ and $l$,

(4.24) $$\frac{\partial^{n,l}J_1(x,t)}{\partial x^n\partial t^l}=\sum_{k=1}^{N}g_2^{(k)}(0)\int_{\Gamma^{+,+}}(i\lambda)^n[-\omega_1(\lambda)]^l e^{i\lambda x-\omega_1(\lambda)t}\frac{1}{(i\lambda)^k}\frac{\lambda d\lambda}{2\rho(\lambda)}$$

$$+\int_{1}^{\infty}(i\lambda)^n[-\omega_1(\lambda)]^l e^{i\lambda x-\omega_1(\lambda)t}\frac{(g_2^{(N)})\hat{}(\lambda)}{(i\lambda)^N}\frac{\lambda d\lambda}{2\rho(\lambda)}, \text{ for } x>t,$$

provided that $N\geq n+l+1$. This, in turn, implies that $J_1(x,t)\in C^\infty(Q\cap\{x>t\})$.
In view of the analysis of the integrals, carried out in section 3, we see that such computations can be made for all the integrals in (1.4) and (1.5), both for $x>t$ and $x<t$, proving the 2nd assertion of the theorem.

**(3)** The 3rd assertion was proved in *Step 1*.

**(4)** The 4th and 5th assertions were proved in *Step 2* and *Step 3*, respectively. The part of the assertions that these limits are uniform in $t$ or $x$, in compact sets, follow easily by examining their proofs.

**(5)** Finally, the 6th and 7th conclusions of the theorem follow from what we proved in *Step 4* and *Step 5*, respectively.

We point out that, besides the existence of the limits of *Step 4*, we have to show also that the extended functions $\psi_j(x,t)$ are continuous on $Q$, as far as the 6th assertion is concerned. This follows easily from the integral representations of these extensions, once the integrals

$$\int_{-\infty}^{-1}e^{i\lambda x-\omega_1(\lambda)t}\frac{d\lambda}{\rho(\lambda)},\ \int_{-\infty}^{-1}e^{i\lambda x-\omega_2(\lambda)t}\frac{d\lambda}{\rho(\lambda)},\ \int_{1}^{\infty}e^{i\lambda x-\omega_1(\lambda)t}\frac{d\lambda}{\rho(\lambda)} \text{ and } \int_{1}^{\infty}e^{i\lambda x-\omega_2(\lambda)t}\frac{d\lambda}{\rho(\lambda)}$$

which are the obstacles for the existence of the limits, do not appear in the formulas, having annihilate their coefficient, namely $[-ig_1(0)+ih_1(0)]$.

An analogous verification has to be made, in order to complete also the proof of the 7th assertion, i.e., to prove that the extended functions $\psi_j(x,t)$ are $C^1$ on $Q$, under the assumption $g_1(0)=h_1(0)$ & $h_1'(0)=-g_2'(0)$. (In this regard, see also the 4th part of Theorem 4, in Section 7, where we express the functions $\Psi_j(x,t)$, $j=1,2$, in a “unified formula”, for $(x,t)$ in $Q$, i.e., including the diagonal « $x=t$ », in the more general case of the inhomogeneous system.)

## 5. Boundary behavior of the derivatives of the solution functions

**Theorem 2** *Assuming (1.3) and letting $\psi_1(x,t)$ and $\psi_2(x,t)$ be the functions defined by (1.4) and (1.5), the following hold:*

*1st For nonnegative integers $n$ and $l$, the limits*

$$\lim_{t\to0^+}\frac{\partial^{n,l}\psi_1(x,t)}{\partial x^n\partial t^l} \text{ and } \lim_{t\to0^+}\frac{\partial^{n,l}\psi_2(x,t)}{\partial x^n\partial t^l}$$

*exist, uniformly for $x$ in compact subsets of $\mathbb{R}^+$ and define $C^\infty$ functions of $x\in\mathbb{R}^+$.*

*2nd For nonnegative integers $n$ and $l$, the limits*

$$\lim_{x\to0^+}\frac{\partial^{n,l}\psi_1(x,t)}{\partial x^n\partial t^l} \text{ and } \lim_{x\to0^+}\frac{\partial^{n,l}\psi_2(x,t)}{\partial x^n\partial t^l}$$

*exist, uniformly for $t$ in compact subsets of $\mathbb{R}^+$ and define $C^\infty$ functions of $t\in\mathbb{R}^+$.*

*3rd The functions $\psi_1(x,t)$ and $\psi_2(x,t)$ extend to $C^\infty$ functions in open neighborhoods of the open half-lines $\{(x,0):x>0\}$ and $\{(0,t):t>0\}$.*

$4^{th}$ $\lim_{t\to 0^+}\frac{\partial^n \psi_1(x,t)}{\partial x^n}=\frac{d^n g_1(x)}{dx^n}$ *and* $\lim_{t\to 0^+}\frac{\partial^n \psi_2(x,t)}{\partial x^n}=\frac{d^n g_2(x)}{dx^n}$, *for* $x\in\mathbb{R}^+$.

$5^{th}$ $\lim_{x\to 0^+}\frac{\partial^l \psi_1(x,t)}{\partial t^l}=\frac{d^l h_1(t)}{dt^l}$, *for* $t\in\mathbb{R}^+$.

**Proof Step 1** It follows from (4.17) that, for fixed $n$ and $l$, and sufficiently large $N$,

$$\text{(5.1)}\quad \lim_{t\to 0^+}\frac{\partial^{n,l} J_1(x,t)}{\partial x^n \partial t^l}=\sum_{k=1}^{N} g_2^{(k)}(0)\int_{\Gamma^{+,+}}(i\lambda)^n[-\omega_1(\lambda)]^l e^{i\lambda x}\frac{1}{(i\lambda)^k}\frac{\lambda d\lambda}{2\rho(\lambda)}$$
$$+\int_1^{\infty}(i\lambda)^n[-\omega_1(\lambda)]^l e^{i\lambda x}\frac{(g_2^{(N)}\hat{)}(\lambda)}{(i\lambda)^N}\frac{\lambda d\lambda}{2\rho(\lambda)},\ \text{for } x>0.$$

It follows that the limit in (5.1) exists, uniformly for $x$ in compact sets of $\mathbb{R}^+$ and defines a $C^\infty$ function of $x\in\mathbb{R}^+$.
In view of the analysis of the integrals, carried out in section 3, we see that such computations can be made for all the integrals in (1.4) and (1.5), for $x>0$, proving the 1st assertion of the theorem.

**Step 2** It follows from (3.9) that, for fixed $n$, $l$, and sufficiently large $N$,

$$\text{(5.2)}\quad \lim_{x\to 0^+}\frac{\partial^{n,l}\mathcal{G}_{1,1}^{1,+}(x,t)}{\partial x^n \partial t^l}(x,t)=\int_{\Gamma^{+,-}}(i\lambda)^n[-\omega_1(\lambda)]^l e^{-\omega_1(\lambda)t}\{\lambda\hat{g}_2(\lambda)+[\rho(\lambda)+m]\hat{g}_1(\lambda)\}\frac{d\lambda}{2\rho(\lambda)}.$$

It follows that the limit in (5.2) exists, uniformly for $t$ in compact sets of $\mathbb{R}^+$, and defines a $C^\infty$ function of $t\in\mathbb{R}^+$.
In view of the analysis of the integrals, carried out in section 3, we see that such computations can be made for all the $\mathcal{G}-integrals$ in (1.4) and (1.5) (working with $x<t$), proving the 2nd assertion of the theorem, as far as the $\mathcal{G}-parts$ of $\psi_1(x,t)$ and $\psi_2(x,t)$ are concerns.
As for the $\mathcal{H}-integrals$, we have to pay attention to a further difficulty which arises when we let $x\to 0^+$ in the derivatives $\frac{\partial^{n,l}\mathcal{H}_{j,s}(x,t)}{\partial x^n \partial t^l}$, $j,s\in\{1,2\}$.
More precisely, according to the decomposition (3.3), in the difference $\mathcal{H}_{1,1}(x,t)-\mathcal{H}_{1,2}(x,t)$, there appears the quantity

$$\mathfrak{D}_1(x,t):=\int_\Delta e^{i\lambda x}\left[\sum_{k=1}^{N}(-1)^{k-1}\frac{h_1^{(k-1)}(t)}{[\omega_1(\lambda)]^k}\right]\frac{\lambda d\lambda}{\rho(\lambda)}-\int_\Delta e^{i\lambda x}\left[\sum_{k=1}^{N}(-1)^{k-1}\frac{h_1^{(k-1)}(t)}{[\omega_2(\lambda)]^k}\right]\frac{\lambda d\lambda}{\rho(\lambda)}.$$

The point here is that, as long as $x$ is kept positive, we can differentiate $\mathfrak{D}_1(x,t)$, obtaining

$$\text{(5.3)}\quad \frac{\partial^{n,l}\mathfrak{D}_1(x,t)}{\partial x^n \partial t^l}=\int_\Delta (i\lambda)^n e^{i\lambda x}\left[\sum_{k=1}^{N}(-1)^{k-1}\frac{h_1^{(k+l-1)}(t)}{[\omega_1(\lambda)]^k}\right]\frac{\lambda d\lambda}{\rho(\lambda)}$$
$$-\int_\Delta (i\lambda)^n e^{i\lambda x}\left[\sum_{k=1}^{N}(-1)^{k-1}\frac{h_1^{(k+l-1)}(t)}{[\omega_2(\lambda)]^k}\right]\frac{\lambda d\lambda}{\rho(\lambda)},\ \text{for } x>0 \text{ and } t>0.$$

However, if we let $x\to 0^+$, we are led to integrals of the form

$$\int_{\Delta}(i\lambda)^{n}\left[\sum_{k=1}^{N}(-1)^{k-1}\frac{h_1^{(k+l-1)}(t)}{[\omega_j(\lambda)]^k}\right]\frac{\lambda d\lambda}{\rho(\lambda)},\ j=1,2,$$

which are highly oscillatory.

In order to overcome this difficulty, we consider the differences

$$\mathfrak{D}_{1,k}(x,t):=\int_{\Delta}e^{i\lambda x}\frac{h_1^{(k-1)}(t)}{[\omega_1(\lambda)]^k}\frac{\lambda d\lambda}{\rho(\lambda)}-\int_{\Delta}e^{i\lambda x}\frac{h_1^{(k-1)}(t)}{[\omega_2(\lambda)]^k}\frac{\lambda d\lambda}{\rho(\lambda)},\ 1\le k\le n,$$

so that

$$\mathfrak{D}_1(x,t)=\sum_{k=1}^{N}(-1)^{k-1}\mathfrak{D}_{1,k}(x,t).$$

Then, in the case $k$ is odd, we have

$$\mathfrak{D}_{1,k}(x,t)=2\int_{\Delta}e^{i\lambda x}\frac{h_1^{(k-1)}(t)}{[\omega_1(\lambda)]^k}\frac{\lambda d\lambda}{\rho(\lambda)}=2h_1^{(k-1)}(t)\int_{C(mi,\varepsilon)}e^{i\lambda x}\frac{\lambda d\lambda}{i^k(\lambda^2+m^2)^{(k+1)/2}},$$

which implies that

$$\frac{\partial^{n,l}\mathfrak{D}_{1,k}(x,t)}{\partial x^n\partial t^l}=2h_1^{(k+l-1)}(t)\int_{C(mi,\varepsilon)}(i\lambda)^n e^{i\lambda x}\frac{\lambda d\lambda}{i^k(\lambda^2+m^2)^{(k+1)/2}}\lambda d\lambda$$

($C(mi,\varepsilon):=\{\lambda\in\mathbb{C}:|\lambda-mi|=\varepsilon\}$) and, therefore,

$$\text{(5.4)}\qquad \lim_{x\to 0^+}\frac{\partial^{n,l}\mathfrak{D}_{1,k}(x,t)}{\partial x^n\partial t^l}=2h_1^{(k+l-1)}(t)\int_{C(mi,\varepsilon)}(i\lambda)^n\frac{\lambda d\lambda}{i^k(\lambda^2+m^2)^{(k+1)/2}}\ \text{(for odd } k\text{)}.$$

In the case $k$ is even, the corresponding limits, as $x\to 0^+$, vanish trivially, since then $\mathfrak{D}_{1,k}(x,t)\equiv 0$.

Similarly, according to the decomposition (3.3), in the difference $\mathcal{H}_{2,1}(x,t)-\mathcal{H}_{2,2}(x,t)$, there appears the quantity

$$\text{(5.5)}\quad \mathfrak{D}_2(x,t):=\int_{\Delta}e^{i\lambda x}\left[\sum_{k=1}^{N}(-1)^{k-1}\frac{h_1^{(k-1)}(t)}{[\omega_1(\lambda)]^k}\right]\frac{[i\omega_1(\lambda)+m]d\lambda}{\rho(\lambda)} -\int_{\Delta}e^{i\lambda x}\left[\sum_{k=1}^{N}(-1)^{k-1}\frac{h_1^{(k-1)}(t)}{[\omega_2(\lambda)]^k}\right]\frac{[i\omega_2(\lambda)+m]d\lambda}{\rho(\lambda)}.$$

If $k$ is even then

$$\begin{aligned}\text{(5.6)}\quad \mathfrak{D}_{2,k}(x,t)&:=\int_{\lambda\in\Delta}e^{i\lambda x}\frac{h_1^{(k-1)}(t)}{[\omega_1(\lambda)]^k}\frac{[i\omega_1(\lambda)+m]d\lambda}{\rho(\lambda)}-\int_{\lambda\in\Delta}e^{i\lambda x}\frac{h_1^{(k-1)}(t)}{[\omega_2(\lambda)]^k}\frac{[i\omega_2(\lambda)+m]d\lambda}{\rho(\lambda)}\\
&=\int_{\Delta}e^{i\lambda x}\frac{h_1^{(k-1)}(t)}{[\omega_1(\lambda)]^k}\frac{i\omega_1(\lambda)d\lambda}{\rho(\lambda)}-\int_{\Delta}e^{i\lambda x}\frac{h_1^{(k-1)}(t)}{[\omega_1(\lambda)]^k}\frac{i\omega_2(\lambda)d\lambda}{\rho(\lambda)}\\
&=2\int_{\lambda\in\Delta}e^{i\lambda x}\frac{h_1^{(k-1)}(t)}{[\omega_1(\lambda)]^k}\frac{i\omega_1(\lambda)d\lambda}{\rho(\lambda)}=-2h_1^{(k-1)}(t)\int_{\lambda\in\Delta}e^{i\lambda x}\frac{d\lambda}{i^k(\lambda^2+m^2)^{k/2}}\\
&=-2h_1^{(k-1)}(t)\int_{C(mi,\varepsilon)}e^{i\lambda x}\frac{d\lambda}{i^k(\lambda^2+m^2)^{k/2}}.\end{aligned}$$

If $k$ is odd then

(5.7) $\mathfrak{D}_{2,k}(x,t) = 2mh_1^{(k-1)}(t)\int\limits_{\Delta} e^{i\lambda x}\dfrac{d\lambda}{[\omega_1(\lambda)]^k\rho(\lambda)}$

$$= 2mh_1^{(k-1)}(t)\int\limits_{\Delta} e^{i\lambda x}\frac{d\lambda}{i^k(\lambda^2+m^2)^{(k+1)/2}} = 2mh_1^{(k-1)}(t)\int\limits_{C(mi,\varepsilon)} e^{i\lambda x}\frac{d\lambda}{i^k(\lambda^2+m^2)^{(k+1)/2}}.$$

Now the 2nd assertion of the theorem, for the $\mathcal{H}-parts$ of $\psi_1(x,t)$ and $\psi_2(x,t)$, follows from (5.3), (5.4), (5.5), (5.6) and (5.7), and completes the proof of 2nd assertion.

**Step 3** *Completion of the proof of the theorem* The 3rd conclusion follows from the 1st and 2nd parts. The 4th conclusion follows from the 3rd one, combined with the fact that

$$\psi_1(x,t)\big|_{(x,t)=(x,0)} = g_1(x) \text{ and } \psi_2(x,t)\big|_{(x,t)=(x,0)} = g_2(x),\ x>0.$$

Finally, the 5th conclusion of the theorem follows from the 3rd one, combined with the fact that

$$\psi_1(x,t)\big|_{(x,t)=(0,t)} = h_1(t),\ t>0.$$

## 6. Behavior of the solution functions as $x \to +\infty$

**Theorem 3** *Assuming (1.3) and letting $\psi_1(x,t)$ and $\psi_2(x,t)$ be the functions defined by (1.4) and (1.5), the following hold: For $T>0$,*

$$\lim_{x\to+\infty}\left[x^E\frac{\partial^{n+l}\psi_1(x,t)}{\partial x^n\partial t^l}\right]=0 \text{ and } \lim_{x\to+\infty}\left[x^E\frac{\partial^{n+l}\psi_2(x,t)}{\partial x^n\partial t^l}\right]=0,$$

*for nonnegative integers $n$, $l$ and $E$, uniformly for $0<t\le T$.*

**Proof** In this proof, we use the formulas of section 3 for $x>T\ge t$.

**Step 1** Let

$$\mathcal{I}(x,t) := \int\limits_{\lambda=-\infty}^{\infty} e^{i\lambda x-\omega_1(\lambda)t}\{\lambda\hat{g}_2(\lambda)\}\frac{d\lambda}{2\rho(\lambda)},$$

which is part of the integral $\mathcal{G}_{1,1}(x,t)$ (the first integral in the RHS of (1.4)).

We claim that

(6.1) $$\lim_{x\to+\infty}\frac{\partial^{n+l}\mathcal{I}(x,t)}{\partial x^n\partial t^l}=0.$$

To prove this we use (3.1) and we write $\mathcal{I}(x,t)$ as follows:

(6.2) $$\mathcal{I}(x,t) = \int\limits_{-\infty}^{-1} e^{i\lambda x-\omega_1(\lambda)t}\left[\lambda\frac{(g_2^{(N)})\hat{}(\lambda)}{(i\lambda)^{\mathrm{N}}}\right]\frac{d\lambda}{2\rho(\lambda)} + \int\limits_{\Gamma^{-,+}} e^{i\lambda x-\omega_1(\lambda)t}\left[\lambda\sum_{k=1}^{N}\frac{g_2^{(k-1)}(0)}{(i\lambda)^k}\right]\frac{d\lambda}{2\rho(\lambda)}$$
$$+\int\limits_{-1}^{1} e^{i\lambda x-\omega_1(\lambda)t}\{\lambda\hat{g}_2(\lambda)\}\frac{d\lambda}{2\rho(\lambda)} + \int\limits_{\Gamma^{+,+}} e^{i\lambda x-\omega_1(\lambda)t}\left[\lambda\sum_{k=1}^{N}\frac{g_2^{(k-1)}(0)}{(i\lambda)^k}\right]\frac{d\lambda}{2\rho(\lambda)}$$
$$+\int\limits_{1}^{\infty} e^{i\lambda x-\omega_1(\lambda)t}\left[\lambda\frac{(g_2^{(N)})\hat{}(\lambda)}{(i\lambda)^N}\right]\frac{d\lambda}{2\rho(\lambda)} =: \mathcal{I}_1+\mathcal{I}_2+\mathcal{I}_3+\mathcal{I}_4+\mathcal{I}_5.$$

By the Riemann-Lebesgue lemma, given $n$ and $l$, and choosing sufficiently large $N$ in (6.2), we obtain

$$\lim_{x\to+\infty}\frac{\partial^{n+l}\mathcal{I}_1(x,t)}{\partial x^n\partial t^l}=0,\quad \lim_{x\to+\infty}\frac{\partial^{n+l}\mathcal{I}_3(x,t)}{\partial x^n\partial t^l}=0,\quad \lim_{x\to+\infty}\frac{\partial^{n+l}\mathcal{I}_5(x,t)}{\partial x^n\partial t^l}=0. \tag{6.3}$$

Moreover, the convergence in (6.3) is uniform for $0<t\le T$. (In this regard, let us notice that these integrals are, uniformly for $0<t\le T$, absolutely convergent.)
On the other hand,

$$\frac{\partial^{n+l}\mathcal{I}_4(x,t)}{\partial x^n\partial t^l}=\int_{\Gamma^{+,+}}(i\lambda)^n[-\omega_1(\lambda)]^l e^{i\lambda x-\omega_1(\lambda)t}\left[\lambda\sum_{k=1}^{N}\frac{g_2^{(k-1)}(0)}{(i\lambda)^k}\right]\frac{d\lambda}{2\rho(\lambda)},$$

which implies that, for $x>2T$,

$$\left|\frac{\partial^{n+l}\mathcal{I}_4(x,t)}{\partial x^n\partial t^l}\right|\preceq\int_{[1,\,1+i]}\left|\lambda^n[\omega_1(\lambda)]^l e^{i\lambda x-\omega_1(\lambda)t}\right|\left|\sum_{k=1}^{N}\frac{g_2^{(k-1)}(0)}{(i\lambda)^k}\right|d|\lambda|$$
$$+e^{-x/4}\int_{[1+i,\,+\infty e^{1\pi/4}]}\left|\lambda^n[\omega_1(\lambda)]^l e^{i\lambda x-\omega_1(\lambda)t}\right|\left|\sum_{k=1}^{N}\frac{g_2^{(k-1)}(0)}{(i\lambda)^k}\right|d|\lambda|.$$

It follows that

$$\lim_{x\to+\infty}\frac{\partial^{n+l}\mathcal{I}_4(x,t)}{\partial x^n\partial t^l}=0,\text{ uniformly for }0<t\le T. \tag{6.4}$$

Similarly we show that

$$\lim_{x\to+\infty}\frac{\partial^{n+l}\mathcal{I}_2(x,t)}{\partial x^n\partial t^l}=0,\text{ uniformly for }0<t\le T. \tag{6.5}$$

Thus, the claim follows from (6.2), (6.3), (6.4) and (6.5).

**Step 2** We claim that

$$\lim_{x\to+\infty}\left[x\frac{\partial^{n+l}\mathcal{I}(x,t)}{\partial x^n\partial t^l}\right]=0,\text{ uniformly for }0<t\le T. \tag{6.6}$$

To prove this, we differentiate (6.2), with sufficiently large $N$, obtaining

$$\frac{\partial^{n+l}\mathcal{I}(x,t)}{\partial x^n\partial t^l}=\int_{-\infty}^{-1}(i\lambda)^n[-\omega_1(\lambda)]^l e^{i\lambda x-\omega_1(\lambda)t}\left[\lambda\frac{(g_2^{(N)})\hat{}(\lambda)}{(i\lambda)^N}\right]\frac{d\lambda}{2\rho(\lambda)} \tag{6.7}$$
$$+\int_{\Gamma^{-,+}}(i\lambda)^n[-\omega_1(\lambda)]^l e^{i\lambda x-\omega_1(\lambda)t}\left[\lambda\sum_{k=1}^{N}\frac{g_2^{(k-1)}(0)}{(i\lambda)^k}\right]\frac{d\lambda}{2\rho(\lambda)}$$
$$+\int_{-1}^{1}(i\lambda)^n[-\omega_1(\lambda)]^l e^{i\lambda x-\omega_1(\lambda)t}\{\lambda\hat{g}_2(\lambda)\}\frac{d\lambda}{2\rho(\lambda)}$$
$$+\int_{\Gamma^{+,+}}(i\lambda)^n[-\omega_1(\lambda)]^l e^{i\lambda x-\omega_1(\lambda)t}\left[\lambda\sum_{k=1}^{N}\frac{g_2^{(k-1)}(0)}{(i\lambda)^k}\right]\frac{d\lambda}{2\rho(\lambda)}$$
$$+\int_{1}^{\infty}(i\lambda)^n[-\omega_1(\lambda)]^l e^{i\lambda x-\omega_1(\lambda)t}\left[\lambda\frac{(g_2^{(N)})\hat{}(\lambda)}{(i\lambda)^{\mathrm{N}}}\right]\frac{d\lambda}{2\rho(\lambda)}.$$

Multiplying (6.7) by $ix$ and integrating by parts, we obtain

$$(6.8)\quad ix\frac{\partial^{n+l}\mathcal{I}(x,t)}{\partial x^n\partial t^l}=\int_{-\infty}^{-1}e^{i\lambda x}\frac{d}{d\lambda}\left\{(i\lambda)^n[-\omega_1(\lambda)]^l e^{-\omega_1(\lambda)t}\left[\lambda\frac{(g_2^{(N)})\hat{}(\lambda)}{(i\lambda)^N}\right]\frac{1}{2\rho(\lambda)}\right\}d\lambda$$
$$+\int_{\Gamma^{-,+}}e^{i\lambda x}\frac{d}{d\lambda}\left\{(i\lambda)^n[-\omega_1(\lambda)]^l e^{-\omega_1(\lambda)t}\left[\lambda\sum_{k=1}^{N}\frac{g_2^{(k-1)}(0)}{(i\lambda)^k}\right]\frac{1}{2\rho(\lambda)}\right\}d\lambda$$
$$+\int_{-1}^{1}e^{i\lambda x}\frac{d}{d\lambda}\left\{(i\lambda)^n[-\omega_1(\lambda)]^l e^{-\omega_1(\lambda)t}\{\lambda\hat{g}_2(\lambda)\}\frac{1}{2\rho(\lambda)}\right\}d\lambda$$
$$+\int_{\Gamma^{+,+}}e^{i\lambda x}\frac{d}{d\lambda}\left\{(i\lambda)^n[-\omega_1(\lambda)]^l e^{-\omega_1(\lambda)t}\left[\lambda\sum_{k=1}^{N}\frac{g_2^{(k-1)}(0)}{(i\lambda)^k}\right]\right\}d\lambda$$
$$+\int_{1}^{\infty}e^{i\lambda x}\frac{d}{d\lambda}\left\{(i\lambda)^n[-\omega_1(\lambda)]^l e^{-\omega_1(\lambda)t}\left[\lambda\frac{(g_2^{(N)})\hat{}(\lambda)}{(i\lambda)^{\mathrm{N}}}\right]\right\}d\lambda\,.$$

(We point out that the intermediate boundary terms, in the integration by parts processes, cancel each other.)
Now, using (6.8) and working as in the proof of (6.1), we can prove (6.6).

**Step 3** Working as in step2 and integrating by parts, repeatedly, we show that

$$(6.9)\qquad \lim_{x\to+\infty}\left[x^E\frac{\partial^{n+l}\mathcal{I}(x,t)}{\partial x^n\partial t^l}\right]=0\text{, uniformly for } 0<t\le T\,.$$

**Step 3** *Completion of the proof of the theorem* Working as in the previous steps, we can show that if $\mathcal{J}(x,t)$ is any of the integrals in the RHSs of (1.4) or (1.5), then

$$\lim_{x\to+\infty}\left[x^E\frac{\partial^{n+l}\mathcal{J}(x,t)}{\partial x^n\partial t^l}\right]=0\text{, uniformly for } 0<t\le T\,,$$

and this completes the proof of the theorem.

## 7. The inhomogeneous system

In this section we address the inhomogeneous version of the problem.

***Problem*** *Solve the system of differential equations*

$$(7.1)\qquad \begin{cases}\partial_t\Psi_1=-\partial_x\Psi_2-im\Psi_1+f_1\\ \partial_t\Psi_2=-\partial_x\Psi_1+im\Psi_2+f_2,\end{cases}$$

*for* $\Psi_1(x,t)$ *and* $\Psi_2(x,t)$, *with* $(x,t)\in\mathbb{R}^+\times\mathbb{R}^+$, *subject to the following initial and boundary conditions:*

$$(7.2)\qquad \begin{cases}\lim\limits_{t\to0^+}\Psi_1(x,t)=g_1(x),\ \ x\in\mathbb{R}^+,\\ \lim\limits_{t\to0^+}\Psi_2(x,t)=g_2(x),\ \ x\in\mathbb{R}^+,\\ \lim\limits_{x\to0^+}\Psi_1(x,t)=h_1(t),\ \ t\in\mathbb{R}^+.\end{cases}$$

***Assumptions*** Throughtout this section, besides (1.3) we make also the following assumptions for the functions $f_1(x,t)$ and $f_2(x,t)$:

(7.3) $\quad f_j(x,t)\in C^\infty(\overline{Q})$ , where $Q=\mathbb{R}^+\times\mathbb{R}^+$, such that $\dfrac{\partial^l f_j}{\partial t^l}(\cdot,t)\in\mathcal{S}([0,\infty))$ $(\forall l)$ $(j=1,2)$,

uniformly for $t$ in compact subsets of $\mathbb{R}^+$. The latter condition imposed on $f_j$ means that, for every $T>0$,

$$\sup\left\{x^E\left|\frac{\partial^{n+l}f_j(x,t)}{\partial x^n\partial t^l}\right|:\ x>0,\ 0\le t\le T\right\}<\infty,\text{ for all nonnegative integers } n,\ l \text{ and } E.$$

***Notation*** For $\lambda\in\mathbb{C}$ with $\operatorname{Im}\lambda\le 0$, we define

$$\hat{f}(\lambda,t)=\int_{y=0}^{\infty}e^{-i\lambda y}f(y,t)dy \text{ and } \tilde{\hat{f}}(\lambda,\omega,t)=\int_{\tau=0}^{t}e^{\omega(\lambda)\tau}\hat{f}(\lambda,\omega,\tau)d\tau \quad (f\in\{f_1,f_2\},\ \omega\in\{\omega_1,\omega_2\}).$$

*The $\mathcal{F}$ – integrals* In order to write down the UTM solution of problem (7.1) & (7.2), we will also use the following notation:

$$\begin{aligned}\mathcal{F}_{1,1}(x,t):=&\int_{-\infty}^{\infty}e^{i\lambda x-\omega_1(\lambda)t}\{\lambda\tilde{\hat{f}}_2(\lambda,\omega_1,t)+[i\omega_2(\lambda)+m]\tilde{\hat{f}}_1(\lambda,\omega_1,t)\}\frac{d\lambda}{2\rho(\lambda)}\\&+\int_{-\infty}^{\infty}e^{i\lambda x-\omega_1(\lambda)t}\{\lambda\tilde{\hat{f}}_2(-\lambda,\omega_1,t)-[i\omega_2(\lambda)+m]\tilde{\hat{f}}_1(-\lambda,\omega_1,t)\}\frac{d\lambda}{2\rho(\lambda)}\\&=:(\mathcal{F}_{1,1})_1(x,t)+(\mathcal{F}_{1,1})_2(x,t),\end{aligned}$$

$$\begin{aligned}\mathcal{F}_{1,2}(x,t):=&-\int_{-\infty}^{\infty}e^{i\lambda x-\omega_2(\lambda)t}\{\lambda\tilde{\hat{f}}_2(\lambda,\omega_2,t)+[i\omega_1(\lambda)+m]\tilde{\hat{f}}_1(\lambda,\omega_2,t)\}\frac{d\lambda}{2\rho(\lambda)}\\&-\int_{-\infty}^{\infty}e^{i\lambda x-\omega_2(\lambda)t}\{\lambda\tilde{\hat{f}}_2(-\lambda,\omega_2,t)-[i\omega_1(\lambda)+m]\tilde{\hat{f}}_1(-\lambda,\omega_2,t)\}\frac{d\lambda}{2\rho(\lambda)}\\&=:-(\mathcal{F}_{1,2})_1(x,t)-(\mathcal{F}_{1,2})_2(x,t),\end{aligned}$$

$$\begin{aligned}\mathcal{F}_{2,1}(x,t):=&\int_{-\infty}^{\infty}e^{i\lambda x-\omega_1(\lambda)t}\{\lambda\tilde{\hat{f}}_1(\lambda,\omega_1,t)-[i\omega_1(\lambda)+m]\tilde{\hat{f}}_2(\lambda,\omega_1,t)\}\frac{d\lambda}{2\rho(\lambda)}\\&-\int_{-\infty}^{\infty}e^{i\lambda x-\omega_1(\lambda)t}\{\lambda\tilde{\hat{f}}_1(-\lambda,\omega_1,t)+[i\omega_1(\lambda)+m]\tilde{\hat{f}}_2(-\lambda,\omega_1,t)\}\frac{d\lambda}{2\rho(\lambda)}\\&=:(\mathcal{F}_{2,1})_1(x,t)-(\mathcal{F}_{2,1})_2(x,t),\end{aligned}$$

$$\begin{aligned}\mathcal{F}_{2,2}(x,t):=&-\int_{-\infty}^{\infty}e^{i\lambda x-\omega_2(\lambda)t}\{\lambda\tilde{\hat{f}}_1(\lambda,\omega_2,t)-[i\omega_2(\lambda)+m]\tilde{\hat{f}}_2(\lambda,\omega_2,t)\}\frac{d\lambda}{2\rho(\lambda)}\\&+\int_{-\infty}^{\infty}e^{i\lambda x-\omega_2(\lambda)t}\{\lambda\tilde{\hat{f}}_1(-\lambda,\omega_2,t)+[i\omega_2(\lambda)+m]\tilde{\hat{f}}_2(-\lambda,\omega_2,t)\}\frac{d\lambda}{2\rho(\lambda)}\\&=:-(\mathcal{F}_{2,2})_1(x,t)+(\mathcal{F}_{2,2})_2(x,t).\end{aligned}$$

**Theorem 4** *Assuming (1.3) and (7.3), we let $\psi_1(x,t)$ and $\psi_2(x,t)$ be as in (1.4) and (1.5), and we define, for $x>0$, $t>0$, $x\neq t$,*

$$\Psi_1(x,t)=\psi_1(x,t)+\frac{1}{2\pi}\mathcal{F}_{1,1}(x,t)+\frac{1}{2\pi}\mathcal{F}_{1,2}(x,t) \quad \text{and} \quad \Psi_2(x,t)=\psi_2(x,t)+\frac{1}{2\pi}\mathcal{F}_{2,1}(x,t)+\frac{1}{2\pi}\mathcal{F}_{2,2}(x,t).$$

*Then the following hold:*

*1[st] The $\mathcal{F}-$integrals, i.e., the ones appearing in the definition of $\mathcal{F}_{1,1}(x,t)$, $\mathcal{F}_{1,2}(x,t)$, $\mathcal{F}_{2,1}(x,t)$ and $\mathcal{F}_{2,2}(x,t)$, converge absolutely, for $x\geq 0$, $t\geq 0$, and define continuous functions for $(x,t)\in\overline{Q}$.*

*2[nd] The functions $\Psi_1(x,t)$ and $\Psi_2(x,t)$ satisfy the analogues of the 2[nd], 4[th], 5[th] and 6[th] conclusions of Theorem 1, i.e., these conclusions hold with $\psi_1(x,t)$ and $\psi_2(x,t)$ replaced by $\Psi_1(x,t)$ and $\Psi_2(x,t)$, respectively.*

*3[rd] The functions $\Psi_1(x,t)$ and $\Psi_2(x,t)$ satisfy the differential equations (7.1) in $Q-\{x=t\}$.*

*4[th] If, in addition to (1.3) and (7.3), we assume*

$$h_1(0)=g_1(0) \quad \text{and} \quad h_1'(0)=-g_2'(0)+f_1(0,0), \tag{7.4}$$

*then the functions $\Psi_1(x,t)$ and $\Psi_2(x,t)$ extend to $C^1$ functions on $Q$ and satisfy, there, the differential equations (7.1). These $C^1$ extensions of $\Psi_1(x,t)$ and $\Psi_2(x,t)$ for $(x,y)\in Q$, are given by the following unified formulas:*

$$(7.5)\quad 2\pi\Psi_1(x,t)=\int_{-\infty}^{\infty}\Big\{\mathcal{I}(\mathcal{G}_{1,1})+\mathcal{I}(\mathcal{G}_{1,2})+\mathcal{I}(\mathcal{H}_{1,1})-\mathcal{I}(\mathcal{G}_{1,3})-\mathcal{I}(\mathcal{G}_{1,4})-\mathcal{I}(\mathcal{H}_{1,2}) \\ +\mathcal{I}[([\mathcal{F}_{1,1})_1]+\mathcal{I}[(\mathcal{F}_{1,1})_2]-\mathcal{I}[(\mathcal{F}_{1,2})_1]-\mathcal{I}[(\mathcal{F}_{1,2})_2]\Big\}(\lambda,x,t)\,d\lambda$$

*and*

$$(7.6)\quad 2\pi\Psi_2(x,t)=\int_{-\infty}^{\infty}\Big\{\mathcal{I}(\mathcal{G}_{2,1})+\mathcal{I}(\mathcal{G}_{2,4})+\mathcal{I}(\mathcal{H}_{2,2})-\mathcal{I}(\mathcal{G}_{2,2})-\mathcal{I}(\mathcal{G}_{2,3})-\mathcal{I}(\mathcal{H}_{2,2}) \\ +\mathcal{I}[(\mathcal{F}_{2,1})_1]-\mathcal{I}[(\mathcal{F}_{2,1})_2]-\mathcal{I}[(\mathcal{F}_{2,2})_1]+\mathcal{I}[(\mathcal{F}_{2,2})_2]\Big\}(\lambda,x,t)\,d\lambda,$$

*where $\mathcal{I}(*)$ stands for the integrand of the indicated integral «$*$».*
*Also, the above integrals, expressing $\Psi_1(x,t)$ and $\Psi_2(x,t)$, converge absolutely and remain absolutely convergent, when we differentiate their integrands – once – with respect to $x$ or $t$, uniformly for $(x,t)$ in compact subsets of $Q$.*

For the proof of this theorem we need the following.

***Preliminary computations***

The formulas derived in this paragraph, are analogues of the formulas (3.1), (3.2) and (3.3), for the transforms $\hat{f}(\lambda,t)$ and $\tilde{\tilde{f}}(\lambda,\omega(\lambda),t)$.

**(1)** Generalizing (3.1), we see that

$$\hat{f}(\lambda,t)=P_{f,M}(\lambda,t)+\Upsilon_{f,M}(\lambda,t), \text{ for } \lambda\in\mathbb{C}-\{0\} \text{ with } \operatorname{Im}\lambda\leq 0, \tag{7.7}$$

where

$$P_{f,M}(\lambda,t)=\sum_{k=1}^{M}\frac{1}{(i\lambda)^k}\frac{\partial^{k-1}f(y,t)}{\partial y^{k-1}}\Bigg|_{y=0} \quad\text{and}\quad \Upsilon_{f,M}(\lambda,t)=\frac{1}{(i\lambda)^M}\int_{y=0}^{\infty}e^{-i\lambda y}\frac{\partial^M f(y,t)}{\partial y^M}dy .$$

It follows that for every $T>0$, uniformly for $t\le T$, and for $\lambda\to\infty$ with $\operatorname{Im}\lambda\le 0$,

$$\text{(7.8)}\qquad \hat{f}(\lambda,t)=\mathrm{O}(1/\lambda) \text{ and } \Upsilon_{f,M}(\lambda,t)=\hat{f}(\lambda,t)-P_{f,M}(\lambda,t)=\mathrm{O}(1/\lambda^{M+1}) .$$

Also,

$$\text{(7.9)}\quad \tilde{\hat{f}}(\lambda,\omega(\lambda),t)=\int_0^t e^{\omega(\lambda)\tau}\hat{f}(\lambda,\tau)d\tau$$

$$=\int_0^t e^{\omega(\lambda)\tau}P_{f,M}(\lambda,\tau)d\tau+\int_0^t e^{\omega(\lambda)\tau}\Upsilon_{f,M}(\lambda,\tau)d\tau=:\tilde{P}_{f,M}(\lambda,\omega(\lambda),t)+\tilde{\Upsilon}_{f,M}(\lambda,\omega(\lambda),t) .$$

Let us notice that, although (7.9) holds for $\lambda\in\mathbb{C}-\{0\}$ with $\operatorname{Im}\lambda\le 0$, the integral defining $\tilde{P}_{f,M}(\lambda,\omega(\lambda),t)$, is defined for every $\lambda\in\mathbb{C}-\{0\}$.

Further integration by parts gives

$$\text{(7.10)}\quad \tilde{P}_{f,M}(\lambda,\omega(\lambda),t)=e^{\omega(\lambda)t}\mu_{f,M,N}(\lambda,\omega(\lambda)t)-\mu_{f,M,N}(\lambda,\omega(\lambda),0)$$
$$-\frac{(-1)^{N-1}}{[\omega(\lambda)]^N}\int_0^t e^{\omega(\lambda)\tau}(P_{f,M})^{(N)}(\lambda,\tau)d\tau ,$$

where

$$\mu_{f,M,N}(\lambda,\omega(\lambda),t):=\frac{P_{f,M}(\lambda,t)}{\omega(\lambda)}-\frac{(P_{f,M})'(\lambda,t)}{\omega^2(\lambda)}+\cdots+(-1)^{N-1}\frac{(P_{f,M})^{(N-1)}(\lambda,t)}{\omega^N(\lambda)}\quad(\lambda\in\mathbb{C}-\{0\}).$$

(The derivatives of $P_{f,M}(\lambda,t)$, in the above quantity, are taken with respect to $t$.)

**(2)** Similarly to (3.4), we have, for $\lambda\to\infty$ with $\lambda\in\mathbb{C}$ and $\operatorname{Re}\omega(\lambda)\ge 0$,

$$\text{(7.11)}\qquad e^{-\omega(\lambda)t}\tilde{P}_{f,M}(\lambda,\omega(\lambda),t)=\mathrm{O}(1/\lambda^2) \text{ and } e^{-\omega(\lambda)t}\tilde{\Upsilon}_{f,M}(\lambda,\omega(\lambda),t)=\mathrm{O}(1/\lambda^{M+2})$$

and

$$\text{(7.12)}\qquad \frac{1}{[\omega(\lambda)]^N}e^{-\omega(\lambda)t}\int_0^t e^{\omega(\lambda)\tau}(P_{f,M})^{(N)}(\lambda,\tau)d\tau=\mathrm{O}(1/\lambda^{N+2}),$$

uniformly for $t$ in compact sets.

**Proof of Theorem 4**

**Step 1** *Proof of 1[st] assertion* The $\mathcal{F}-$integrals converge absolutely, for every $(x,t)\in\overline{Q}$, since their integrands are $O(1/\lambda^2)$, as $|\lambda|\to+\infty$ with $\lambda\in\mathbb{R}$. Indeed, for example, for the integrand of the first integral of $\mathcal{F}_{1,1}(x,t)$, we have

$$\text{(7.13)}\qquad \frac{e^{i\lambda x-\omega_1(\lambda)t}\{\lambda\tilde{\hat{f}}_2(\lambda,\omega_1,t)+[i\omega_2(\lambda)+m]\tilde{\hat{f}}_1(\lambda,\omega_1,t)\}}{2\rho(\lambda)}=O(1/\lambda^2),$$

as it follows from (7.11).

Furthermore, since (7.13) and all of their analogues are uniform for $(x,t)$ in bounded subsets of $\overline{Q}$, the functions $\mathcal{F}_{1,1}(x,t)$, $\mathcal{F}_{1,2}(x,t)$, $\mathcal{F}_{2,1}(x,t)$ and $\mathcal{F}_{2,2}(x,t)$, are, indeed, continuous for $(x,t)\in\overline{Q}$.

**Step 2** *Analysis of the* $\mathcal{F}-$*integrals* As in section 3, we can rewrite these integrals by splitting them and deforming the contours of the various pieces in such a way, so that we can conclude that the functions $\mathcal{F}_{1,1}(x,t)$, $\mathcal{F}_{1,2}(x,t)$, $\mathcal{F}_{2,1}(x,t)$ and $\mathcal{F}_{2,2}(x,t)$, are $C^{\infty}$ for $(x,t)\in Q-\{x=t\}$. The main tools in this part of the proof are the relations (7.7), (7.8), (7.9) and (7.10).
As an example of these interpretations, let us consider the integral

$$(\mathcal{F}_{1,2})_{2,2}(x,t):=\int_{-\infty}^{\infty} e^{i\lambda x-\omega_2(\lambda)t}[i\omega_1(\lambda)+m]\tilde{\hat{f}}_1(\lambda,\omega_2,t)\frac{d\lambda}{2\rho(\lambda)},$$

which is part of $(\mathcal{F}_{1,2})_2(x,t)$, and let us split it as follows:

$$(\mathcal{F}_{1,2})_{2,2}(x,t)=\int_{\lambda=-\infty}^{-1}+\int_{\lambda=-1}^{1}+\int_{\lambda=1}^{\infty}=:\mathfrak{I}^-(x,t)+\mathfrak{I}_0(x,t)+\mathfrak{I}^+(x,t).$$

Then, recalling that

$$e^{i\lambda x-\omega_2(\lambda)t}=e^{i\lambda x-i\lambda t}e^{i\lambda t+i\rho(\lambda)} \text{ and } \sup_{\substack{|\lambda|\to+\infty\\ \operatorname{Re}\lambda<0}}\left|e^{i\lambda t+i\rho(\lambda)}\right|<+\infty,$$

we have:

$$\begin{aligned}
(7.14)\ \text{For } x>t,\ \mathfrak{I}^-(x,t)&=\int_{-\infty}^{-1}e^{i\lambda x-\omega_2(\lambda)t}[i\omega_1(\lambda)+m]\tilde{\hat{f}}_1(\lambda,\omega_2,t)d\lambda\\
&=\int_{\Gamma^{-,+}}e^{i\lambda x}[i\omega_1(\lambda)+m]\mu_{f_1,M,N}(\lambda,\omega_2,t)d\lambda\\
&\quad-\int_{\Gamma^{-,+}}e^{i\lambda x-\omega_2(\lambda)t}[i\omega_1(\lambda)+m]\mu_{f_1,M,N}(\lambda,\omega_2,0)d\lambda\\
&\quad-\int_{-\infty}^{-1}e^{i\lambda x}[i\omega_1(\lambda)+m]\left\{\frac{e^{-\omega_2(\lambda)t}}{[\omega_2(\lambda)]^N}\int_{\tau=0}^{t}e^{\omega_2(\lambda)\tau}P_{f_1,M}^{(N)}(\lambda,\omega_2(\lambda),t)d\tau\right\}d\lambda\\
&\quad+\int_{-\infty}^{-1}e^{i\lambda x}[i\omega_1(\lambda)+m]\left\{\frac{e^{-\omega_2(\lambda)t}}{(i\lambda)^M}\int_{\tau=0}^{t}\left[\int_{y=0}^{\infty}e^{-i\lambda y}\frac{\partial^M f_1(y,\tau)}{\partial y^M}dy\right]e^{\omega_2(\lambda)\tau}\right\}d\lambda,
\end{aligned}$$

while,

$$\begin{aligned}
(7.15)\qquad \text{for } x<t,\ \mathfrak{I}^-(x,t)&=\int_{-\infty}^{-1}e^{i\lambda x-\omega_2(\lambda)t}[i\omega_1(\lambda)+m]\tilde{\hat{f}}_1(\lambda,\omega_2,t)d\lambda\\
&=\int_{\Gamma^{-,-}}e^{i\lambda x-\omega_2(\lambda)t}[i\omega_1(\lambda)+m]\tilde{\hat{f}}_1(\lambda,\omega_2,t)d\lambda.
\end{aligned}$$

Applying (7.14) and (7.15) with sufficient large $M$ and $N$, we conclude that

$$\mathfrak{I}^-(x,t)\in C^{\infty}(Q-\{x=t\}).$$

Also (7.14) and (7.15) lead to formulas for the derivatives of $\mathfrak{I}^-(x,t)$. For example,

$$(7.16)\ \text{For } x>t,\ \frac{\partial^{n+l}\mathfrak{I}^-(x,t)}{\partial x^n\partial t^l}(x,t)=\int_{\Gamma^{-,+}}(i\lambda)^n e^{i\lambda x}[i\omega_1(\lambda)+m]\mu_{f_1,M,N}(\lambda,\omega_2,t)d\lambda$$

$$- \int_{\Gamma^{-,+}} (i\lambda)^n [-\omega_2(\lambda)]^l e^{i\lambda x - \omega_2(\lambda)t} [i\omega_1(\lambda) + m] \mu_{f_1,M,N}(\lambda, \omega_2, 0) d\lambda$$

$$- \int_{-\infty}^{-1} (i\lambda)^n e^{i\lambda x} [i\omega_1(\lambda) + m] \frac{\partial^l}{\partial t^l} \left\{ \frac{e^{-\omega_2(\lambda)t}}{[\omega_2(\lambda)]^N} \int_0^t e^{\omega_2(\lambda)\tau} P_{f_1,M}^{(N)}(\lambda, \omega_2(\lambda), t) d\tau \right\} d\lambda$$

$$+ \int_{-\infty}^{-1} (i\lambda)^n e^{i\lambda x} [i\omega_1(\lambda) + m] \frac{\partial^l}{\partial t^l} \left\{ \frac{e^{-\omega_2(\lambda)t}}{(i\lambda)^M} \int_0^t \left[ \int_0^\infty e^{-i\lambda y} \frac{\partial^M f_1(y,\tau)}{\partial y^M} dy \right] e^{\omega_2(\lambda)\tau} d\tau \right\} d\lambda,$$

with the all the above integrals being absolutely convergent.

Similar computations can be made also for the integral $\Im^+(x,t)$, while the contour of the integral $\Im_0(x,t)$ is finite. Thus, $(\mathcal{F}_{1,2})_{2,2}(x,t) \in C^\infty(Q - \{x = t\})$.

Similar analysis can be carried out for all the parts of the $\mathcal{F}-$integrals, and this easily leads to the conclusion that, indeed, the functions $\Psi_1(x,t)$ and $\Psi_2(x,t)$ satisfy the analogues of the 2nd, 4th, 5th and 6th conclusions of Theorem 1.

**Step 3** *Proof of 3rd assertion* Firstly, let us prove that

$$\partial_t^* [\mathcal{F}_{1,1}(x,t) + \mathcal{F}_{1,2}(x,t)] = f_1(x,t), \tag{7.17}$$

where $\partial_t^*$ is the differentiation $\partial / \partial t$, acting ***only*** on $t$ in the terms $\tilde{\hat{f}}_1(\lambda, \omega_j, t)$ and $\tilde{\hat{f}}_2(\lambda, \omega_j, t)$ of the $\mathcal{F}-$integrals. Having at our disposal the analysis of the $\mathcal{F}-$integrals, carried out in step 2, we can compute:

(7.18) $\partial_t^* [(\mathcal{F}_{1,1})_2(x,t) - (\mathcal{F}_{1,2})_2(x,t)]$

$$= \partial_t^* \left\{ \int_{-\infty}^{\infty} e^{i\lambda x - \omega_1(\lambda)t} \{\lambda \tilde{\hat{f}}_2(-\lambda, \omega_1, t) - [i\omega_2(\lambda) + m] \tilde{\hat{f}}_1(-\lambda, \omega_1, t)\} \frac{d\lambda}{2\rho(\lambda)} \right\}$$

$$- \partial_t^* \left\{ \int_{-\infty}^{\infty} e^{i\lambda x - \omega_2(\lambda)t} \{\lambda \tilde{\hat{f}}_2(-\lambda, \omega_2, t) - [i\omega_1(\lambda) + m] \tilde{\hat{f}}_1(-\lambda, \omega_2, t)\} \frac{d\lambda}{2\rho(\lambda)} \right\}$$

$$= \int_{-\infty}^{\infty} e^{i\lambda x} \{\lambda \hat{f}_2(-\lambda, t) - [i\omega_2(\lambda) + m] \hat{f}_1(-\lambda, t)\} \frac{d\lambda}{2\rho(\lambda)} - \int_{-\infty}^{\infty} e^{i\lambda x} \{\lambda \hat{f}_2(-\lambda, t) - [i\omega_1(\lambda) + m] \hat{f}_1(-\lambda, t)\} \frac{d\lambda}{2\rho(\lambda)}$$

$$= - \int_{-\infty}^{\infty} e^{i\lambda x} \{[i\omega_2(\lambda) - i\omega_1(\lambda)] \hat{f}_1(-\lambda, t)\} \frac{d\lambda}{2\rho(\lambda)} = - \int_{-\infty}^{\infty} e^{i\lambda x} \hat{f}_1(-\lambda, t)\} d\lambda = 0,$$

where the last equation follows from Cauchy's theorem and Jordan's lemma.

*Note* We point out that in order to justify the interchange of the operator $\partial_t^*$ with the integrations, in the proof of (7.18), we rely on the analysis of the $\mathcal{F}-$integrals carried out in Step 2. Indeed, this is because, after the interchange of the operator $\partial_t^*$ with the integrations, the resulting intagrals do not converge – in general – absolutely. However, they exist in the generalized sense, and this justifies the computation leading to (7.18).

Similarly,

(7.19) $\partial_t^* [(\mathcal{F}_{1,1})_1(x,t) - (\mathcal{F}_{1,2})_1(x,t)]$

$$= \partial_t^* \left\{ \int_{-\infty}^{\infty} e^{i\lambda x - \omega_1(\lambda)t} \{\lambda \tilde{\hat{f}}_2(\lambda, \omega_1, t) + [i\omega_2(\lambda) + m] \tilde{\hat{f}}_1(\lambda, \omega_1, t)\} \frac{d\lambda}{2\rho(\lambda)} \right\}$$

$$- \partial_t^* \left\{ \int_{-\infty}^{\infty} e^{i\lambda x - \omega_2(\lambda)t} \{\lambda \tilde{\hat{f}}_2(\lambda, \omega_2, t) + [i\omega_1(\lambda) + m] \tilde{\hat{f}}_1(\lambda, \omega_2, t)\} \frac{d\lambda}{2\rho(\lambda)} \right\}$$

$$= \int_{-\infty}^{\infty} e^{i\lambda x} \{\lambda \hat{f}_2(\lambda, t) + [i\omega_2(\lambda) + m] \hat{f}_1(-\lambda, t)\} \frac{d\lambda}{2\rho(\lambda)} - \int_{-\infty}^{\infty} e^{i\lambda x} \{\lambda \hat{f}_2(\lambda, t) + [i\omega_1(\lambda) + m] \hat{f}_1(-\lambda, t)\} \frac{d\lambda}{2\rho(\lambda)}$$

$$= \int_{-\infty}^{\infty} e^{i\lambda x} \{[i\omega_2(\lambda) - i\omega_1(\lambda)] \hat{f}_1(\lambda, t)\} \frac{d\lambda}{2\rho(\lambda)} = \int_{-\infty}^{\infty} e^{i\lambda x} \hat{f}_1(\lambda, t)\} d\lambda = 2\pi f_1(x, t).$$

Thus, (7.18) and (7.19) imply (7.17).

Similarly, we show that

$$\partial_t^* [-(\mathcal{F}_{2,1})_2(x,t) + (\mathcal{F}_{2,2})_2(x,t)] = 0 \text{ and } \partial_t^* [(\mathcal{F}_{2,1})_1(x,t) - (\mathcal{F}_{2,2})_1(x,t)] = 2\pi f_2(x,t),$$

which implies

$$\partial_t^* [\mathcal{F}_{2,1}(x,t) + \mathcal{F}_{2,2}(x,t)] = 2\pi f_2(x,t). \tag{7.20}$$

Finally, combining (7.16) and (7.19), with the fact that the basic exponential solutions (4.1) satisfy the homogeneous equations of the system, we easily complete the proof of the 3[rd] assertion.

**Step 4** *Proof of 4[th] assertion* Working as in Step 5, of the proof of Theorem 1, we find that if $h_1(0) = g_1(0)$ then, for example,

$$\frac{\partial}{\partial x} [\mathcal{G}_{1,3}^-(x,t) + \mathcal{G}_{1,4}^-(x,t) + \mathcal{H}_{1,2}^-(x,t) + (\mathcal{F}_{1,2})_1^-(x,t) + (\mathcal{F}_{1,2})_2^-(x,t)]$$
$$\approx [-ig_2'(0) - ih_1'(0) + if_1(0,0)] \int_{-\infty}^{-1} e^{i\lambda x - \omega_2(\lambda)t} \frac{d\lambda}{\rho(\lambda)} \tag{7.21}$$

and

$$\frac{\partial}{\partial t} [\mathcal{G}_{2,1}^+(x,t) - \mathcal{G}_{2,2}^+(x,t) + \mathcal{H}_{2,1}^+(x,t) + (\mathcal{F}_{2,1})_1^+(x,t) - (\mathcal{F}_{2,1})_2^+(x,t)]$$
$$\approx [ig_2'(0) + ih_1'(0) - if_1(0,0)] \int_{1}^{\infty} e^{i\lambda x - \omega_1(\lambda)t} \frac{d\lambda}{\rho(\lambda)}. \tag{7.22}$$

Similar computations can be carried out for all the analogous combinations which are parts of the terms

$$\frac{\partial \Psi_j(x,t)}{\partial x}, \; \frac{\partial \Psi_j(x,t)}{\partial t}, \; j = 1, 2,$$

and these computations show that the obstacles for the existence of the limits

$$\lim_{\substack{(x,t)\to(a,a) \\ x\neq t}} \frac{\partial \Psi_j(x,t)}{\partial x}, \; \lim_{\substack{(x,t)\to(a,a) \\ x\neq t}} \frac{\partial \Psi_j(x,t)}{\partial t}, \; a > 0, \; j = 1, 2,$$

which are the integrals

$$\int_{-\infty}^{-1} e^{i\lambda x-\omega_1(\lambda)t}\,\frac{d\lambda}{\rho(\lambda)},\ \int_{-\infty}^{-1} e^{i\lambda x-\omega_2(\lambda)t}\,\frac{d\lambda}{\rho(\lambda)},\ \int_{1}^{\infty} e^{i\lambda x-\omega_1(\lambda)t}\,\frac{d\lambda}{\rho(\lambda)} \text{ and } \int_{1}^{\infty} e^{i\lambda x-\omega_2(\lambda)t}\,\frac{d\lambda}{\rho(\lambda)},$$

are annihilated, if $g_1(0)-h_1(0)=0$ and $g_2'(0)+h_1'(0)-f_1(0,0)$.

Thus, if these obstacles cease to exist, we are led to the formulas for $\Psi_1(x,t)$ and $\Psi_2(x,t)$, valid also for $x=t$, as stated in the $4^{\text{th}}$ part of the theorem. Now we can easily complete the proof of $4^{\text{th}}$ assertion. It suffices, for example, to notice that, if $g_1(0)-h_1(0)=0$ and $g_2'(0)+h_1'(0)-f_1(0,0)=0$, the computations which yielded (7.21) and (7.22), imply, furthermore, that the sums of the integrands

$$\{\mathcal{I}(\mathcal{G}_{1,3}^-)+\mathcal{I}(\mathcal{G}_{1,4}^-)+\mathcal{I}(\mathcal{H}_{1,2}^-)+\mathcal{I}[(\mathcal{F}_{1,2})_1^-]+\mathcal{I}[(\mathcal{F}_{1,2})_2^-]\}(\lambda,x,t)$$

and

$$\{\mathcal{I}(\mathcal{G}_{2,1}^+)-\mathcal{I}(\mathcal{G}_{2,2}^+)+\mathcal{I}(\mathcal{H}_{2,1}^+)+\mathcal{I}[(\mathcal{F}_{2,1})_1^+]-\mathcal{I}[\mathcal{F}_{2,1})_2^+]\}(\lambda,x,t),$$

are $\mathrm{O}(1/\lambda^3)$, as $\lambda\to-\infty$ or $\lambda\to+\infty$, respectively, uniformly for $(x,t)$ in compact subsets of $Q$. (Analogous conclusions have to be drawn for all those appropriate combinations, until we exhaust the integrands in the RHSs of (7.5) and (7.6).)

The following theorem is the analogue of Theorem 2, in the inhomogeneous case.

**Theorem 5** *Assuming (1.3) and (7.3), and letting $\Psi_1(x,t)$ and $\Psi_2(x,t)$ be as in Theorem 4, the following hold:*

*$1^{st}$ For nonnegative integers $n$ and $l$, the limits*

$$\lim_{t\to0^+}\frac{\partial^{n,l}\Psi_1(x,t)}{\partial x^n\partial t^l} \ \textit{and}\ \lim_{t\to0^+}\frac{\partial^{n,l}\Psi_2(x,t)}{\partial x^n\partial t^l}$$

*exist, uniformly for $x$ in compact subsets of $\mathbb{R}^+$ and define $C^\infty$ functions of $x\in\mathbb{R}^+$.*

*$2^{nd}$ For nonnegative integers $n$ and $l$, the limits*

$$\lim_{x\to0^+}\frac{\partial^{n,l}\Psi_1(x,t)}{\partial x^n\partial t^l} \ \textit{and}\ \lim_{x\to0^+}\frac{\partial^{n,l}\Psi_2(x,t)}{\partial x^n\partial t^l}$$

*exist, uniformly for $t$ in compact subsets of $\mathbb{R}^+$ and define $C^\infty$ functions of $t\in\mathbb{R}^+$.*

*$3^{rd}$ The functions $\Psi_1(x,t)$ and $\Psi_2(x,t)$ extend to $C^\infty$ functions in open neighborhoods of the open half-lines $\{(x,0):x>0\}$ and $\{(0,t):t>0\}$.*

*$4^{th}$* $\lim_{t\to0^+}\frac{\partial^n\Psi_1(x,t)}{\partial x^n}=\frac{d^n g_1(x)}{dx^n}$ *and* $\lim_{t\to0^+}\frac{\partial^n\Psi_2(x,t)}{\partial x^n}=\frac{d^n g_2(x)}{dx^n}$, *for* $x\in\mathbb{R}^+$.

*$5^{th}$* $\lim_{x\to0^+}\frac{\partial^l\Psi_1(x,t)}{\partial t^l}=\frac{d^l h_1(t)}{dt^l}$, *for* $t\in\mathbb{R}^+$.

**Proof Step 1** Using the analysis of the $\mathcal{F}-$integrals, carried out in Step 2, of the proof of Theorem 4, the proof of $1^{\text{st}}$ assertion is similar to the proof of the corresponding part of Theorem 2. However, we have some difficulty to address for the proof of the $2^{\text{nd}}$ assertion.

**Step 2** *Proof of $2^{nd}$ assertion* In view of formulas (7.9) and (7.10), and recalling that $\omega_1=-\omega_2=i\rho$, we find that

(7.23) $$\frac{e^{i\lambda x-\omega_1(\lambda)t}\{\lambda\widetilde{\widetilde{f}}_2(\lambda,\omega_1,t)-e^{i\lambda x-\omega_2(\lambda)t}\{\lambda\widetilde{\widetilde{f}}_2(\lambda,\omega_2,t)}{2\rho(\lambda)}$$

$$=\frac{e^{i\lambda x-\omega_1(\lambda)t}\lambda\widetilde{P}_{f_2,M}(\lambda,\omega_1(\lambda),t)-e^{i\lambda x-\omega_2(\lambda)t}\lambda\widetilde{P}_{f_2,M}(\lambda,\omega_2(\lambda),t)}{2\rho(\lambda)}$$

$$+\frac{e^{i\lambda x-\omega_1(\lambda)t}\lambda\widetilde{\Upsilon}_{f_2,M}(\lambda,\omega_1(\lambda),t)-e^{i\lambda x-\omega_2(\lambda)t}\lambda\widetilde{\Upsilon}_{f_2,M}(\lambda,\omega_2(\lambda),t)}{2\rho(\lambda)}$$

$$=e^{i\lambda x}\lambda\sum_{s=1}^{N}\frac{(P_{f_2,M})^{(2s-2)}(\lambda,t)}{\omega_1^{2s-1}(\lambda)\rho(\lambda)}-\lambda\frac{e^{i\lambda x-\omega_1(\lambda)t}\mu_{f_2,M,2N}(\lambda,\omega_1(\lambda),0)-e^{i\lambda x-\omega_2(\lambda)t}\mu_{f_2,M,2N}(\lambda,\omega_2(\lambda),0)}{2\rho(\lambda)}$$

$$-\frac{\lambda}{2\rho(\lambda)}\left\{\frac{e^{i\lambda x-\omega_1(\lambda)t}}{[\omega_1(\lambda)]^{2N}}\int_0^t e^{\omega_1(\lambda)\tau}P_{f_2,M}^{(2N)}(\lambda,t)d\tau-\frac{e^{i\lambda x-\omega_2(\lambda)t}}{[\omega_2(\lambda)]^{2N}}\int_{\tau=0}^t e^{\omega_2(\lambda)\tau}P_{f_2,M}^{(2N)}(\lambda,t)d\tau\right\}$$

$$+\frac{\lambda e^{i\lambda x}}{2\rho(\lambda)}\left\{\frac{e^{-\omega_1(\lambda)t}}{(i\lambda)^M}\int_0^t\left[e^{\omega_1(\lambda)\tau}\int_0^\infty e^{-i\lambda y}\frac{\partial^M f_2(y,\tau)}{\partial y^M}dy\right]d\tau-\frac{e^{-\omega_2(\lambda)t}}{(i\lambda)^M}\int_0^t\left[e^{\omega_2(\lambda)\tau}\int_0^\infty e^{-i\lambda y}\frac{\partial^M f_2(y,\tau)}{\partial y^M}dy\right]d\tau\right\}.$$

Next we observe – and this is crucial – that

$$\lambda\sum_{s=1}^{N}\frac{(P_{f_2,M})^{(2s-2)}(\lambda,t)}{\omega_1^{2s-1}(\lambda)\rho(\lambda)}=\lambda\sum_{s=1}^{N}\frac{(P_{f_2,M})^{(2s-2)}(\lambda,t)}{i^{2s-1}(\lambda^2+m^2)^s}\text{ is analytic in }\mathbb{C}-\{mi,-mi,0\},$$

and, therefore,

(7.24) $$\int_{\Gamma^{-,+}+\Gamma^{+,+}}e^{i\lambda x}\lambda\sum_{s=1}^{N}\frac{(P_{f_2,M})^{(2s-2)}(\lambda,t)}{\omega_1^{2s-1}(\lambda)\rho(\lambda)}d\lambda=\int_{[1,1+i]+\gamma+[-1+i,-1]}e^{i\lambda x}\lambda\sum_{s=1}^{N}\frac{(P_{f_2,M})^{(2s-2)}(\lambda,t)}{i^{2s-1}(\lambda^2+m^2)^s}d\lambda,$$

where $\gamma$ is a simple contour starting at the point $1+i$, going above the point $mi$, and ending at the point $-1+i$, as depicted in fig.8.

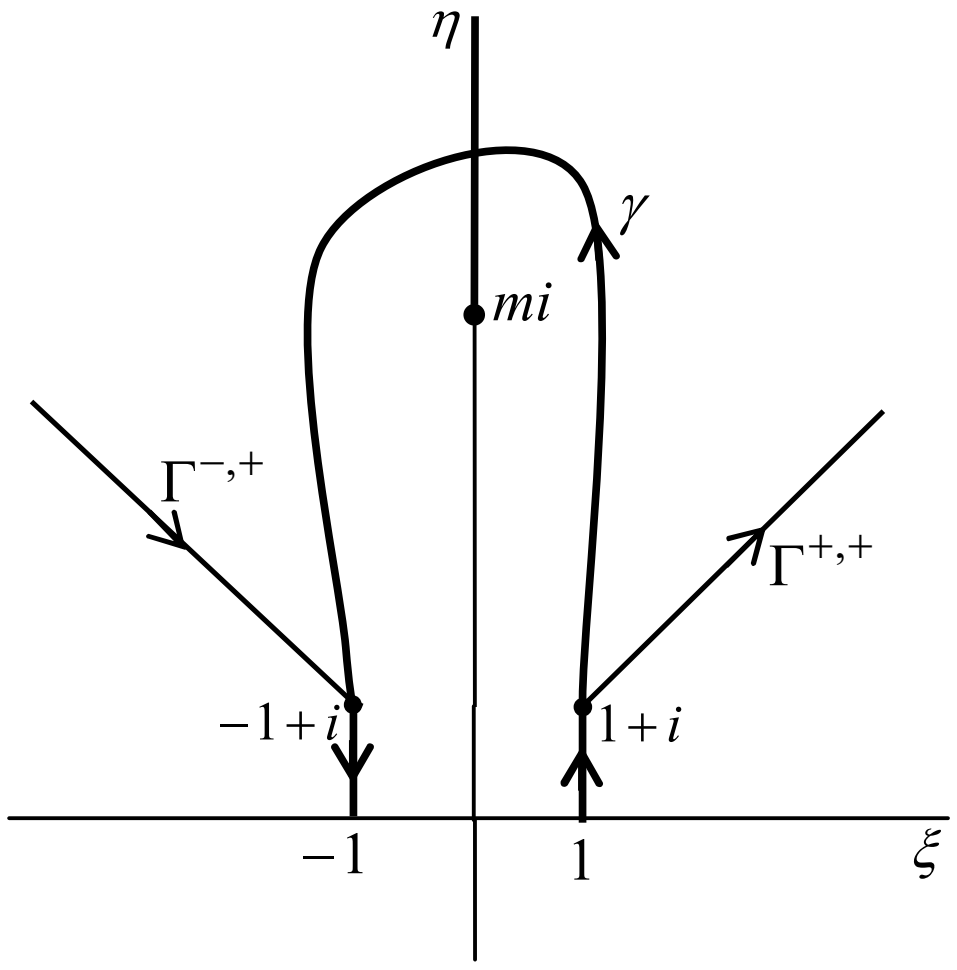

**Fig.8** *A typical contour* $\gamma$ *and the corresponding contour* $[1,1+i]+\gamma+[-1+i,-1]$

Recalling that, in this part of the proof, we work for $x<t$, we obtain, by (7.23) and (7.24), that

$$(7.25)\quad \left(\int_{-\infty}^{-1}+\int_{1}^{\infty}\right)\frac{e^{i\lambda x-\omega_1(\lambda)t}\{\lambda\tilde{\hat{f}}_2(\lambda,\omega_1,t)-e^{i\lambda x-\omega_2(\lambda)t}\{\lambda\tilde{\hat{f}}_2(\lambda,\omega_2,t)}{2\rho(\lambda)}d\lambda$$

$$=\int_{[-1+i,\,-1]+\gamma+[1,\,1+i]} e^{i\lambda x}\lambda\sum_{s=1}^{N}\frac{(P_{f_2,M})^{(2s-2)}(\lambda,t)}{i^{2s-1}(\lambda^2+m^2)^s}d\lambda$$

$$-\int_{\Gamma^{-,+}+\Gamma^{+,-}}\lambda\frac{e^{i\lambda x-\omega_1(\lambda)t}\mu_{f_2,M,2N}(\lambda,\omega_1(\lambda),0)}{2\rho(\lambda)}d\lambda+\int_{\Gamma^{-,-}+\Gamma^{+,+}}\lambda\frac{e^{i\lambda x-\omega_2(\lambda)t}\mu_{f_2,M,2N}(\lambda,\omega_2(\lambda),0)}{2\rho(\lambda)}d\lambda$$

$$-\left(\int_{-\infty}^{-1}+\int_{1}^{\infty}\right)\left[\frac{\lambda}{2\rho(\lambda)}\left\{\frac{e^{i\lambda x-\omega_1(\lambda)t}}{[\omega_1(\lambda)]^{2N}}\int_{\tau=0}^{t}e^{\omega_1(\lambda)\tau}P^{(2N)}_{f_2,M}(\lambda,t)d\tau-\frac{e^{i\lambda x-\omega_2(\lambda)t}}{[\omega_2(\lambda)]^{2N}}\int_{\tau=0}^{t}e^{\omega_2(\lambda)\tau}P^{(2N)}_{f_2,M}(\lambda,t)d\tau\right\}\right]d\lambda$$

$$+\left(\int_{-\infty}^{-1}+\int_{1}^{\infty}\right)\left[\frac{\lambda}{2\rho(\lambda)}\frac{e^{i\lambda x-\omega_1(\lambda)t}}{(i\lambda)^M}\int_0^t\left\{e^{\omega_1(\lambda)\tau}\int_0^\infty e^{-i\lambda y}\frac{\partial^M f_2(y,\tau)}{\partial y^M}dy\right\}d\tau\right]d\lambda$$

$$-\left(\int_{-\infty}^{-1}+\int_{1}^{\infty}\right)\left[\frac{\lambda}{2\rho(\lambda)}\frac{e^{i\lambda x-\omega_2(\lambda)t}}{(i\lambda)^M}\int_0^t\left\{e^{\omega_2(\lambda)\tau}\int_0^\infty e^{-i\lambda y}\frac{\partial^M f_2(y,\tau)}{\partial y^M}dy\right\}d\tau\right]d\lambda.$$

Now, there is no difficulty in applying the differentiation $\partial^{n,l}/\partial x^n\partial t^l$ to both sides of (7.25), interchanging it with every integral in the RHS, and then letting $x\to 0^+$, proving that the limit

$$\lim_{x\to0^+}\frac{\partial^{n,l}}{\partial x^n\partial t^l}\left\{\left(\int_{-\infty}^{-1}+\int_{1}^{\infty}\right)\frac{e^{i\lambda x-\omega_1(\lambda)t}\{\lambda\tilde{\hat{f}}_2(\lambda,\omega_1,t)-e^{i\lambda x-\omega_2(\lambda)t}\{\lambda\tilde{\hat{f}}_2(\lambda,\omega_2,t)}{2\rho(\lambda)}d\lambda\right\}$$

exists, uniformly for $t$ in compact subsets of $\mathbb{R}^+$ and that it defines a $C^\infty$ function of $t\in\mathbb{R}^+$. (We point out that we have to choose, *gradually*, $M$ and $N$, sufficiently large.)

The main goal of the computation, which yielded (7.25), was to reduce the integral in the LHS of (7.24) to the integral in the RHS, whose contour, namely $[-1+i,-1]+\gamma+[1,1+i]$, is finite, so that the interchange of the process

$$\frac{\partial^\ell}{\partial t^\ell}\left(\lim_{x\to0^+}\frac{\partial^{n,l}}{\partial x^n\partial t^l}\right)$$

with the integration, is justified, by Lebesgue's dominated convergence theorem. To achieve this reduction, we paired up two parts of $\mathcal{F}_{1,1}(x,t)+\mathcal{F}_{1,2}(x,t)$, namely

$$\left(\int_{-\infty}^{-1}+\int_{1}^{\infty}\right)\frac{e^{i\lambda x-\omega_1(\lambda)t}\{\lambda\tilde{\hat{f}}_2(\lambda,\omega_1,t)\}}{2\rho(\lambda)}d\lambda\quad\text{with}\quad-\left(\int_{-\infty}^{-1}+\int_{1}^{\infty}\right)\frac{e^{i\lambda x-\omega_2(\lambda)t}\{\lambda\tilde{\hat{f}}_2(\lambda,\omega_2,t)\}}{2\rho(\lambda)}d\lambda.$$

Such combinations can be made, until we exhaust all of $\mathcal{F}_{1,1}(x,t)+\mathcal{F}_{1,2}(x,t)$, and similarly all of $\mathcal{F}_{2,1}(x,t)+\mathcal{F}_{2,2}(x,t)$, leading to the 2nd conclusion.

**Step 3** *Completion of the proof of the theorem* Having at our disposal the 1st and 2nd conclusions, the proofs of 3rd, 4th and 5th assertions are similar to the corresponding parts of Theorem 2.

The following theorem is the analogue of Theorem 3, in the inhomogeneous case, and its proof is similar.

**Theorem 6** *Assuming (1.3) and (7.3), and letting* $\Psi_1(x,t)$ *and* $\Psi_2(x,t)$ *be as in Theorem 4, the following hold: For* $T>0$,

$$\lim_{x\to+\infty}\left[x^E\,\frac{\partial^{n+l}\Psi_1(x,t)}{\partial x^n\partial t^l}\right]=0 \ \ and \ \ \lim_{x\to+\infty}\left[x^E\,\frac{\partial^{n+l}\Psi_2(x,t)}{\partial x^n\partial t^l}\right]=0,$$

*for nonnegative integers* $n$, $l$ *and* $E$, *uniformly for* $0<t\le T$.

## 8. The behaviour of the solution, as $t\to\infty$, with periodic data

In this section, we consider the solution of problem (7.1) & (7.2), given by Theorem 4, in the case the data $h_1(t)$, $f_1(x,t)$ and $f_2(x,t)$ are periodic in $t$, and study its asymptotic behaviour as $t\to\infty$. More precisely, we will assume, in addition to (1.3) and (7.3), that, for a fixed $T>0$,

(8.1) $\quad h_1(t+T)=h_1(t)$, $f_1(x,t)=f_1(x,t+T)$ and $f_2(x,t)=f_2(x,t+T)$ for every $t\ge 0$ and $x\ge 0$,

and we will prove the following theorem.

**Theorem 7** *Assuming (1.3), (7.3) and (8.1), the functions* $\Psi_1(x,t)$ *and* $\Psi_2(x,t)$, *defined in Theorem 4, satisfy the following: For fixed* $X>0$,

$$\lim_{t\to+\infty}\left\{t^E\,\frac{\partial^{n+l}[\Psi_1(x,t+T)-\Psi_1(x,t)]}{\partial x^n\partial t^l}\right\}=0 \ \ and \ \ \lim_{t\to+\infty}\left\{t^E\,\frac{\partial^{n+l}[\Psi_2(x,t+T)-\Psi_2(x,t)]}{\partial x^n\partial t^l}\right\}=0,$$

*for nonnegative integers* $n$, $l$ *and* $E$, *uniformly for* $0<x\le X$.

The conclusion of Theorem 7 follows immediately from the following lemmas. Let us point out that in this section we will be working with $t>>X$ and $0<x\le X$. In particular, for the deformation of certain integrals we will have $t>>x$.

**Lemma 1** *Let* $\mathcal{G}_1(x,t)$ *be the part of* $\Psi_1(x,t)$ *which involves* $g_1$ *and* $g_2$, *i.e., with the notation of (1.4),* $2\pi\mathcal{G}_1=\mathcal{G}_{1,1}+\mathcal{G}_{1,2}-\mathcal{G}_{1,3}-\mathcal{G}_{1,4}$. *Then, for fixed* $X>0$,

(8.2)
$$\lim_{t\to+\infty}[t^E\mathcal{G}_1(x,t)]=0,$$

*for every nonnegative integer* $E$, *uniformly for* $0<x\le X$.
*Similarly, if* $\mathcal{G}_2(x,t)$ *is the part of* $\Psi_2(x,t)$ *which involves* $g_1$ *and* $g_2$, *then*

(8.3)
$$\lim_{t\to+\infty}[t^E\mathcal{G}_2(x,t)]=0,$$

*for every nonnegative integer* $E$, *uniformly for* $0<x\le X$.

**Proof** Let us consider the integral

$$I_1(x,t)=\int_{-\infty}^{\infty}e^{i\lambda x-\omega_1(\lambda)t}\{\lambda\hat g_2(\lambda)\}\frac{d\lambda}{2\rho(\lambda)},$$

which is a part of $\mathcal{G}_1(x,t)$. Then, deforming the contours and integrating by parts, we obtain

(8.4)
$$I_1(x,t)=\left(\int_{-\infty}^{-1}+\int_{-1}^{1}+\int_{1}^{\infty}\right)e^{i\lambda x-\omega_1(\lambda)t}\{\lambda\hat g_2(\lambda)\}\frac{d\lambda}{2\rho(\lambda)}$$

$$= \int_{\Gamma^{-,+}} e^{i\lambda x-\omega_1(\lambda)t}\left\{\lambda\frac{g_2(0)}{i\lambda}\right\}\frac{d\lambda}{2\rho(\lambda)}+\int_{-\infty}^{-1} e^{i\lambda x-\omega_1(\lambda)t}\left\{\lambda\frac{(g_2^{(1)})\hat{}(\lambda)}{i\lambda}\right\}\frac{d\lambda}{2\rho(\lambda)}$$

$$+\int_{-1}^{1} e^{i\lambda x-\omega_1(\lambda)t}\{\lambda\hat{g}_2(\lambda)\}\frac{d\lambda}{2\rho(\lambda)}+\int_{\Gamma^{+,-}} e^{i\lambda x-\omega_1(\lambda)t}\{\lambda\hat{g}_2(\lambda)\}\frac{d\lambda}{2\rho(\lambda)}$$

$$=-\frac{1}{t}\int_{\Gamma^{-,+}}\frac{1}{\omega_1'(\lambda)}\frac{d[e^{-\omega_1(\lambda)t}]}{d\lambda}e^{i\lambda x}\left\{\lambda\frac{g_2(0)}{i\lambda}\right\}\frac{d\lambda}{2\rho(\lambda)}-\frac{1}{t}\int_{-\infty}^{-1}\frac{1}{\omega_1'(\lambda)}\frac{d[e^{-\omega_1(\lambda)t}]}{d\lambda}e^{i\lambda x}\left\{\lambda\frac{(g_2^{(1)})\hat{}(\lambda)}{i\lambda}\right\}\frac{d\lambda}{2\rho(\lambda)}$$

$$-\frac{1}{t}\int_{-1}^{1}\frac{1}{\omega_1'(\lambda)}\frac{d[e^{-\omega_1(\lambda)t}]}{d\lambda}e^{i\lambda x}\{\lambda\hat{g}_2(\lambda)\}\frac{d\lambda}{2\rho(\lambda)}-\frac{1}{t}\int_{\Gamma^{+,-}}\frac{1}{\omega_1'(\lambda)}\frac{d[e^{-\omega_1(\lambda)t}]}{d\lambda}e^{i\lambda x}\{\lambda\hat{g}_2(\lambda)\}\frac{d\lambda}{2\rho(\lambda)}$$

$$=\frac{1}{t}\int_{\Gamma^{-,+}}\frac{1}{\omega_1'(\lambda)}e^{-\omega_1(\lambda)t}\frac{d}{d\lambda}\left[e^{i\lambda x}\left\{\lambda\frac{g_2(0)}{i\lambda}\right\}\frac{1}{2\rho(\lambda)}\right]d\lambda$$

$$+\frac{1}{t}\int_{-\infty}^{-1}\frac{1}{\omega_1'(\lambda)}e^{-\omega_1(\lambda)t}\frac{d}{d\lambda}\left[e^{i\lambda x}\left\{\lambda\frac{(g_2^{(1)})\hat{}(\lambda)}{i\lambda}\right\}\frac{1}{2\rho(\lambda)}\right]d\lambda$$

$$+\frac{1}{t}\int_{-1}^{1}\frac{1}{\omega_1'(\lambda)}e^{-\omega_1(\lambda)t}\frac{d}{d\lambda}\left[e^{i\lambda x}\{\lambda\hat{g}_2(\lambda)\}\frac{1}{2\rho(\lambda)}\right]d\lambda$$

$$+\frac{1}{t}\int_{\Gamma^{+,-}}\frac{1}{\omega_1'(\lambda)}e^{-\omega_1(\lambda)t}\frac{d}{d\lambda}\left[e^{i\lambda x}\{\lambda\hat{g}_2(\lambda)\}\frac{1}{2\rho(\lambda)}\right]d\lambda.$$

Now, observing that the last four integrals, in the above formula (8.4), are absolutely convergent, uniformly for $t>>X$ (and $0<x\le X$), we easily obtain that $\lim_{t\to+\infty} I_1(x,t)=0$. Similar computations can be carried out for all the integrals which are parts of $\mathcal{G}_1(x,t)$, and this implies (8.2) when $E=0$.
Integrating by parts once more – in the last four integrals of (8.4) – shows that

$$I_1(x,t)=\frac{1}{t^2}\times\{a\ sum\ of\ \ four\ integrals\}, \tag{8.5}$$

where the integrals in (8.5) are again absolutely convergent, uniformly for $t>>X$ (and $0<x\le X$). This implies that $\lim_{t\to+\infty}[tI_1(x,t)]=0$. Thus, working in this way, we obtain (8.2) when $E=1$.
Continuing in this way, by further integration by parts, we obtain (8.2) and (8.3) for every $E$.

**Lemma 2** *With notation as in Lemma 1, for fixed* $X>0$*,*

$$\lim_{t\to+\infty}\left[t^E\frac{\partial^{n+l}\mathcal{G}_1(x,t)}{\partial x^n\partial t^l}\right]=0\ \ and\ \ \lim_{t\to+\infty}\left[t^E\frac{\partial^{n+l}\mathcal{G}_2(x,t)}{\partial x^n\partial t^l}\right]=0, \tag{8.6}$$

*for nonnegative integers* $n$*,* $l$ *and* $E$*, uniformly for* $0<x\le X$*.*

**Proof** Let us consider the integral

$$I_2(x,t)=\int_{-\infty}^{\infty} e^{i\lambda x-\omega_2(\lambda)t}\{[i\omega_1(\lambda)+m]\hat{g}_1(\lambda)\}\frac{d\lambda}{2\rho(\lambda)},$$

which is a part of $\mathcal{G}_1(x,t)$, and let us write it as follows:

$$I_2(x,t)=\int_{\Gamma^{-,-}} e^{i\lambda x-\omega_2(\lambda)t}\{[i\omega_1(\lambda)+m]\hat{g}_1(\lambda)\}\frac{d\lambda}{2\rho(\lambda)}$$

$$+\int_{-1}^{1} e^{i\lambda x-\omega_2(\lambda)t}\{[i\omega_1(\lambda)+m]\hat{g}_1(\lambda)\}\frac{d\lambda}{2\rho(\lambda)}$$
$$+\sum_{k=1}^{N} g_1^{(k)}(0)\int_{\Gamma^{+,+}} e^{i\lambda x-\omega_2(\lambda)t}\{[i\omega_1(\lambda)+m]/(i\lambda)^k\}\frac{d\lambda}{2\rho(\lambda)}$$
$$+\int_{1}^{\infty} e^{i\lambda x-\omega_2(\lambda)t}\{[i\omega_1(\lambda)+m](g_1^{(N)}\hat{)}(\lambda)/(i\lambda)^N\}\frac{d\lambda}{2\rho(\lambda)}.$$

It follows that, for sufficiently large $N$,

$$\text{(8.7)}\quad \frac{\partial^{n+l}I_2(x,t)}{\partial x^n\partial t^l}=\int_{\Gamma^{-,-}}(i\lambda)^n[-\omega_2(\lambda)]^l e^{i\lambda x-\omega_2(\lambda)t}\{[i\omega_1(\lambda)+m]\hat{g}_1(\lambda)\}\frac{d\lambda}{2\rho(\lambda)}$$
$$+\int_{-1}^{1}(i\lambda)^n[-\omega_2(\lambda)]^l e^{i\lambda x-\omega_2(\lambda)t}\{[i\omega_1(\lambda)+m]\hat{g}_1(\lambda)\}\frac{d\lambda}{2\rho(\lambda)}$$
$$+\sum_{k=1}^{N} g_1^{(k)}(0)\int_{\Gamma^{+,+}}(i\lambda)^n[-\omega_2(\lambda)]^l e^{i\lambda x-\omega_2(\lambda)t}\{[i\omega_1(\lambda)+m]/(i\lambda)^k\}\frac{d\lambda}{2\rho(\lambda)}$$
$$+\int_{1}^{\infty}(i\lambda)^n[-\omega_2(\lambda)]^l e^{i\lambda x-\omega_2(\lambda)t}\{[i\omega_1(\lambda)+m](g_1^{(N)}\hat{)}(\lambda)/(i\lambda)^N\}\frac{d\lambda}{2\rho(\lambda)}.$$

Next, using the identity

$$\text{(8.8)}\qquad e^{i\lambda x-\omega_2(\lambda)t}=-\frac{1}{t}e^{i\lambda x}\frac{1}{\omega_2'(\lambda)}\frac{d[e^{-\omega_2(\lambda)t}]}{d\lambda}$$

and substituting the quantity $e^{i\lambda x-\omega_2(\lambda)t}$ in the integrals in the RHS of (8.7), we integrate by parts, as in the proof of Lemma 1. The conclusion is that

$$\text{(8.9)}\qquad \frac{\partial^{n+l}I_2(x,t)}{\partial x^n\partial t^l}=\frac{1}{t}\times\{\textit{a sum of four integrals}\},$$

where the integrals in (8.8) are absolutely convergent, uniformly for $t>>X$ and $0<x\le X$ (for a fixed $X>0$). This implies that

$$\lim_{t\to+\infty}\frac{\partial^{n+l}I_2(x,t)}{\partial x^n\partial t^l}=0,\ \text{uniformly for } 0<x\le X.$$

Working in this way with all the parts of $\mathcal{G}_1(x,t)$ and $\mathcal{G}_2(x,t)$, we prove (8.6) in the case $E=0$. The general case, $E\ge 1$, requires further integration by parts and can be easily treated, using (8.8).

**Lemma 3** *Assuming that the function $h_1(t)$ satisfies (8.1), let*

$$\mathcal{H}_1(x,t):=\mathcal{H}_{1,1}(x,t)-\mathcal{H}_{1,2}(x,t)\ \textit{and}\ \mathcal{H}_2(x,t):=-\mathcal{H}_{2,1}(x,t)+\mathcal{H}_{2,2}(x,t)$$

*be the parts of $\Psi_1(x,t)$ and $\Psi_2(x,t)$, repectively, which involve $h_1(t)$. Then, for fixed $X>0$,*

$$\text{(8.10)}\quad \lim_{t\to+\infty}\left\{t^E\frac{\partial^{n+l}[\mathcal{H}_1(x,t+T)-\mathcal{H}_1(x,t)]}{\partial x^n\partial t^l}\right\}=0\ \ \textit{and}\ \ \lim_{t\to+\infty}\left\{t^E\frac{\partial^{n+l}[\mathcal{H}_2(x,t+T)-\mathcal{H}_2(x,t)]}{\partial x^n\partial t^l}\right\}=0,$$

*for nonnegative integers $n$, $l$ and $E$, uniformly for $0<x\le X$.*

**Proof** Firstly, it follows from (8.1) that

(8.11) $$e^{-\omega(\lambda)(t+T)}\int_{\tau=0}^{t+T} h_1(\tau)e^{\omega(\lambda)\tau}d\tau - e^{-\omega(\lambda)t}\int_{\tau=0}^{t} h_1(\tau)e^{\omega(\lambda)\tau}d\tau$$

$$= e^{-\omega(\lambda)t}\int_0^t e^{\omega(\lambda)\tau}[h_1(\tau+T)-h_1(\tau)]d\tau + e^{-\omega(\lambda)t}\int_{-T}^{0} e^{\omega(\lambda)\tau}h_1(\tau+T)d\tau = e^{-\omega(\lambda)t}\int_{-T}^{0} e^{\omega(\lambda)\tau}h_1(\tau+T)d\tau\,,$$

for $\omega \in \{\omega_1, \omega_2\}$.

It follows that

(8.12) $$\mathcal{H}_{1,1}(x,t+T)-\mathcal{H}_{1,1}(x,t)$$

$$= \int_{-\infty}^{\infty} e^{i\lambda x-\omega_1(\lambda)(t+T)}\left[\int_{\tau=0}^{t+T} h_1(\tau)e^{\omega_1(\lambda)\tau}d\tau\right]\frac{\lambda d\lambda}{\rho(\lambda)} - \int_{-\infty}^{\infty} e^{i\lambda x-\omega_1(\lambda)t}\left[\int_{\tau=0}^{t} h_1(\tau)e^{\omega_1(\lambda)\tau}d\tau\right]\frac{\lambda d\lambda}{\rho(\lambda)}$$

$$= \int_{-\infty}^{\infty} e^{i\lambda x-\omega_1(\lambda)t}\left[\int_{-T}^{0} e^{\omega_1(\lambda)\tau}h_1(\tau+T)d\tau\right]\frac{\lambda d\lambda}{\rho(\lambda)}.$$

Working with $t >> X$ and $0 < x \le X$, and deforming the contour in the last integral of (8.12), we obtain

$$\mathcal{H}_{1,1}(x,t+T)-\mathcal{H}_{1,1}(x,t) = \int_{L^{-,+}+L^{+,-}} e^{i\lambda x-\omega_1(\lambda)t}\left[\int_{-T}^{0} e^{\omega_1(\lambda)\tau}h_1(\tau+T)d\tau\right]\frac{\lambda d\lambda}{\rho(\lambda)}.$$

Therefore,

(8.13) $$\frac{\partial^{n+l}[\mathcal{H}_{1,1}(x,t+T)-\mathcal{H}_{1,1}(x,t)]}{\partial x^n \partial t^l}$$

$$= \int_{L^{-,+}+L^{+,-}} (i\lambda)^n[-\omega_1(\lambda)]^l e^{i\lambda x-\omega_1(\lambda)t}\left[\int_{-T}^{0} e^{\omega_1(\lambda)\tau}h_1(\tau+T)d\tau\right]\frac{\lambda d\lambda}{\rho(\lambda)}$$

$$= -\frac{1}{t}\int_{L^{-,+}+L^{+,-}} (i\lambda)^n[-\omega_1(\lambda)]^l e^{i\lambda x}\frac{1}{\omega_1'(\lambda)}\frac{d[e^{-\omega_1(\lambda)t}]}{d\lambda}\left[\int_{-T}^{0} e^{\omega_1(\lambda)\tau}h_1(\tau+T)d\tau\right]\frac{\lambda d\lambda}{\rho(\lambda)}$$

$$= \frac{1}{t}\int_{L^{-,+}+L^{+,-}} e^{i\lambda x-\omega_1(\lambda)t}\frac{d}{d\lambda}\left\{(i\lambda)^n[-\omega_1(\lambda)]^l\frac{1}{\omega_1'(\lambda)}\left[\int_{-T}^{0} e^{\omega_1(\lambda)\tau}h_1(\tau+T)d\tau\right]\frac{\lambda}{\rho(\lambda)}\right\}d\lambda\,,$$

where the last equation is obtained by integration by parts.

Since the last integral in (8.13) is absolutely convergent, uniformly for $t >> X$ and $0 < x \le X$, it follows that

$$\lim_{t\to+\infty}\frac{\partial^{n+l}[\mathcal{H}_{1,1}(x,t+T)-\mathcal{H}_{1,1}(x,t)]}{\partial x^n \partial t^l} = 0.$$

Equations, analogous to (8.13), can also be obtained for the differences

(8.14) $$\mathcal{H}_{1,2}(x,t+T)-\mathcal{H}_{1,2}(x,t)\,,\ \mathcal{H}_{2,1}(x,t+T)-\mathcal{H}_{2,1}(x,t)\,,\ \mathcal{H}_{2,2}(x,t+T)-\mathcal{H}_{2,2}(x,t)\,,$$

and this implies (8.10),in the case $E = 0$.

The general case, $E \ge 1$, requires further integration by parts in the last integral of (8.13), and the analogous formulas for the differences (8.14).

**Lemma 4** *Assuming that the functions* $f_1(x,t)$ *and* $f_2(x,t)$ *satisfy (8.1), let*

$$\mathcal{F}_1(x,t) := \frac{1}{2\pi}[\mathcal{F}_{1,1}(x,t) + \mathcal{F}_{1,2}(x,t)] \text{ and } \mathcal{F}_2(x,t) := \frac{1}{2\pi}[\mathcal{F}_{2,1}(x,t) + \mathcal{F}_{2,2}(x,t)]$$

*be the parts of* $\Psi_1(x,t)$ *and* $\Psi_2(x,t)$, *repectively, which involve* $f_1(x,t)$ *and* $f_2(x,t)$.

*Then, for fixed* $X > 0$,

$$(8.15)\quad \lim_{t\to+\infty}\left\{t^E \frac{\partial^{n+l}[\mathcal{F}_1(x,t+T) - \mathcal{F}_1(x,t)]}{\partial x^n \partial t^l}\right\} = 0 \text{ and } \lim_{t\to+\infty}\left\{t^E \frac{\partial^{n+l}[\mathcal{F}_2(x,t+T) - \mathcal{F}_2(x,t)]}{\partial x^n \partial t^l}\right\} = 0,$$

*for nonnegative integers* $n$, $l$ *and* $E$, *uniformly for* $0 < x \le X$.

**Proof** It follows from (8.1) that

$$(8.16)\quad e^{-\omega(\lambda)(t+T)}\tilde{\hat{f}}(\lambda,\omega(\lambda),t+T) - e^{-\omega(\lambda)t}\tilde{\hat{f}}(\lambda,\omega(\lambda),t)$$

$$= e^{-\omega(\lambda)(t+T)}\int_0^{t+T} e^{-\omega(\lambda)\tau}\hat{f}(\lambda,\omega(\lambda),\tau)d\tau - e^{-\omega(\lambda)t}\int_0^{t} e^{-\omega(\lambda)\tau}\hat{f}(\lambda,\omega(\lambda),\tau)d\tau$$

$$= e^{-\omega(\lambda)t}\int_0^{t} e^{\omega(\lambda)\tau}[\hat{f}(\lambda,\tau+T) - \hat{f}(\lambda,\tau)]d\tau + e^{-\omega(\lambda)t}\int_{-T}^{0} e^{\omega(\lambda)\tau}\hat{f}(\lambda,\tau+T)d\tau$$

$$= e^{-\omega(\lambda)t}\int_{-T}^{0} e^{\omega(\lambda)\tau}\hat{f}(\lambda,\tau+T)d\tau,$$

for $f \in \{f_1, f_2\}$ and $\omega \in \{\omega_1, \omega_2\}$.

Using (8.16), we can easily prove (8.15) by combining computations which we used in the proofs of both Lemmas 2 and 3. The point here is that we have to invoke also the formula (7.7).

## 9. Further results

### (1) Behavior of the solution functions as $\overline{Q} \ni (x,t) \to (0,0)$

Let us recall that, as it follows from Theorem 1 (6th part), Theorem 4 (1st part) and Theorem 5 (3rd part),

$$h_1(0) = g_1(0) \;\Rightarrow\; \Psi_1(x,t) \text{ and } \Psi_2(x,t) \text{ are continuous on } \overline{Q} - \{(0,0)\}.$$

(The notation is as in section 7.)

***Proposition* 1** *Assuming in addition to (1.3) and (7.3), that* $h_1(0) = g_1(0)$, *we have that the limits*

$$\lim_{\substack{(x,t)\to(0,0)\\ (x,y)\in\overline{Q}-\{(0,0)\}}} \Psi_1(x,t) \quad \text{and} \quad \lim_{\substack{(x,t)\to(0,0)\\ (x,y)\in\overline{Q}-\{(0,0)\}}} \Psi_2(x,t)$$

*exist, extendind* $\Psi_1(x,t)$ *and* $\Psi_2(x,t)$ *to continuous functions on* $\overline{Q}$.

***Proof*** It follows immediately from the 6th conclusion of Theorem 1, the 1st conclusion of Theorem 4 and the computations leading to their proof.

***Proposition* 2** *Assuming (1.3), (7.3) and (7.4), we have that the limits*

$$\lim_{\substack{(x,t)\to(0,0)\\(x,y)\in\overline{Q}-\{(0,0)\}}} \frac{\partial\Psi_j(x,t)}{\partial x} \quad and \quad \lim_{\substack{(x,t)\to(0,0)\\(x,y)\in\overline{Q}-\{(0,0)\}}} \frac{\partial\Psi_j(x,t)}{\partial t}, \quad j=1,2,$$

*exist, extendind the derivatives* $\dfrac{\partial\Psi_j(x,t)}{\partial x}$ *and* $\dfrac{\partial\Psi_j(x,t)}{\partial t}$ *to continuous functions on* $\overline{Q}$.

***Proof*** It follows from the 4th conclusion of Theorem 4 and the computations leading to the proof of (7.5) and (7.6), carried out in Step 4 of the proof of Theorem 4.

**(2) A uniqueness theorem**

***Proposition* 3** *Assuming (1.3), (7.3) and (7.4), the solution* $(\Psi_1,\Psi_2)$ *of (7.1) & (7.2), defined in Theorem 4, is unique in the following sense:*

*If* $\mathcal{Y}_1(x,t)$ *&* $\mathcal{Y}_2(x,t)$ *are* $C^1$ *functions for* $(x,t)\in\overline{Q}-\{(0,0)\}$, *satisfy (7.1) & (7.2) (with* $(\Psi_1,\Psi_2)$ *replaced by* $(\mathcal{Y}_1,\mathcal{Y}_2)$*), and*

(9.1) $$\lim_{x\to+\infty}\mathcal{Y}_1(x,t)=0 \ \& \ \lim_{x\geq 1}\left|\mathcal{Y}_2(x,t)\right|<+\infty \ (\forall t),$$

*and, in addition, for every* $T>0$, *the functions* $\left|\mathcal{Y}_1(x,t)\right|^2$, $\left|\mathcal{Y}_2(x,t)\right|^2$, $\left|\partial\mathcal{Y}_1(x,t)/\partial t\right|^2$, $\left|\partial\mathcal{Y}_2(x,t)/\partial t\right|^2$, *are, uniformly for* $0<t\leq T$, *integrable with respect to* $x\in[0,+\infty)$, *i.e., there exists a positive function* $B_T(x)$ *such that* $\int_0^\infty B_T(x)dx<+\infty$ *and, for* $0<t\leq T$,

(9.2) $$\left|\mathcal{Y}_1(x,t)\right|^2\leq B_T(x),\ \left|\mathcal{Y}_2(x,t)\right|^2\leq B_T(x),\ \left|\partial\mathcal{Y}_1(x,t)/\partial t\right|^2\leq B_T(x),\ \left|\partial\mathcal{Y}_2(x,t)/\partial t\right|^2\leq B_T(x)\ (\forall x>0),$$

*then* $(\Psi_1(x,t),\Psi_2(x,t))\equiv(\mathcal{Y}_1(x,t),\mathcal{Y}_2(x,t))$.

***Proof*** Firstly, we point out that, in view of properties of the solution $(\Psi_1,\Psi_2)$, stated in Theorems 4, 5, and 6, and Propositions 1 and 2, (9.1) and (9.2) hold with $(\mathcal{Y}_1,\mathcal{Y}_2)$ replaced by $(\Psi_1,\Psi_2)$, by making $B_T(x)$ larger, if necessary.

Next, let us consider the differences

$$u_1(x,t):=\Psi_1(x,t)-\mathcal{Y}_1(x,t) \ \text{ and } \ u_2(x,t):=\Psi_2(x,t)-\mathcal{Y}_2(x,t),\ (x,t)\in\overline{Q}-\{(0,0)\}.$$

Then, for $(x,t)\in\overline{Q}-\{(0,0)\}$,

(9.3) $$\frac{\partial u_1(x,t)}{\partial t}=-\frac{\partial u_2(x,t)}{\partial x}-imu_1(x,t)$$

and

(9.4) $$\frac{\partial u_2(x,t)}{\partial t}=-\frac{\partial u_1(x,t)}{\partial x}+imu_2(x,t).$$

Also,

(9.5) $$u_1(x,0)=0 \ \text{ and } \ u_2(x,0)=0,\ \text{for } x>0,$$

and

(9.6) $$u_1(0,t)=0,\ \text{for } t>0.$$

Then, (9.3) gives

$$\overline{u_1}(x,t)\frac{\partial u_1(x,t)}{\partial t} = -\overline{u_1}(x,t)\frac{\partial u_2(x,t)}{\partial x} - im|u_1(x,t)|^2$$

and

$$u_1(x,t)\frac{\partial \overline{u_1}(x,t)}{\partial t} = -u_1(x,t)\frac{\partial \overline{u_2}(x,t)}{\partial x} + im|u_1(x,t)|^2.$$

Adding the above equations, we obtain

$$\frac{\partial[|u_1(x,t)|^2]}{\partial t} = -\overline{u_1}(x,t)\frac{\partial u_2(x,t)}{\partial x} - u_1(x,t)\frac{\partial \overline{u_2}(x,t)}{\partial x}. \tag{9.7}$$

Similarly, (9.4) gives

$$\frac{\partial[|u_2(x,t)|^2]}{\partial t} = -\overline{u_2}(x,t)\frac{\partial u_1(x,t)}{\partial x} - u_2(x,t)\frac{\partial \overline{u_1}(x,t)}{\partial x}. \tag{9.8}$$

Addind (9.7) and (9.8), we obtain

$$\frac{\partial[|u_1(x,t)|^2 + |u_2(x,t)|^2]}{\partial t} = -\frac{\partial[\overline{u_1}(x,t)u_2(x,t) + u_1(x,t)\overline{u_2}(x,t)]}{\partial x}.$$

It follows that, for $\mathrm{A} > 0$,

$$\int_{x=0}^{\mathrm{A}} \frac{\partial[|u_1(x,t)|^2 + |u_2(x,t)|^2]}{\partial t}dx = -\int_{x=0}^{\mathrm{A}} \frac{\partial[\overline{u_1}(x,t)u_2(x,t) + u_1(x,t)\overline{u_2}(x,t)]}{\partial x}dx$$

$$= -[\overline{u_1}(\mathrm{A},t)u_2(\mathrm{A},t) + u_1(\mathrm{A},t)\overline{u_2}(\mathrm{A},t)] + [\overline{u_1}(0,t)u_2(0,t) + u_1(0,t)\overline{u_2}(0,t)]$$

$$= -[\overline{u_1}(\mathrm{A},t)u_2(\mathrm{A},t) + u_1(\mathrm{A},t)\overline{u_2}(\mathrm{A},t)], \text{ for } t > 0,$$

where we used (9.6).

Therefore, letting $\mathrm{A} \to \infty$, we have, also in view of (9.1) (which holds with $(\mathcal{Y}_1, \mathcal{Y}_2)$ replaced by $(u_1, u_2)$),

$$\int_{x=0}^{\infty} \frac{\partial[|u_1(x,t)|^2 + |u_2(x,t)|^2]}{\partial t}dx = 0, \text{ for } t > 0. \tag{9.9}$$

In view of (9.2) (which holds with $(\mathcal{Y}_1, \mathcal{Y}_2)$ replaced by $(u_1, u_2)$, by making $B_T(x)$ larger, if necessary), we can switch the order of integration and differentiation in (9.9), obtaining

$$\frac{\partial}{\partial t}\left\{\int_{x=0}^{\infty}[|u_1(x,t)|^2 + |u_2(x,t)|^2]dx\right\} = 0, \text{ for } t > 0.$$

It follows that

$$\int_{x=0}^{\infty}[|u_1(x,t)|^2 + |u_2(x,t)|^2]dx = \int_{x=0}^{\infty}[|u_1(x,s)|^2 + |u_2(x,s)|^2]dx, \text{ for } 0 < s < t.$$

Letting $s \to 0^+$, we obtain, also in view of (9.5), that

$$\int_{x=0}^{\infty}[|u_1(x,t)|^2 + |u_2(x,t)|^2]dx = 0,$$

i.e., $u_1(x,t) \equiv 0$ and $u_2(x,t) \equiv 0$. This completes the proof.

**(3) Smoothness of the solution functions on the diagonal** $\{(x,t)\in Q : x=t\}$

As we have seen in Theorem 4 and Proposition 1, if $h_1(0)=g_1(0)$ then the functions $\Psi_1(x,t)$ and $\Psi_2(x,t)$ are continuous on $\overline{Q}$. Also the assumption (7.4) implies that $\Psi_1(x,t)$ and $\Psi_2(x,t)$ are $C^1$ in $\overline{Q}-\{0,0)\}$. Higher order of smoothness on the diagonal $\{x=t\}$ requires further combatibility conditions of the data at the origin. The following two propositions give examples in this direction. Their proofs are similar.

***Proposition* 4** *If*

$$h_1(0)=g_1(0),\ h_1'(0)=-g_2'(0)+f_1(0,0)\ \textit{and}\ \ h_1''(0)=g_1''(0)+\frac{\partial^2 f_1}{\partial t^2}(0,0)-\frac{\partial f_2}{\partial x}(0,0),$$

*then the functions* $\Psi_1(x,t)$ *and* $\Psi_2(x,t)$ *are* $C^2$ *in* $\overline{Q}-\{(0,0)\}$, *and the limits*

$$\lim_{\substack{(x,t)\to(0,0)\\ (x,y)\in\overline{Q}-\{(0,0)\}}} \frac{\partial^{n+l}\Psi_j(x,t)}{\partial x^n \partial t^l},\ \ n+l\le 2,$$

*exist.*

***Proposition* 5** *If*

$$h_1(0)=g_1(0),\ h_1'(0)=-g_2'(0)+f_1(0,0),\ h_1''(0)=g_1''(0)+\frac{\partial^2 f_1}{\partial t^2}(0,0)-\frac{\partial f_2}{\partial x}(0,0)$$

$$\textit{and}\ \ h_1'''(0)=-g_2'''(0)+\frac{\partial^2 f_1}{\partial t^2}(0,0)-\frac{\partial^2 f_2}{\partial x\partial t}(0,0)-im\frac{\partial f_1}{\partial t}(0,0),$$

*then the functions* $\Psi_1(x,t)$ *and* $\Psi_2(x,t)$ *are* $C^3$ *in* $\overline{Q}-\{(0,0)\}$, *and the limits*

$$\lim_{\substack{(x,t)\to(0,0)\\ (x,y)\in\overline{Q}-\{(0,0)\}}} \frac{\partial^{n+l}\Psi_j(x,t)}{\partial x^n \partial t^l},\ \ n+l\le 3,$$

*exist.*

## References

[1] P. A. M. Dirac, The quantum theory of the electron, Proc. Roy. Soc. London Ser. A 117 (1928), 610–624.
[2] B. Thaller, The Dirac Equation, Springer, 1992.
[3] A. Das, Lectures on Quantum Field Theory, World Scientific, 2008.
[4] M. E. Peskin and D. V. Schroeder, An Introduction to Quantum Field Theory, Addison-Wesley, 1995.
[5] A.S. Fokas, A unified transform method for solving linear and certain nonlinear PDEs, *Proc. Roy. Soc. London Ser. A* (1997).
[6] A.S. Fokas, On the integrability of certain linear and nonlinear partial differential equations, *J. Math. Phys.* (2000).
[7] Fokas, A. S. (2001). Two–dimensional linear partial differential equations in a convex polygon. *Proc. Roy. Soc. London Ser. A*, *457*, 371-393.
[8] A.S. Fokas, A new transform method for evolution partial differential equations, *IMA J. Appl. Math.* (2002).

[9] A.S. Fokas, “From Green to Lax via Fourier”, In: Recent Advances in Nonlinear Partial Differential Equations and Applications, Proc. Sympos. Appl. Math. 65 (Dedicated to P.D. Lax and L. Nirenberg), Amer. Math. Soc. (2007).
[10] A.S. Fokas, “Lax Pairs: A novel type of separability”, *Inv. Probl.* (2009).
[11] A.S. Fokas, A Unified Approach to Boundary Value Problems, CBMS-NSF Series Appl. Math. 78, SIAM (2008).
[12] A. S. Fokas, E. A. Spence, Synthesis, as opposed to separation, of variables, *SIAM Rev.* (2012).
[13] B. Deconinck, T. Trogdon, V. Vasan, The method of Fokas for solving linear partial differential equations, *SIAM Rev.* (2014).
[14] Fokas, A. S., & Pelloni, B. (Eds.). (2014). *Unified transform for boundary value problems: applications and advances*. Society for Industrial and Applied Mathematics.
[15] A Chatziafratis, AS Fokas, K Kalimeris, The Fokas method for evolution partial differential equations, *Partial Differ. Equ. Appl. Math.* (2025).
[16] C. Garner, W. Green, Solving the Dirac equation with the unified transform method, SIAM Online (2019).
[17] W. E. Thirring, A soluble relativistic field theory, Ann. Physics 3 (1958), 91–112.
[22] S. Machihara, K. Nakanishi and T. Ozawa, Small global solutions and the nonrelativistic limit for the nonlinear Dirac equation, *Rev. Mat. Iberoamericana* 19 (2003), 179–194.
[18] I. P. Naumkin, Cubic nonlinear Dirac equation in a quarter plane, J. Math. Anal. Appl. 434 (2016), 1633–1664.
[19] I. P. Naumkin, Initial-boundary value problem for the one dimensional Thirring model, J. Differential Equations 261 (2016), 4486–4523.
[20] B. Xia, The massive Thirring system in the quarter plane, Proc. Roy. Soc. Edinburgh Sect. A 150 (2020), 2387–2416.
[21] F. Liu, K. Tian, Y. Zhang, Z. Gong and S. Lyv, Low regularity for nonlinear Dirac equation in the half line, in Proceedings of the International Academic Summer Conference on Number Theory and Information Security (NTIS 2023), Advances in Physics Research 9, Atlantis Press, 2024, 88–103.
[23] J.L. Bona, A. Chatziafratis, H. Chen, S. Kamvissis, The linearised BBM equation on the half-line revisited, *Lett. Math. Phys.* (2024).
[24] M. M. Cavalcanti, A. Chatziafratis, and L. Rosier. A unified approach to lack of controllability for linear evolution equations on quarter-planes. *preprint* (2025).
[25] A. Chatziafratis, E.C. Aifantis, A. Carbery, A.S. Fokas, Integral representations for the double-diffusion system on the half-line, *Z. Angew. Math. Phys.* (2024).
[26] A. Chatziafratis, A.S. Fokas, E.C. Aifantis, On Barenblatt’s pseudoparabolic equation on the half-line via the Fokas method, *Z. Angew. Math. Mech.* (2024).
[27] A. Chatziafratis, C. Giorgi, A. Miranville, and F. Zullo, On generalised d’Alembert-type integral representations for the damped wave equation on the quarter-plane, *preprint* (2025).
[28] A. Chatziafratis, L. Grafakos, S. Kamvissis, Long-range instabilities for linear evolution PDE on semi-bounded domains, *Dyn. PDE* (2024).
[29] Chatziafratis, A., Kamvissis, S., Infinity of solutions to initial-boundary value problems for linear constant-coefficient evolution PDEs on semi-infinite intervals. *Bull. London Math. Soc.*, *57*(12), 4111-4121 (2025).
[30] A. Chatziafratis, S. Kamvissis, A. Miranville. Rigorous qualitative theory for evolution equations with variable coefficients on the quarter-plane. *preprint* (2026).
[31] A. Chatziafratis, S. Kamvissis, I. G. Stratis, Boundary behavior of the solution to the linear KdV equation on the half-line, *Stud. Appl. Math*. (2023).
[32] A. Chatziafratis, D. Mantzavinos, Boundary behavior for the heat equation on the half-line, *Math. Meth. Appl. Sci.* (2022).
[33] A. Chatziafratis, A. Miranville, G. Karali, A.S. Fokas, E.C. Aifantis, Higher-order diffusion and Cahn–Hilliard-type models revisited on the half-line, *Math. Models Meth. Appl. Sci.* (2025).
[34] A. Chatziafratis and T. Ozawa. New instability, blow-up and break-down effects for Sobolev-type evolution PDE: asymptotic analysis for a celebrated pseudo-parabolic model on the quarter-plane. *Partial Differ. Equ. Appl.*, 5(30), 2024.

[35] A. Chatziafratis, T. Ozawa, S.-F. Tian, Rigorous analysis of the unified transform method and long-range instabilities for the inhomogeneous time-dependent Schrödinger equation on the quarter-plane, *Math. Ann.* (2024).
[36] A. Chatziafratis, On the Fokas Method for Evolution PDE on the Half-Space, Thesis (in Greek), Advisors: N. Alikakos, G. Barbatis, I.G. Stratis, National and Kapodistrian University of Athens (2019).
[37] A. Chatziafratis, G. Karali, C. E. Synolakis, Rigorous analytical solution for the linearized Whitham-Broer-Kaup system on the half-line, *Stud Appl Math* (2026).
[38] A. Chatziafratis, E. Luca and S. L. Waters, On the Milne–Barenblatt pseudo-parabolic equation with time-dependent coefficients in biological fluid mechanics and regenerative medicine, *Discrete Contin. Dyn. Syst. S* (2026).
[39] A. Chatziafratis, A. Miranville, T. Ozawa, The linear Cahn-Hilliard equation with an interface, arxiv *preprint* (2026).
[40] A. Chatziafratis and T. Ozawa, Non-standard initial-boundary-value and interface problems for the linear BBM equation. *Discrete Contin. Dyn. Syst. S* (2026).